\documentclass[lettersize,journal]{IEEEtran}
\usepackage{amsmath,amsfonts}
\usepackage{algorithmic}
\usepackage{algorithm}
\usepackage{array}
\usepackage[caption=false,font=small,labelfont=sf,textfont=sf]{subfig}
\usepackage{textcomp}
\usepackage{stfloats}
\usepackage{url}
\usepackage{verbatim}
\usepackage{graphicx}
\usepackage{cite}
\usepackage{color,array}
\usepackage{graphicx}

\usepackage{caption}
\usepackage{float}
\usepackage{booktabs}
\usepackage{multirow}
\usepackage{makecell}

\usepackage{booktabs}
\usepackage{diagbox}

\usepackage{cases}

\newtheorem{theorem}{Theorem}

\newtheorem{lemma}{Lemma}

\newtheorem{remark}{Remark}

\graphicspath{{./images/}}
\DeclareGraphicsExtensions{.eps,.ps,.jpg,.bmp}

\makeatletter % Force all sequence numbers to be the same size
\renewcommand{\maketag@@@}[1]{\hbox{\m@th\normalsize\normalfont#1}}%
\makeatother

\begin{document}
	
	\title{Intelligent acceleration adaptive control of linear $2\times2$ hyperbolic PDE systems}
	
	\author{Xianhe~Zhang, Yu~Xiao, Xiaodong~Xu,\IEEEmembership{ Member,~IEEE}, Biao Luo,\IEEEmembership{ Senior Member,~IEEE} \vspace{0mm}
	
	\thanks{
			Xianhe Zhang, Yu Xiao, Xiaodong Xu, Biao Luo, Zhiyong Chen are with the School of Automation, Central South University, Changsha, China (e-mail:
			zxh\_csu@csu.edu.cn, yu\_xiao@csu.edu.cn, xx1@ualberta.ca, biao.luo@hotmail.com).}
%		\thanks{Zhiyong Chen is with the School of Electrical Engineering and Computing, University of Newcastle, Australia  (e-mail:
%			zhiyong.chen@newcastle.edu.au).}
	}


	
	\maketitle
	
	\begin{abstract}
	
	Traditional approaches to stabilizing hyperbolic PDEs, such as PDE backstepping, often encounter challenges when dealing with high-dimensional or complex nonlinear problems. Their solutions require high computational and analytical costs. 
	Recently, neural operators (NOs) for the backstepping design of first-order hyperbolic partial differential equations (PDEs) have been introduced, which rapidly generate gain kernel without requiring online numerical solution.
	In this paper we apply neural operators to a more complex class of 2 × 2 hyperbolic PDE systems for adaptive stability control. 
	Once the NO has been well-trained offline on a sufficient training set obtained using a numerical solver, the kernel equation no longer needs to be solved again, thereby avoiding the high computational cost during online operations.
	Specifically, we introduce the deep operator network (DeepONet), a neural network framework, to learn the nonlinear operator of the system parameters to the kernel gain. 
	The approximate backstepping kernel is obtained by utilizing the network after learning, instead of numerically solving the kernel equations in the form of PDEs, to further derive the approximate controller and the target system. 
	We analyze the existence and approximation of DeepONet operators and provide stability and convergence proofs for the closed-loop systems with NOs. 
	Finally, the effectiveness of the proposed NN-adaptive control scheme is verified by comparative simulation, which shows that the NN operator achieved up to three orders of magnitude faster compared to conventional PDE solvers, significantly improving real-time control performance.
	\end{abstract}

	\begin{IEEEkeywords}
		Adaptive control, Intelligent acceleration, Hyperbolic PDE, DeepONet, Neural operator.
	\end{IEEEkeywords}

\section{Introduction}
	\subsection{Background and related works}
	\IEEEPARstart{I}{n} the realm of computational science and engineering, the applications of neural networks have emerged as a transformative paradigm, revolutionizing how we approach complex problem-solving across various domains \cite{lecun2015deep}, increasing the speed and efficiency of problem solving \cite{bongard2007automated}. Within this domain, the study of PDEs, particularly hyperbolic PDE systems, holds significant relevance, as they underpin dynamic phenomena across diverse fields such as fluid dynamics \cite{holl2020learning, zhang2020data}, electromagnetism \cite{wang2021regulation, wang2021deep}, and elasticity \cite{paranjape2013pde, cuong2023robust}.
	
	Traditional methods for solving PDE systems face challenges when dealing with high-dimensional or nonlinear systems, often leading to computational inefficiency and limitations in capturing intricate dynamics accurately \cite{krstic1995nonlinear}. Significant advancements have recently been witnessed with the integration of neural networks into the computational framework, offering approaches that exhibit enhanced efficiency, scalability, and adaptability \cite{liu2021active, abdelaty2021daics, niu2023adaptive, sirignano2018dgm, li2021muscle, becattini2023transformer}. But many of these methods still demand a significant amount of training data, which can be expensive to acquire. Therefore, we introduce neural operators to improve the efficiency of data utilization.
	
	\subsubsection{Neural operator}
	The neural operator approach to solving parametric PDE equations requires the learning of infinite dimensional function space mappings (operators), in contrast to neural networks that aim to learn  vector space mappings (functions). NOs enable us to learn the mapping of the entire set of system parameters without the need for retraining as the parameters change. But some of these methods \cite{li2020neural, li2020multipole} require high costs for data acquisition and training. 
	To this end, recent research \cite{lu2021learning} proposes the DeepONet network, a popular framework for learning nonlinear operators. This network comprises two sub-networks: a trunk network and a branch network. It offers a straightforward and intuitive model framework that can be trained rapidly. For instance,  \cite{wang2021learning} introduces a physics-informed DeepONet for fast prediction of solutions to various parametric PDEs. And \cite{lu2022multifidelity} develops the Multifidelity DeepONet and applies it to learn the Phonon Boltzmann transport equation, achieving rapid solutions.
	According to the universal approximation theorem \cite{chen1995universal, chen1993approximations}, a neural network with a single hidden layer can accurately approximate any nonlinear operator. Therefore, DeepONet can theoretically obtain neural approximation operators with arbitrary accuracy. This provides inspiration to incorporate the DeepONet operator into the backstepping-based adaptive control of hyperbolic PDE systems.
	
	\subsubsection{Adaptive control of hyperbolic PDEs}
	Over the decades, adaptive control techniques have become crucial in solving complex control problems for hyperbolic PDE systems \cite{slotine1991applied}. These methods provide robust and effective control strategies by flexibly adapting to the uncertainties, disturbances, and parameter variations inherent in real-world systems \cite{li2023learning,cui2023event}. However, designing adaptive controllers for hyperbolic PDEs in practical applications is often
	challenging.
	Related research \cite{vazquez2008control} reviews various control strategies for linear hyperbolic systems, including adaptive control methods. An adaptive control scheme for stabilizing a specific class of hyperbolic PDE systems via boundary induction only is developed in \cite{anfinsen2017adaptive}. Additionally, \cite{coron2013local} stabilizes a pair of coupled hyperbolic PDEs with a single boundary input, and \cite{di2013stabilization} extends this approach to $n+1$ hyperbolic PDEs.
	In \cite{wang2021regulation}, an adaptive boundary control scheme based on deterministic equivalence is proposed. \cite{anfinsen2018adaptive} proposes two state feedback adaptive control laws for linear hyperbolic PDEs systems with uncertain domain and boundary parameters. 
	These studies achieved adaptive control via PDE backstepping, but due to the real-time variation of the system parameters, the backstepping kernel equations have to be solved repeatedly, imposing a heavy computational burden.
%	We consider adaptive stabilization control of a class of $2\times2$ hyperbolic PDE system by using the backstepping method, 
	Therefore we introducing neural operators to approximate the kernel gain to address this issue.
	
	\subsubsection{NO approximations for PDE control}
	Recently, neural operators combined with backstepping control have been employed to approximate the kernel gain. \cite{krstic2024neural} develops neural operators for backstepping design of first-order hyperbolic PDEs. In \cite{bhan2023operator}, the utilization of neural operators to approximate parameter identifiers and control gains for adaptive control is investigated. Additionally, \cite{bhan2023neural} employs deep neural networks to learn nonlinear operators from object parameters to control gains, thereby eliminating the need for gain function computations. Furthermore, \cite{lamarque2024gain} introduces a neural operator to accelerate the gain scheduling process of a go-up PDE system. 
	These articles use sufficiently trained neural network operators instead of the nonlinear operators of the system parameters to the backstepping kernel of the PDE, thereby avoiding the cumbersome computation of the numerical solution and speeding up the control process. As long as the NO can generate the approximate kernel function with sufficient accuracy, the stability of the closed-loop system can be guaranteed.
	
	\subsection{Paper structure and contributions}
	
	Previous studies \cite{krstic2024neural},\cite{bhan2023neural},\cite{vazquez2024gain} have employed NOs for approximating kernel gains of first-order PDE systems. Typically, these works involve estimating the backstepping kernel only once.
	% the backstepping kernel needs to be estimated only once. 
	In contrast, our research focuses on adaptive control of hyperbolic PDE systems, which contain unknown parameters requiring repeated solving of the backstepping kernel.
	Additionally, compared with \cite{bhan2023operator},\cite{lamarque2024adaptive}, which emphasize the methodology and analysis of adaptive control schemes with neural operators, we are committed to applying NOs to practical PDE systems to accelerate backstepping kernel solving and improve the real-time control performance of the system.
	Specifically, this paper combines neural operators with backstepping design, applied to adaptive stabilization control of a class of 2$\times$2 hyperbolic PDE systems.
	However, employing the DeepONet operator to approximate the numerical solution of the backstepping kernel equation introduces gain errors, which in turn add perturbation terms to the controlled system. These perturbations complicate the theoretical analysis of the stability of the PDE system.
	To this end, we provide a rigorous theoretical proof of the stability of closed-loop systems under approximate controllers. Simulation experiments verify the effectiveness of the backstepping controller with approximated NOs gain, exhibiting a three-order-of-magnitude acceleration in kernel gain solving.
	
%	Based on the above research objectives, the technical route of this paper is illustrated in Fig.\ref{accelerated control}. 
%	The left-hand side represents the traditional backstepping adaptive control scheme, where the original system is decomposed using the swapping design and backstepping methods, then the kernel equations, control laws, and the transformed target system are subsequently derived. 
%	The right-hand side represents the intelligent acceleration scheme using neural operators. In this approach, the offline well-trained DeepONet operators are employed to estimate the kernel gain, replacing the cumbersome analytical computation process and resulting in a significant speed-up in control.
	
	To summarize, the contributions of this paper are outlined below.
	\begin{itemize}
		\item   We extend the neural network approximation method to a more complex task of adaptive stabilization control for a class of 2$\times$2 hyperbolic PDE systems, enriching the theory of PDE adaptive control with neural operators.
		\item   This paper theoretically verifies the operator approximation property of backstepping kernel. The neural operator kernel can accelerate the control process by three orders of magnitude and significantly improves the real-time performance of the system online control.
		\item   A detailed analytical proof of the stability of the closed-loop system with an approximate kernel is provided based on Lyapunov theory. Simulation experiments are conducted to validate the effectiveness of the DeepONet operators, demonstrating the efficacy of the NOs in accelerating control.
	\end{itemize}	
	
	\textit{Paper structure}: Section II describes the control system and its objective. Section III illustrates the basic adaptive control based on swapping design. Section IV describes the operators network and the properties of approximate kernel. Section V analyzes the application of the approximate operators and theoretically proves the stability and convergence of the system. Section VI verifies the validity of the designed approximate controller and acceleration effect of NOs through simulation experiments. Finally, Section VII concludes the paper.

	\subsection{Notation}
	Throughout the paper, $\mathbb{R}^n$ denotes the  $n$ dimensional space, and  the corresponding Euclidean norms are denoted $\left|\cdot\right|$.
	For a signal $\psi (x,t)$ defined on $0 \le x \le 1,t \ge 0$, $\left\| \psi \right\|$ denotes the $L_2$-norm 
		$\left\| {\psi(t)} \right\| = \sqrt {\int_0^1 {{\psi^2}(x,t)dx} }$.
	For a time-varying signal $\psi\left(x,t\right)\in \mathbb{R}$, the $L_p$-norm is denoted by
    $\psi  \in {L_p} \Leftrightarrow {\left\| \psi  \right\|_p} = {(\int_0^\infty  {{{\left| {\psi (t)} \right|}^p}dt} )^{\frac{1}{p}}} < \infty $.
	For the vector signal $\psi(x)$ for all $x\in\left[0,1\right]$,
	we introduce the following operator 
	${I_\delta }[\psi] = \int_0^1 {{e^{\delta x}}\psi(x)dx} $
	with the derived norm 
	${I_\delta }[\psi^T{\psi}] = {\left\| \psi \right\|_\delta^2 } = \int_0^1 {{e^{\delta x}}{\psi^T}(x)\psi(x)dx} $.
	The following properties can be derived
	\begin{equation} \label{property 1}
		\setlength{\abovedisplayskip}{2pt}
		\setlength{\belowdisplayskip}{1pt}
	{I_\delta }[\psi{\psi_x}] = \frac{1}{2}({e^\delta }{\psi^2}(1) - {\psi^2}(0) - \delta {\left\| \psi \right\|_\delta^2 })
	\end{equation}
	\begin{equation} \label{property 2}
		\setlength{\abovedisplayskip}{1pt}
		\setlength{\belowdisplayskip}{2pt}
		\left\|  \cdot  \right\|_{ - \delta }^2 \le {\left\|  \cdot  \right\|^2} \le {e^\delta }\left\|  \cdot  \right\|_{ - \delta }^2 \le \left\|  \cdot  \right\|_\delta ^2,\delta  \ge 1
	\end{equation}
	
	Moreover, use the shorthand notation ${\psi _x} = {\partial _x}\psi  = \frac{{\partial \psi }}{{\partial x}}$ and ${\psi _t} = {\partial _t}\psi  = \frac{{\partial \psi }}{{\partial t}}$.
	
	\begin{figure}[!t] 
		\centering
		%\subfloat[parameter estimation]
		{\includegraphics[width=0.7\linewidth]{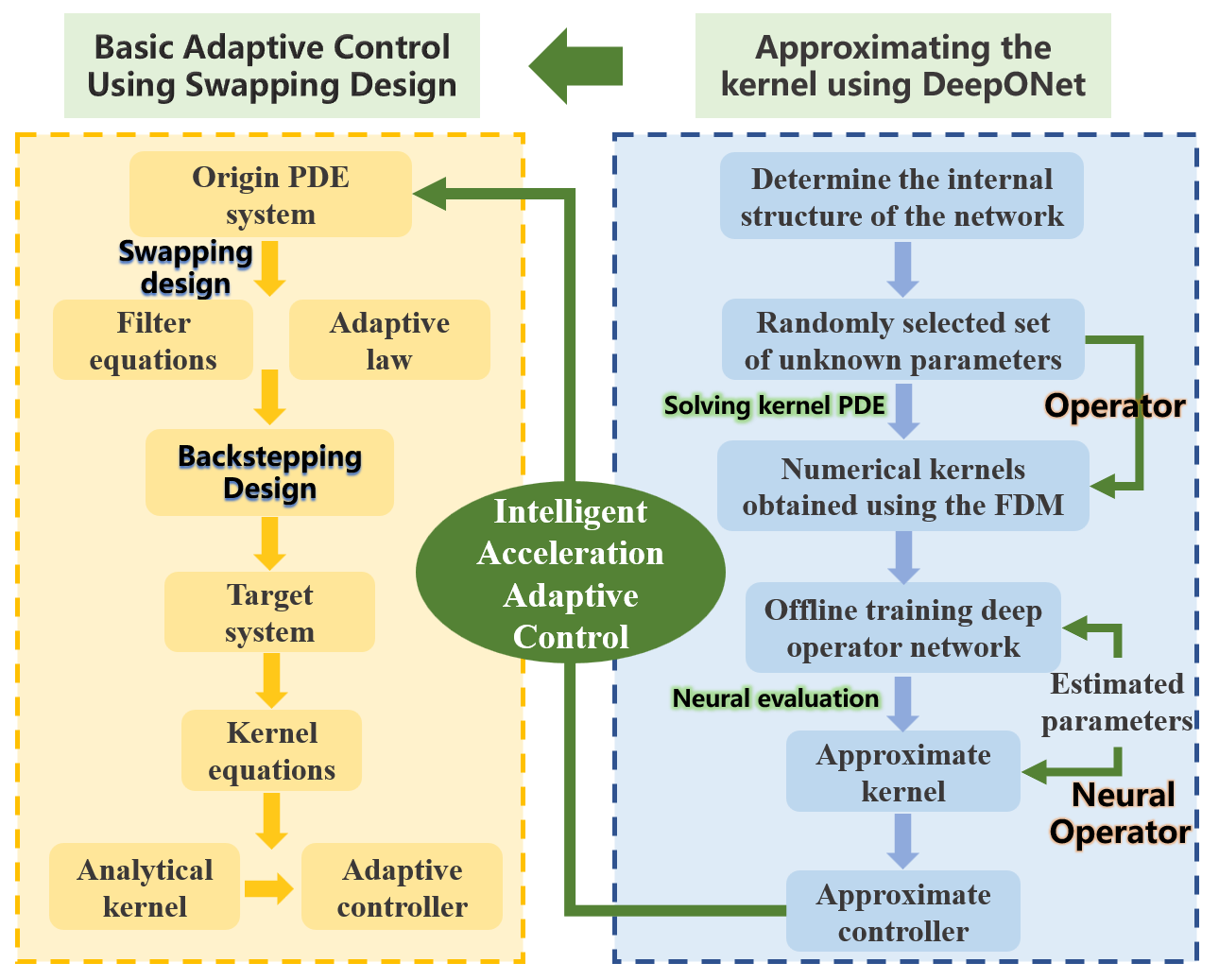}} 
		\caption{The intelligent acceleration adaptive control scheme combined with DeepONet. The approximated kernel obtained from the off-line trained operator network instead of the analytic kernel in numerical computation.}
		\label{accelerated control}
	\end{figure}

\section{Problem statement}
    Consider the following linear $2\times2$ hyperbolic PDE systems with constant in-domain coefficients
    \begin{subequations}\label{original_system}
    	\setlength{\abovedisplayskip}{2pt}
    	\setlength{\belowdisplayskip}{2pt}
    	\begin{align}
     		{\psi _t}(x,t) + \lambda {\psi_x}(x,t) &= {m_1}\psi(x,t) + {m_2}\varphi(x,t)\\
    		{\varphi_t}(x,t) - \mu {\varphi_x}(x,t) &= {m_3}\psi(x,t) + {m_4}\varphi(x,t)\\
    		\psi(0,t) &= q\varphi(0,t)\\
    		\varphi (1,t) &= U(t)
    	\end{align}
    \end{subequations}
    defined for $0 \le x \le 1,t \ge 0$, where $\psi(x,t)$,$\varphi(x,t)$ are the system states, $U$ is the control input, and $0 < \lambda  \in \mathbb{R}$, $0 < \mu  \in \mathbb{R}$ are known transport speeds. The boundary parameter $q$ is known constant,  while the coupling coefficient ${m_1}, {m_2}, {m_3}, {m_4}$ are unknown. The initial condition $\psi (x,0) = {\psi _0}(x)$, $\varphi (x,0) = {\varphi _0}(x)$ are assumed to satisfy ${\psi _0}(x),{\varphi _0}(x) \in {L_2}$. Since the system bounds are not strictly defined, we assume that there have some bounds on ${m_i}$ defined as ${{{\bar m}_i}}$ so that 
    \begin{equation} \label{parameter bounds}
    	\setlength{\abovedisplayskip}{1pt}
    	\setlength{\belowdisplayskip}{1pt}
    	\left| {{m_i}} \right| \le {{\bar m}_i},i = 1,2,3,4 .
    \end{equation}
    
    The control goal is to design an adaptive state controller that achieves regulation of the system states $\psi$ and $\varphi$ to zero equilibrium profile, while ensuring that all additional signals remain within certain bounds. 
    
    Furthermore, 
    this paper applies neural operators to the hyperbolic PDE system, significantly accelerates the control process. 
    When using the DeepONet approximation kernel instead of the backstepping kernel by numerical solutions, it is important to analyze the accuracy of the neural operator. It is also necessary to consider the boundedness and convergence of the controlled system with NOs.
 
\section{Basic adaptive control using swapping design} \label{basic adaptive control}
    In this section, by using the swapping design, the system states are transformed into a linear combination of boundary parameter filters, system parameters, and error terms. Therefore, the adaptive controller design constitutes the following three subsections.
    %The kernel gain is obtained by backstepping the design and further deriving the controller and target system
    
	\subsection{Filter equation}
	Firstly, we introduce the parameter filters that are defined on $0 \le x \le 1,t \ge 0$:
	\begin{subequations}\label{filter1}
		\setlength{\abovedisplayskip}{2pt}
		\setlength{\belowdisplayskip}{0pt}
		\begin{align}
		{\partial_t}{n_1}(x,t)+\lambda {\partial_x}{n_1}(x,t) &= \psi (x,t)&&{n_1}(0,t) = 0  \label{filter1 a} \\
		{\partial_t}{n_2}(x,t)+\lambda {\partial_x}{n_2}(x,t) &= \varphi (x,t)&&{n_2}(0,t) = 0 \label{filter1 b} \\
		{\partial_t}{n_3}(x,t)+\lambda {\partial_x}{n_3}(x,t) &= 0, &&{n_3}(0,t) = \varphi (0,t) \label{filter1 c}
		\end{align}
	\end{subequations}
	\begin{subequations}\label{filter2}
		\setlength{\abovedisplayskip}{0pt}
		\setlength{\belowdisplayskip}{2pt}
		\begin{align}
		{\partial _t}{z_1}(x,t)-\mu {\partial_x}{z_1}(x,t)&= \psi (x,t)&&{z_1}(1,t) = 0  \label{filter2 a} \\
		{\partial _t}{z_2}(x,t)-\mu {\partial_x}{z_2}(x,t)&= \varphi (x,t)&&{z_2}(1,t) = 0 \label{filter2 b} \\
		{\partial _t}{z_3}(x,t)-\mu {\partial_x}{z_3}(x,t)&= 0, &&{z_3}(1,t) = U(t) \label{filter2 c}
		\end{align}
	\end{subequations}
	where the initial states conditions ${n_i}(x,0) = n_i^0(x),{z_i}(x,0) = z_i^0(x),i = 1,2,3$ satisfy $n_1^0,n_2^0,n_3^0,z_1^0,z_2^0,z_3^0 \in {L_2}$.
	
    The following aggregated symbols are defined for  analysis easily
    \begin{subequations}\label{zn}
    	\setlength{\abovedisplayskip}{2pt}
    	\setlength{\belowdisplayskip}{1pt}
	    \begin{align}
	    n(x,t) = {[{n_1}(x,t),{n_2}(x,t)]^T}\\
	    z(x,t) = {[{z_1}(x,t),{z_2}(x,t)]^T},
		\end{align}
    \end{subequations}
    \begin{equation}\label{MN}
    	\setlength{\abovedisplayskip}{0pt}
   		\setlength{\belowdisplayskip}{2pt}
   		M = {[{m_1},{m_2}]^T},N = {[{m_3},{m_4}]^T}.
    \end{equation}
    Then we introduce the non-adaptive estimate of the system states
	\begin{subequations}\label{non-adaptive states}
		\setlength{\abovedisplayskip}{2pt}
		\setlength{\belowdisplayskip}{2pt}
		\begin{align}
			\bar \psi (x,t) &= {n^T}(x,t)M + {n_3}(x,t)q\\
			\bar \varphi (x,t) &= {z^T}(x,t)N + {z_3}(x,t).
		\end{align}
	\end{subequations}
	The non-adaptive system states estimation errors
	\begin{subequations}\label{non-adaptive errors}
		\setlength{\abovedisplayskip}{2pt}
		\setlength{\belowdisplayskip}{2pt}
		\begin{align}
			e(x,t) &= \psi (x,t) - \bar \psi (x,t)\\
			\tau (x,t) &= \varphi (x,t) - \bar \varphi (x,t).
		\end{align}
	\end{subequations}
	From \eqref{original_system}, we can straightforwardly obtain
	\begin{subequations}\label{non-adaptive errors system}
		\setlength{\abovedisplayskip}{2pt}
		\setlength{\belowdisplayskip}{2pt}
		\begin{align}
			{e_t}(x,t) + \lambda {e_x}(x,t) &= 0, \quad e(0,t) = 0, \label{real error bound 1}\\
			{\tau _t}(x,t) - \mu {\tau _x}(x,t) &= 0, \quad \tau (1,t) = 0
			\label{real error bound 2}.
		\end{align}
	\end{subequations}
	
	\subsection{Parameter estimation law}
	Then, we introduce the adaptive estimate
	\begin{subequations}\label{adaptive states}
		\setlength{\abovedisplayskip}{2pt}
		\setlength{\belowdisplayskip}{2pt}
		\begin{align}
			\hat \psi (x,t) &= {n^T}(x,t)\hat M + {n_3}(x,t)q\\
			\hat \varphi (x,t) &= {z^T}(x,t)\hat N + {z_3}(x,t),
		\end{align}
	\end{subequations}
	where $\hat M = {[{{\hat m}_1},{{\hat m}_2}]^T}$,$\hat N = {[{{\hat m}_3},{{\hat m}_4}]^T}$ denote the estimates of coupling coefficient. So that the prediction errors ${\hat e}$,${\hat \tau }$ denotes as
	\begin{subequations}\label{adaptive errors}
		\setlength{\abovedisplayskip}{1pt}
		\setlength{\belowdisplayskip}{2pt}
		\begin{align}
			&\hat e(x,t) = \psi (x,t) - \hat \psi (x,t)\\
			&\hat \tau (x,t) = \varphi (x,t) - \hat \varphi (x,t).
		\end{align}
	\end{subequations}
	Introducing the following adaptive laws
	\begin{subequations}\label{adaptive laws}
		\setlength{\abovedisplayskip}{2pt}
		\setlength{\belowdisplayskip}{2pt}
		\begin{small}
		\begin{align}
			&{{\dot {\hat m}}_1} \!=\! {\mathcal{P}_{{{\bar m}_1}}}\left\{ {{\rho _1}\frac{{\int_0^1 {{n_1}(x)\hat e(x)dx} }}{{1 + {{\left\| {{n_1}} \right\|}^2} + {{\left\| {{n_2}} \right\|}^2} + {{\left\| {{n_3}} \right\|}^2}}}} \right\}\\
			&{{\dot {\hat m}}_2} \!=\! {\mathcal{P}_{{{\bar m}_2}}}\left\{ {{\rho _2}\frac{{\int_0^1 {{n_2}(x)\hat e(x)dx} }}{{1 + {{\left\| {{n_1}} \right\|}^2} + {{\left\| {{n_2}} \right\|}^2} + {{\left\| {{n_3}} \right\|}^2}}}} \right\}\\
			&{{\dot {\hat m}}_3} \!=\! {\mathcal{P}_{{{\bar m}_3}}} \!\! \left\{ \! {{\rho _3} \!\! \left[ \! {\frac{{\int_0^1 {{z_1}(x)\hat \tau (x)dx} }}{{1 \!+\! {{\left\| {{z_1}} \right\|}^2} \!+\! {{\left\| {{z_2}} \right\|}^2}}} \!+\! \frac{{{z_1}(0)\hat \tau (0)}}{{1 \!+\! {{\left| {{z_1}(0)} \right|}^2} \!+\! {{\left| {{z_2}(0)} \right|}^2}}}} \! \right]} \!\! \right\}\\
			&{{\dot {\hat m}}_4} \!=\! {\mathcal{P}_{{{\bar m}_4}}} \!\! \left\{ \! {{\rho _4} \!\! \left[ \! {\frac{{\int_0^1 {{z_2}(x)\hat \tau (x)dx} }}{{1 \!+\! {{\left\| {{z_1}} \right\|}^2} \!+\! {{\left\| {{z_2}} \right\|}^2}}} \!+\! \frac{{{z_2}(0)\hat \tau (0)}}{{1 \!+\! {{\left| {{z_1}(0)} \right|}^2} \!+\! {{\left| {{z_2}(0)} \right|}^2}}}} \! \right]} \!\! \right\}
		\end{align}
		\end{small}
	\end{subequations}
	where ${\rho _1},{\rho _2},{\rho _3},{\rho _4}$ are the given positive gains, 
	we assume that the parameters estimated by the above equation also satisfy bounds \eqref{parameter bounds}, that expressed as
	\vspace{-1mm}
	\begin{equation} \label{hat parameter bounds}
		\left| {{\hat{m}_i}} \right| \le {{\bar m}_i},i = 1,2,3,4 .
		\vspace{-1mm}
	\end{equation}
	The adaptive laws operator $\mathcal{P}$  defined as 
	\vspace{-1mm}
	\begin{equation}
		{\mathcal{P}_{{{\bar m}_i}}}({\theta _i},{{\hat m}_i}) = \left\{ \begin{array}{l}
			0\;\;\; \text{if}\;\left| {{{\hat m}_i}} \right| \ge {{\bar m}_i}\;\& \;{\theta _i}{{\hat m}_i} \ge 0\\
			{\theta _i}\;\;\; \text{others}.
		\end{array} \right.
		\vspace{-1mm}
	\end{equation}
	
	\subsection{Adaptive control}
	Associative initial system \eqref{original_system}, parameter filters \eqref{filter1} and \eqref{filter2}, the dynamic states estimation system can be formulated as
	\begin{small}
	\begin{subequations}\label{dynamic system}
		\setlength{\abovedisplayskip}{5pt}
		\setlength{\belowdisplayskip}{5pt}
		\begin{align}
			{{\hat \psi }_t}(x,t) \!+\! \lambda {{\hat \psi }_x}(x,t) &={{\hat m}_1}\psi (x,t) \!+\! {{\hat m}_2}\varphi (x,t) \!+\! {n^T}(x,t)\dot {\hat M} \label{dynamic system a} \\ 
			{{\hat \varphi }_t}(x,t) \!-\! \mu{{\hat \varphi }_x}(x,t)&={{\hat m}_3}\psi (x,t)\!+\!{{\hat m}_4}\varphi(x,t) \!+\! {z^T}(x,t)\dot {\hat N} \label{dynamic system b} \\
			\hat \psi (0,t) &= q\varphi (0,t) \label{dynamic system c} \\
			\hat \varphi (1,t) &= U(t). \label{dynamic system d}
		\end{align}
	\end{subequations}
	\end{small}
	Then we apply the backstepping transformation to map the dynamic into target system, using the boundary control law
	\begin{small}
		\begin{equation}
			\setlength{\abovedisplayskip}{3pt}
			\setlength{\belowdisplayskip}{3pt}
			\begin{aligned}
				U(t) \!=\! \int_0^1 \!{\alpha (1,\xi ,t)} \hat \psi (\xi ,t)d\xi  \!+\! \int_0^1 \!{\beta (1,\xi ,t)} \hat \varphi (\xi ,t)d\xi ,
			\end{aligned}
		\end{equation}
	\end{small}
	where $\alpha$ and $\beta$ represent the resolvent kernel of the Volterra integral transformations, and satisfy the following kernel equation (\textit{Here t is omitted for brevity}). Note that the numerical kernel is obtained by solving the following PDE, contrasted with the approximated kernel.
	\begin{small}
	\begin{subequations}\label{kernel equation}
		\setlength{\abovedisplayskip}{2pt}
		\setlength{\belowdisplayskip}{2pt}
		\begin{align}
			\mu {\alpha _x}(x,\xi ) - \lambda {\alpha _\xi }(x,\xi ) &= ({{\hat m}_1} - {{\hat m}_4})\alpha (x,\xi ) + {{\hat m}_3}\beta (x,\xi )\\
			\mu {\beta _x}(x,\xi ) + \mu {\beta _\xi }(x,\xi ) &= {{\hat m}_2}\alpha (x,\xi )\\
			\beta (x,0) &= q\frac{\lambda }{\mu }\alpha (x,0)\\
			\alpha (x,x) &=  - \frac{{{{\hat m}_3}}}{{\lambda  + \mu }},
		\end{align}
	\end{subequations}
	\end{small}
	 	%\resizebox{0.9\hsize}{!}{$
	defined on
	\vspace{-1mm}
	\begin{align}\label{region T}
		T=\left\{ {(x,\xi )\left| {0 \le \xi  \le x \le 1} \right.} \right\} \times \left\{ {t \ge 0} \right\}.
		\vspace{-1mm}
	\end{align} 
	Moreover, the backstepping kernel is invertible with inverse
	\begin{small}
		\begin{subequations}\label{inverse kernel}
			\setlength{\abovedisplayskip}{5pt}
			\setlength{\belowdisplayskip}{5pt}
			\begin{align}
				{\alpha _I}(x,\xi,t ) = \alpha (x,\xi,t ) + \int_\xi ^x {\beta (x,s,t){\alpha _I}(s,\xi,t )ds} \\
				{\beta _I}(x,\xi,t ) = \beta (x,\xi,t ) + \int_\xi ^x {\beta (x,s,t){\beta _I}(s,\xi,t )ds} .
			\end{align}
		\end{subequations}
	\end{small}
	With the above kernel function, we define the following backstepping transformation 
	\begin{small}
	\begin{subequations}\label{backstepping transformation}
		\setlength{\abovedisplayskip}{5pt}
		\setlength{\belowdisplayskip}{5pt}
		\begin{align}
			&f(x,t) = \hat \psi (x,t)\\
			&\begin{aligned}
					h(x,t) = \Gamma [\hat \psi ,\hat \varphi ](x) &= \hat \varphi (x,t)- \int_0^x {\alpha (x,\xi ,t)} \hat \psi (\xi ,t)d\xi  \\
					& - \int_0^x {\beta (x,\xi ,t)} \hat \varphi (\xi ,t)d\xi ,
			\end{aligned}
		\end{align}
	\end{subequations}
	\end{small}
	the inverse form of this transformation is expressed as
	\begin{small}
		\begin{subequations}\label{backstepping transformation inverse}
			\setlength{\abovedisplayskip}{5pt}
			\setlength{\belowdisplayskip}{5pt}
			\begin{align}
				&\hat \psi (x,t) = f(x,t)\\
				&\begin{aligned}
					\hat \varphi (x,t) = {\Gamma ^{ - 1}}[f,h](x) &\!=\! h(x,t) \!+\! \int_0^x \! {{\alpha _I}(x,\xi ,t)} f(\xi ,t)d\xi \\
					&\!+ \int_0^x {{\beta _I}(x,\xi ,t)} h(\xi ,t)d\xi ,
				\end{aligned}
			\end{align}
		\end{subequations}
	\end{small}
	where $\Gamma$ and ${\Gamma^{ -1}}$ denote the transformation operators.
	
    To illustrate the boundedness of the operators defined above, we introduce a lemma as follows.
    \begin{lemma} \label{operator boundary}
    	For every time $t \ge 0$, we introduce the following properties for the operators defined in \eqref{region T}
    	\begin{small}
    		\begin{subequations}\label{kernel boundary}
    			\begin{align}
    				&\left| {\alpha (x,\xi ,t)} \right| \le \bar \alpha ,   \;\;\;\;\ 
    				\left| {\beta (x,\xi ,t)} \right| \le \bar \beta, \label{lemma2a} \\
    				&\left| {{\alpha _I}(x,\xi ,t)} \right| \le {{\bar \alpha }_I}, \;\;
    				\left| {{\beta _I}(x,\xi ,t)} \right| \le {{\bar \beta }_I}, \label{lemma2b} \\
    				&\left| {{\alpha _t}(x,\xi ,t)} \right|,  \   
    				\left| {{\beta _t}(x,\xi ,t)} \right| \in {L_2}. \label{lemma2d}
    			\end{align}
    		\end{subequations}
    	\end{small}
    \end{lemma}  
    \begin{IEEEproof}
    	For every time $t$, it has been proved in \cite{coron2013local} that $\alpha\left(x,\xi,t\right)$ and $\beta\left(x,\xi,t\right)$ are bounded and unique, and meet the following inequalities for all $(x,\xi ,t) \in T$:
    	\begin{small}
    		\begin{subequations}
    			\setlength{\abovedisplayskip}{3pt}
    			\setlength{\belowdisplayskip}{3pt}
    			\begin{align}
    				\left| {\alpha\left( {x,\xi ,t} \right)} \right|& \le {F_1}\left( {{{\hat m}_1}, \ldots {{\hat m}_4}} \right),  \label{proof2a}\\
    				\left| {\beta\left( {x,\xi ,t} \right)} \right| &\le {F_2}\left( {{{\hat m}_1}, \ldots {{\hat m}_4}} \right), \label{proof2b}\\
    				\left| {{\alpha_t}\left( {x,\xi ,t} \right)} \right| &\le {G_1}\left( {\left| {{{\dot {\hat m}}_1}} \right| + \left| {{{\dot{ \hat m}}_2}} \right| + \left| {{{\dot {\hat m}}_3}} \right| + \left| {{{\dot {\hat m}}_4}} \right|} \right), \label{proof2c}\\
    				\left| {{\beta_t}\left( {x,\xi ,t} \right)} \right| &\le {G_2}\left( {\left| {{{\dot {\hat m}}_1}} \right| + \left| {{{\dot {\hat m}}_2}} \right| + \left| {{{\dot {\hat m}}_3}} \right| + \left| {{{\dot {\hat m}}_4}} \right| } \right), \label{proof2d}
    			\end{align}
    		\end{subequations}
    	\end{small}
    	where $F_1\left(\cdot\right)$ and $F_2\left(\cdot\right)$ are continuous functions concerning the variables therein. Since $\hat m_i,i=1,\ldots,4$ are compact due to the projection ${\mathcal{P}_{{{\bar m}_i}}}$, $F_1\left(\cdot\right)$ and $F_2\left(\cdot\right)$ are bounded. Let $\bar \alpha$ and $\bar \beta$ be the maximum values of $F_1\left(\cdot\right)$ and $F_2\left(\cdot\right)$, respectively,  then \eqref{lemma2a} must hold. 
    	With the successive approximation method, the solution $\left({\alpha}_I,\beta_I\right)$ to equation \eqref{inverse kernel} is bounded whose bounds $\bar {\alpha}_I$ and $\bar {\beta}_I$ are determined by $\bar {\alpha}$, $\bar {\beta}$.
    	Again using the results in \cite{coron2013local}, the time derivative $\left(\alpha_t,\beta_t\right)$ satisfies \eqref{proof2c} and \eqref{proof2d}. Note  $\dot{ \hat{m}}_i , i=1...4$ is bounded convergent, hence $\left| {\alpha_t} \right|$, $\left| {\beta_t} \right|$ are bounded and convergent as well, \eqref{lemma2d} is established.
    \end{IEEEproof}
    
    By differentiating \eqref{backstepping transformation} with respect to time, inserting the dynamic system equation \eqref{dynamic system}, the kernel function transformation \eqref{kernel equation} and the boundary filters \eqref{filter1}, \eqref{filter2}, we can derive the target system equation.     
    
    In brief, this section introduces the above backstepping adaptive control process using the swapping design. The target system equations and the detailed derivation, as well as system stability proofs, are outlined in \cite[Section 4]{anfinsen2018adaptive} and will not be repeated here.
	
	It's worth noting that conventional backstepping adaptive controllers typically rely on numerical solutions to obtain the kernel gain.
	Undoubtedly, due to real-time estimation of unknown parameters, kernel functions represented as PDEs require iterative solving, leading to a rapid increase in computational cost as the control time domain extends. To address this, we will introduce a neural operator that streamlines kernel computation as a straightforward function evaluation, thereby accelerating the control process.

\section{Approximated kernel by neural operator}
    In this section, the DeepONet network is introduced, inspired by the universal approximation theorem for operators. This network can provide an approximate spatial gain kernel with arbitrary accuracy.
    We analyze its structure and study the existence and accuracy of the approximate backstepping kernel for adaptive control.
    \subsection{Operator learning with DeepONet}
    \subsubsection{The network structure of DeepONet}
    DeepONet is a general neural network framework for learning continuous nonlinear operators, which consists of an offline training stage followed by an online inference stage. In the offline stage, classical numerical methods are employed to solve for the target operator within the appropriate input space, and then train our network. In the online stage,the trained network is used to perform real-time inference, significantly enhancing the speed of inference.
    
    As depicted in Fig.\ref{deeponet structure}, DeepONet consists of two basic sub-networks, the branch net is used to extract the potential representation of the input function $u({x_i})$, corresponding to a set of sensor locations $\{{x_1},{x_2},...,{x_m}\}$, and the trunk net is used to extract the potential representation of the input coordinates $y$ at which the output functions are evaluated. The continuous and differentiable representations of the output function are then obtained by merging the potential representations extracted by each sub-network through dot product.
	
	Here ${\Phi}$ as the operator of the input function $u$, $\Phi(u)(\cdot)$ represents the corresponding output function, $\Phi(u)(y)$ denotes the corresponding neural operator value, and for any point $y$ in the defined domain of ${\Phi}(u)(\cdot)$, the DeepONet neural operator ${\Phi_{\cal N}}$ for approximating a nonlinear operator ${\Phi}$ is defined as follows
	\begin{equation}
		\setlength{\abovedisplayskip}{3pt}
		\setlength{\belowdisplayskip}{3pt}
		\begin{aligned}\label{deeponet}
			{\Phi _{\cal N}}(u)(y): = \sum\limits_{k = 1}^p {{g^{\cal N}}({u_m};\theta _1^{(k)}){f^{\cal N}}(y;\theta _2^{(k)})},
		\end{aligned}
	\end{equation}
   	where $p$ is the number of basis components, ${\theta _1^{(k)}}$, ${\theta _2^{(k)}}$ denote the weight coefficients of the trunk and branch net respectively. Here ${{g^{\cal N}}}$, ${{f^{\cal N}}}$ are NNs termed branch and trunk networks, which can be any neural network satisfying the classical universal approximation theorem, such as FNNs(fully-connected neural networks) and CNNs(convolutional neural networks).
   	
   	\subsubsection{Accuracy of approximation of neural operators}
   	It is proven that DeepONet satisfies the universal approximation of continuous operators as long as the branch and trunk neural networks satisfy the universal approximation theorem of compact set continuous functions. Now let's introduce this theorem.
   	\begin{theorem}\label{universal approximation theorem}(Universal approximation of continuous operators \cite[Theorem 2.1]{lu2021learning}).
	Let ${D_x} \subset {R^{{d_x}}}$, ${D_y} \subset {R^{{d_y}}}$ be compact sets of vectors $x \in {D_x}$ and $y \in {D_y}$. Let $D_u:{D_x} \to {D_u} \subset {R^{{d_u}}}$ is a compact set and $u \in {D_u}$, which be sets of continuous functions $u(x_i)$. let ${D_\Phi }:{D_y} \to {D_\Phi } \subset {R^{{d_\Phi }}}$ be sets of continuous functions of $\Phi (u)(y)$. Assume that $\Phi :{D_u} \to {D_\Phi }$ is a continuous operator. Then, for any $\varepsilon  > 0$, there exist ${m^*},{p^*} \in \mathbb{N}$ such that for each $m \ge {m^*},p \ge {p^*}$, there exist weight coefficients ${\theta _1^{(k)}}$, ${\theta _2^{(k)}}$, neural networks ${{g^{\cal N}}}$, ${{f^{\cal N}}}$, $k = 1,...,p$ and ${x_i} \in {D_x},i = 1,2,...,m$, there exists a neural operator ${\Phi _N}$ such that the following holds:
	\begin{equation}
		\setlength{\abovedisplayskip}{3pt}
		\setlength{\belowdisplayskip}{3pt}
		\begin{aligned}\label{universal approximation}
			\left| {\Phi (u)(y) - {\Phi _N}(u)}(y) \right| < \varepsilon, \;\;\;\;\forall u \in {D_u},y \in {D_y}
		\end{aligned}	
	\end{equation}
 
   	\end{theorem}
   	where $\textbf{u}_m = {\left[ {u\left( {{x_1}} \right){\rm{, }}u\left( {{x_2}} \right){\rm{, }}...{\rm{, }}u\left( {{x_m}} \right)} \right]^T}$.
   	
   	This means that by adjusting the input grid points, the size of the training samples or the structure of the sub-network, etc., it is theoretically possible to obtain neural operators that can provide approximate kernel gains with arbitrary accuracy.

    \subsection{Backstepping kernel with neural operators} 
    In equation \eqref{kernel equation}, there are four time-varying constants, which require the kernel operators to be solved repeatedly to account for these variations. We consider employing an approximate operators rather than seeking an exact solution, introducing the DeepONet network to learn gain function mapping as an alternative to PDE numerical solvers. By formalizing neural operators, we can effectively learn the mapping of various system parameter values without necessitating retraining. The gain kernel can be obtained by evaluating at certain parameter values, thereby significantly accelerating the solution process.
    
    Define the following notations, $ \varpi  = \{ {m_1},{m_2},{m_3},{m_4}\} $, its boundary $\bar \varpi \ge \left| \varpi \right|$, then denote by ${\Phi _\alpha }$, ${\Phi _\beta }$ the operators that maps $\varpi$ to the kernel $\alpha$, $\beta$ that satisfies \eqref{kernel equation}. Expressed as 
    \begin{subequations}\label{kernel maps}
    	\setlength{\abovedisplayskip}{5pt}
    	\setlength{\belowdisplayskip}{5pt}
    	\begin{align}
    		&{\Phi _\alpha }:\varpi \to \alpha ,\;\; \alpha = {\Phi _\alpha }(\varpi),\\
    		&{\Phi _\beta }:\varpi \to \beta ,\;\; \beta = {\Phi _\beta }(\varpi).
    	\end{align}
    \end{subequations}
    According to the proof of Lemma \ref{operator boundary}, the backstepping kernel operator $\alpha$ and $\beta$ satisfy \eqref{proof2a}-\eqref{proof2b}, with their boundaries $\bar{\alpha}$ and $\bar{\beta}$ determined by the continuous functions $F_1\left(\cdot\right)$ and $F_2\left(\cdot\right)$ of the variable $\varpi$. 
    It was shown in \cite{bhan2023neural} that when the boundary of $\varpi$ exists and is Lipschitz continuous, we obtain:
    \begin{subequations}\label{Lipschitz}
    	\setlength{\abovedisplayskip}{5pt}
    	\setlength{\belowdisplayskip}{5pt}
    	\begin{align}
    		&{\left\| {{\alpha _1} - {\alpha _2}} \right\|_ \infty } < {C_{\bar \alpha }}{\left\| {{\varpi _1} - {\varpi _2}} \right\|_ \infty }\\
    		&{\left\| {{\beta _1} - {\beta _2}} \right\|_ \infty } < {C_{\bar \beta }}{\left\| {{\varpi _1} - {\varpi _2}} \right\|_ \infty }
    	\end{align}
    \end{subequations}
    This leads to the operator ${\Phi _\alpha }$, ${\Phi _\beta }$ possessing the Lipschitz property
    \begin{subequations}\label{operator Lipschitz}
    	\setlength{\abovedisplayskip}{5pt}
    	\setlength{\belowdisplayskip}{5pt}
    	\begin{align}
    		{\left\| {{\Phi _\alpha }({\varpi _1}) - {\Phi _\alpha }({\varpi _2})} \right\|_\infty } < {C_{\bar \alpha }}{\left\| {{\varpi _1} - {\varpi _2}} \right\|_\infty }\\
    		{\left\| {{\Phi _\beta }({\varpi _1}) - {\Phi _\beta }({\varpi _2})} \right\|_\infty } < {C_{\bar \beta }}{\left\| {{\varpi _1} - {\varpi _2}} \right\|_\infty },
    	\end{align}
    \end{subequations}
    where ${C_{\bar \alpha }}$ and ${C_{\bar \beta }}$ denote the Lipschitz constants,  which can be expressed as expressions for the boundaries $\bar \alpha$, $\bar \beta$ of the variable function.

    Then combined with DeepONet universal approximation theorem (The theorem \ref{universal approximation theorem}), we obtain
    \begin{subequations}\label{kernel maps approximation}
    	\setlength{\abovedisplayskip}{5pt}
    	\setlength{\belowdisplayskip}{5pt}
    	\begin{align}
    			\left| {{\Phi _\alpha }(\varpi ) - {\Phi _{\mathcal{N}\alpha }}(\varpi )} \right| &< \varepsilon, \\
    			\left| {{\Phi _\beta }(\varpi ) - {\Phi _{\mathcal{N}\beta }}(\varpi )} \right| &< \varepsilon,
    	\end{align}
    \end{subequations}
    where ${\varpi _m} = {(\varpi ({x_1}),\varpi ({x_2}),...,\varpi ({x_m}))^T}$ with corresponding ${x_i} \in {D_x}$, for $i \!=\! 1,2,...,m$. 
    
    \begin{remark}
    	The above was inspired by \cite{lamarque2024gain}.
    	Moreover, for target systems containing derivative terms, the differential forms  of their backstepping transformation kernels also correspond to specific DeepONet operators, whose existence and Lipschitz continuity is analyzed in \cite[Theorem 2]{lamarque2024gain}.
    \end{remark}
    
    The inequalities \eqref{kernel maps approximation} demonstrates that the backstepping kernels $\alpha$ and $\beta$ are approximable, and that the error between the approximated kernel and the exact kernel is within a small domain. Furthermore, for the stability of the system under kernel approximation with DeepONet, we will provide rigorous proof in the following section.
    
    The approximated kernels here are denoted as
    \vspace{-1mm}
    \begin{equation}\label{hat kernel maps}
    	\begin{aligned}
    		{\hat{\alpha} } = {\Phi _{\mathcal{N}\alpha }}(\varpi),
    		\quad \quad 
    		{\hat{ \beta} } = {\Phi _{\mathcal{N}\beta }}(\varpi).
    	\end{aligned}
    \end{equation}
%    From \eqref{inverse kernel}, its inverse
%    \begin{subequations}\label{approximate inverse kernel}
%    	 \begin{align}
%    		{\hat{\alpha} _I}(x,\xi ) = \hat{\alpha } (x,\xi ) + \int_\xi ^x {\hat{\beta} (x,s){\hat{\alpha} _I}(s,\xi )ds} \\
%    		{\hat{\beta} _I}(x,\xi ) = \hat{\beta} (x,\xi ) + \int_\xi ^x {\hat{\beta} (x,s){\hat{\beta} _I}(s,\xi )ds} .
%    	 \end{align}
%	\end{subequations}
	
\section{Approximate kernels are applied in \\ adaptive control}
	In this section, we utilize approximate kernels \eqref{hat kernel maps} instead of analytical computations for control using swapping design. (Some formulas have been defined in section \ref{basic adaptive control}). 
	\subsection{Adaptive control}
	
	There will be an approximation error between the approximate kernel and the exact kernel, which we define as
	\vspace{-1mm}
	\begin{equation}\label{kernel error}
		\begin{aligned}
			\tilde \alpha  = \alpha  - \hat \alpha, \quad \tilde \beta  = \beta  - \hat \beta .
		\end{aligned}
		\vspace{-1mm}
	\end{equation} 
	The sum of the boundary integrals of the operator kernel forms the approximate controller, denoted as
	\begin{small}
		\begin{equation}
			\begin{aligned}\label{approximate controller}
				\hat U(t) \!=\! \int_0^1 \!{\hat \alpha (1,\xi ,t)} \hat \psi (\xi ,t)d\xi  \!+\! \int_0^1 \!{\hat \beta (1,\xi ,t)} \hat \varphi (\xi ,t)d\xi.
			\end{aligned}
		\end{equation}
	\end{small}
	
	Obviously, the DeepONet operator kernel brings the approximation error term, resulting in an accumulation of errors between the approximate controller and the original controller. Then, the new transformation terms and transformation operators are obtained through PDE backstepping design.
    
	As follows, we consider the new backstepping transformation with the approximate kernel 
	\begin{small}
		\begin{subequations}\label{approximate backstepping transformation}
			\begin{align}
				&{f_2}(x,t) = \hat \psi (x,t) \label{approximate backstepping transformation a}\\
				&\begin{aligned}\label{approximate backstepping transformation b}
				 	{h_2}(x,t) = {\Gamma _2}[\hat \psi ,\hat \varphi ](x) &= \hat \varphi(x,t) - 		\int_0^x  {\hat \alpha (x,\xi ,t)} \hat \psi (\xi ,t)d\xi \\
				    &- \int_0^x {\hat \beta (x,\xi ,t)} \hat \varphi (\xi ,t)d\xi ,
				\end{aligned}	
			\end{align}
		\end{subequations}
	\end{small}
	and its inverse transformation
		\begin{subequations}\label{approximate inverse transformation}
		\setlength{\abovedisplayskip}{2pt}
		\setlength{\belowdisplayskip}{2pt}
		\begin{footnotesize}
			\begin{align}
				&\hat \psi (x,t) = {f_2}(x,t) \label{approximate inverse transformation a}\\
				&\begin{aligned}\label{approximate inverse transformation b}
					\hat \varphi (x,t) = \Gamma _2^{ - 1}[{f_2},{h_2}](x) &= {h_2}(x,t) \!+\! \int_0^x \! {{{\hat \alpha }_I}(x,\xi ,t)} {f_2}(\xi ,t)d\xi \\
					&+ \int_0^x {{{\hat \beta }_I}(x,\xi ,t)} {h_2}(\xi ,t)d\xi  ,
				\end{aligned}
			\end{align}
		\end{footnotesize}
		\end{subequations}
	where ${\Gamma _2}$ and $\Gamma _2^{ - 1}$ are defined as new conversion operators.
	Similar to \eqref{inverse kernel}, the approximate inversion kernel ${{\hat \alpha }_I}$, ${{\hat \alpha }_I}$ satisfies the following form 
	\begin{subequations}\label{approximate inverse kernel}
		\setlength{\abovedisplayskip}{2pt}
		\setlength{\belowdisplayskip}{2pt}
		\begin{small}
			\begin{align}
				{\hat{\alpha} _I}(x,\xi,t ) = \hat{\alpha} (x,\xi,t ) + \int_\xi ^x {\hat{\beta} (x,s,t){\hat{\alpha} _I}(s,\xi,t )ds} \\
				{\hat{\beta} _I}(x,\xi,t ) = \hat{\beta} (x,\xi,t ) + \int_\xi ^x {\hat{\beta} (x,s,t){\hat{\beta} _I}(s,\xi,t )ds} .
			\end{align}
		\end{small}
	\end{subequations}
	
	The error introduced by the approximate kernel adds further perturbation terms to the target system, making its composition more complex and posing challenges to theoretical stability analysis.
	Considering the property \eqref{kernel maps approximation} obtained from the universal approximation theorem, the error in the operator kernel can be confined to a small domain. This suggests that the approximation operator has a similar boundary as in Eqs. \eqref{kernel boundary}. Formulated as the following lemma.
	\begin{lemma}
		For all $(x,\xi ,t) \in T$ at \eqref{region T}, the approximate kernel is bounded and the approximate errors are square integrable.
		\begin{subequations}\label{approximate kernel boundary}
			\begin{small}
				\begin{align}
					&\left| {\tilde \alpha } \right| \le \bar {\tilde \alpha} , \quad \left| {{{\tilde \alpha }_I}} \right| \le {{\bar {\tilde \alpha} }_I}, \quad \left| {\hat \alpha } \right| \le \bar {\hat \alpha} , \quad \left| {{{\hat \alpha }_I}} \right| \le {{\bar {\hat \alpha} }_I}\\
					&\left| {\tilde \beta } \right| \le \bar {\tilde \beta} , \quad \left| {{{{\tilde \beta} }_I}} \right| \le {{\bar {\tilde \beta} }_I}, \quad \left| {\hat \beta } \right| \le \bar {\hat \beta} , \quad \left| {{{\hat \beta }_I}} \right| \le {{\bar {\hat \beta} }_I}\\
					&\left| {{\partial _t}\tilde \alpha } \right|, \quad \left| {{\partial _t}{{\tilde \alpha }_I}} \right|, \quad \left| {{\partial _t}\hat \alpha } \right|, \quad \left| {{\partial _t}{{\hat \alpha }_I}} \right| \in {L_2}\\
					&\left| {{\partial _t}\tilde \beta } \right|, \quad \left| {{\partial _t}{{\tilde \beta }_I}} \right|, \quad \left| {{\partial _t}\hat \beta } \right|, \quad \left| {{\partial _t}{{\hat \beta }_I}} \right| \in {L_2}.
				\end{align}
			\end{small}
		\end{subequations}
	\end{lemma}
	\begin{IEEEproof}
		Similar to the proof of Lemma \ref{operator boundary}.
	\end{IEEEproof}
	
%	\begin{lemma}\label{lemma target system}
	The backstepping transformation \eqref{approximate backstepping transformation} and the approximate controller \eqref{approximate controller}, with the approximate backstepping kernels satisfying \eqref{hat kernel maps} and \eqref{approximate inverse kernel}, map the dynamics \eqref{dynamic system} into the new target system as follows (\textit{Here t is omitted for brevity}).
	\begin{small} 
		\begin{subequations}\label{target system2}
			\begin{align}
				&\begin{aligned}\label{target system2 a}
					&{\partial _t}{f_2}(x) + \lambda {\partial _x}{f_2}(x) = {n^T}(x)\dot {\hat M}\\
					&\ +{{\hat m}_1}{f_2}(x) + {{\hat m}_1}\hat e(x) + {{\hat m}_2}\hat \tau (x) + {{\hat m}_2}h_2(x)\\
					&\ +{{\hat m}_2}\int_0^x {{{\hat \alpha }_I}(x,\xi )} {f_2}(\xi )d\xi + {{\hat m}_2}\int_0^x {{{\hat \beta }_I}(x,\xi )} {h_2}(\xi )d\xi  ,
				\end{aligned}\\
				&\begin{aligned}\label{target system2 b}
					&{\partial _t}{h_2}(x) - \mu {\partial _x}{h_2}(x)\\
					&\ = [{{\hat m}_3}\hat e(x) + {{\hat m}_4}\hat \tau (x) + {{\hat m}_4}\Gamma _2^{ - 1}[{f_2},{h_2}](x)] + [{z^T}(x) \dot {\hat N}]\\
					&\ - \int_0^x {\left[ {{{\hat \alpha }_t}(x,\xi ) - \lambda {{\tilde \alpha }_\xi }(x,\xi ) + \mu {{\tilde \alpha }_x}(x,\xi ) + {{\hat m}_4}\alpha (x,\xi )} \right.} \\
						&\;\;\;\;\;\;\;\;\; \left. { - {{\hat m}_1}\tilde \alpha (x,\xi ) - {{\hat m}_3}\tilde \beta (x,\xi )} \right]{f_2}(\xi )d\xi  \\
					&\ - \int_0^x {\left[ {{{\hat \beta }_t}(x,\xi ) + \mu {{\tilde \beta }_\xi }(x,\xi ) + \mu {{\tilde \beta }_x}(x,\xi ) + {{\hat m}_4}\hat \beta (x,\xi )} \right.} \\
						&\;\;\;\;\;\;\;\;\; \left. { - {{\hat m}_2}\tilde \alpha (x,\xi)} \right]\Gamma _2^{ - 1}[{f_2},{h_2}](\xi)d\xi \\
					&\ - \int_0^x {\hat \alpha (x,\xi)} [{{\hat m}_1}\hat e(\xi) + {{\hat m}_2}\hat \tau (\xi)]d\xi \\
					&\ - \int_0^x {\hat \beta (x,\xi)} [{{\hat m}_3}\hat e(\xi) + {{\hat m}_4}\hat \tau (\xi)]d\xi \\
					&\ - \int_0^x {\hat \alpha (x,\xi)[} {n^T}(\xi)\dot {\hat M}]d\xi  - \int_0^x {\hat \beta (x,\xi)} [{z^T}(\xi)\dot {\hat N}]d\xi \\
					&\ - (\lambda  + \mu )\tilde \alpha (x,x){f_2}(x) - \lambda \hat \alpha (x,0)[q{h_2}(0) + q\hat \tau (0)] \\
					&\ + \mu \hat \beta (x,0){h_2}(0)  ,
				\end{aligned} \\
				& {f_2}(0) = q{h_2}(0) + q\hat \tau (0)  \label{target system2 c}, \\
				& {h_2}(1) = U - \hat U  \label{target system2 d},
			\end{align}
		\end{subequations}
	\end{small}
	where $n$ and $z$ represent the boundary filters, $\hat e$ and $\hat \tau$ denote the state estimation error as defined in \eqref{adaptive errors}. ${\dot {\hat M}}$ and ${\dot {\hat N}}$ signify the derivatives of the time-varying parameters of the system, satisfying
	$	\dot {\hat M}(t) = {[{{\dot {\hat m}}_1}(t),{{\dot {\hat m}}_2}(t)]^T},\dot {\hat N}(t) = {[{{\dot {\hat m}}_3}(t),{{\dot {\hat m}}_4}(t)]^T} $.
%	\end{lemma}
	\begin{IEEEproof}
		The detailed derivation of the target system \eqref{target system2} obtained using the approximate backstepping kernel is given in the Appendix \ref{appendix B}.
	\end{IEEEproof}
	
	To summarize, we use the approximation gain kernel for the adaptive control, where the kernel error adds perturbation terms to the controller and target system, and thus the stability and convergence of the controlled system needs to be further analyzed, which we describe in detail in the next subsection. 
	
	On the other hand, we introduce neural operators for offline learning, once such solution operators are learned, they can directly evaluate any new input queries with different system parameters without solving the kernel differential equations. This makes the NN solver orders of magnitude faster than conventional numerical solvers, rendering it more suitable for real-time adaptive control applications.
	
	\subsection{Convergence analysis}
%	In the previous subsection, we apply the approximate kernels as control gain to the adaptive control process. 
	By bounding the operator error, we can then treat the NN as a disturbance effect on the adaptive controller and analyze the resulting stability properties of the closed-loop system that using the approximate controller.
	% The resulting kernel error may affect the controller and system state, 
	Specifically, we establish Theorem \ref{converge} based on the Lyapunov stability principle, as demonstrated in proof \ref{converge proof}. Further details are provided in Appendix \ref{appendix C}.
	\begin{theorem}\label{converge}
		Consider the origin system \eqref{original_system} and the state estimates ${\hat \psi }$,${\hat \varphi }$ defined from \eqref{adaptive states} using the boundary filters \eqref{filter1}-\eqref{filter2} and the adaptive control law \eqref{adaptive laws}. Consider the approximated controller $\hat U(t) = \int_0^1 {\hat \alpha (1,\xi ,t)} \hat \psi (\xi ,t)d\xi  + \int_0^1 {\hat \beta (1,\xi ,t)} \hat \varphi (\xi ,t)d\xi$, where $\hat\alpha,\hat\beta$ are the approximated kernels \eqref{hat kernel maps} obtained from DeepONet. Then all signals within the closed-loop system exhibit boundedness and integrability in the L2-sense. Furthermore, $\psi (x, \cdot ),\varphi (x, \cdot ) \in {L_\infty } \cap {L_2}$ and $\psi (x, \cdot ),\varphi (x, \cdot ) \to 0$ for all $x \in [0,1]$.
	\end{theorem}
                                                                                                                                                                        
	\begin{IEEEproof} \label{converge proof}
	Consider the following Lyapunov candidate sub-functions as
	\begin{subequations} \label{Lyapunov sub-function}
		\begin{align}\label{V11}
			&{V_1} = \left\| {{n_1}} \right\|_{ - a}^2 \; {V_2} = \left\| {{n_2}} \right\|_{ - a}^2 \; {V_3} = \left\| {{n_3}} \right\|_{ - a}^2\\
			&{V_4} = \left\| {{z_1}} \right\|_b^2  \;\;\;\; {V_5} = \left\| {{z_2}} \right\|_b^2 \;\;\;\; {V_6} = \left\| {{z_3}} \right\|_b^2\\
			&{V_7} = \left\| {{f_2}} \right\|_{ - a}^2 \;\; {V_8} = \left\| {{h_2}} \right\|_b^2\\
			&{V_{9}} = \left\| e \right\|_{ - a}^2 \;\;\;\;  {V_{10}} = \left\| \tau  \right\|_b^2,
		\end{align}
	\end{subequations}
	where we set $a$ and $b$ as constants greater than one for ease of analysis below. Then, we consider Lyapunov candidate functions that are linear combinations of sub-functions
	\begin{equation}
		\setlength{\abovedisplayskip}{3pt}
		\setlength{\belowdisplayskip}{3pt}
		\begin{aligned}
			{V_{11}} &= {d_1}[{V_1} + {V_2}] + {d_2}{V_3} + {d_3}[{V_4} + {V_5}]\\
			&+ {d_4}{V_6} + {d_5}{V_7} + {d_6}{V_8} + {d_7}{V_9} + {d_8}{V_{10}},
		\end{aligned}
	\end{equation}
	where $d_i, i\!=\!1,...,8$ are positive constants. 
	From Appendix \ref{appendix C}, 
	the differentiation of all Lyapunov sub-functions is transformed into linear combinations of candidate functions. Their coefficient multipliers are transformed into constant ($c_i, \, i=1...10$) or bounded convergence ($l_i, \, i=1...14$) terms.
	Further integration leads to the following form
	\begin{equation} \label{V11 de}
		\setlength{\abovedisplayskip}{5pt}
		\setlength{\belowdisplayskip}{5pt}
		\begin{small}
			\begin{aligned}  
				&{{\dot V}_{11}} \le  - {d_1}(\lambda a - 1)[{V_1} + {V_2}] - {d_2}\lambda a{V_3}\\
				&- {d_3}(\mu b - 1)[{V_4} + {V_5}] - {d_4}\mu b{V_6}\\
				&- [{d_5}(\lambda a - {c_7}) - {d_1}{c_1}{e^a} - {d_3}{c_3}{e^{b + a}} - {d_4}{c_5}{e^{b + a}}]{V_7}\\
				&- [{d_6}(\mu b - {c_8}) - {d_1}{c_2} - {d_3}{c_4}{e^b} - {d_4}{c_6}{e^b} - 2{d_5}]{V_8}\\
				&- [{d_7}\lambda a - 4{d_1}{e^a} - 4{d_3}{e^{b + a}} - 2{d_5}{e^a} - {d_6}{c_9}{e^{b + a}}]{V_9}\\
				&- [{d_8}\mu b - 4{d_1} - 4{d_3}{e^b} - 2{d_5} - {d_6}{c_{10}}{e^b}]{V_{10}}\\
				&+ \left( {{d_1}{l_1} + {d_3}{l_4} + {d_5}{l_6} + {d_6}{l_9}} \right)[{V_1} + {V_2}] + {d_6}{l_{11}}{V_7}\\
				&+ ({d_1}{l_2} + {d_3}{l_5} + {d_5}{l_7} + {d_6}{l_{10}})[{V_4} + {V_5}] 
				+ {d_6}{l_{12}}{V_8}\\
				&- ( - 2{d_2}\lambda  - 2{d_5}\lambda {q^2} + {d_6}\mu )h_2^2(0)\\
				&+ (4{d_2}\lambda  + 4{d_5}\lambda {q^2} + 4{d_6}{e^b}){\tau ^2}(0)+ {d_6}{l_{14}}h_2^2(0)\\
				&- ({d_3}\mu  - {d_2}{l_3} - {d_5}{l_8} - {d_6}{l_{13}}){z^2}(0) ,
			\end{aligned}
		\end{small}
	\end{equation}
	the constant coefficient multipliers for these subsystems are set to
	\vspace{-1mm}
	\begin{equation} \label{d set}
		\begin{aligned}
			&{d_1} = {d_2} = {e^{ - a}}, \;\;\;\;\;\; {d_3}= {d_4} = {e^{ - b - a}}, \;\;\;\;\;\; {d_5} = 1,\\
			&{d_6} = \frac{{2{d_2}\lambda  + 2{d_5}\lambda {q^2}}}{\mu }, \;\;\;\;\;\; {d_7} = {e^{b + a}}, \;\;\;\;\;\; {d_8} = {e^b},
		\end{aligned}
		\vspace{-1mm}
	\end{equation}                                               
	and we choose
    \begin{subequations}
    	\setlength{\abovedisplayskip}{5pt}
    	\setlength{\belowdisplayskip}{5pt}
	\begin{footnotesize}
		 \begin{align}
			& \begin{aligned}\label{b set}
				\setlength{\abovedisplayskip}{1pt}
				\setlength{\belowdisplayskip}{1pt}
				b \!> \! \max \! \left\{ {\!1,\! \frac{1}{\mu },\! \frac{{{e^{\! - \!a}}{(c_2+c_4+c_6)} \!\!+\! {d_6}{c_8}}}{{\mu {d_6}}},\!} \!\right.\left. \! {\frac{{8{e^{ \!-\! b \!-\! a}} \!\!+\! 2{e^{ \!-\! b}} \!\!+\! {d_6}{c_{10}}}}{\mu }} \!\! \right\} ,
			\end{aligned} \\
			& \begin{aligned}\label{a set}
				\setlength{\abovedisplayskip}{1pt}
				\setlength{\belowdisplayskip}{1pt}
				a \!>\! \max \left\{ {1,\frac{1}{\lambda },\frac{{{c_1} \!+\! {c_3} \!+\! {c_5} \!+\! {c_7}}}{\lambda },} \right.\left. {\frac{{8{e^{ - \delta  - \alpha }} \!+\! 2{e^{ - \alpha }} \!+\! {d_6}{h_9}}}{\lambda }} \right\},
			\end{aligned} 
		\end{align}
	\end{footnotesize}
    \end{subequations}
    then we can obtain
    \begin{small}
    	\setlength{\abovedisplayskip}{5pt}
    	\setlength{\belowdisplayskip}{5pt}
    	\begin{equation}\label{V11 simp}
    		\begin{aligned}
    			{{\dot V}_{11}}\left( t \right) &\le  - C{V_{11}}\left( t \right) + {F}\left( t \right){V_{11}}\left( t \right) + {G}\left( t \right)\\
    			&- \left( {{d_3}\mu  - {H}\left( t \right)} \right){z^2}(0).
    		\end{aligned}
    	\end{equation}
    \end{small}
    The above formula corresponds to equation \eqref{V11 de}, $C$ is a constant greater than 0, determined by the first eight coefficients of equation \eqref{V11 de}, ${F}\left( t \right)$ is integrable function corresponding to terms 9-12, and ${G}\left( t \right) = (4{d_2}\lambda  + 4{d_5}\lambda {q^2} + 4{d_6}{e^b}){\tau ^2}(0)+{d_6}{l_{14}}h_2^2(0)$ converges in finite time, ${H}\left( t \right) = {d_2}{l_3} + {d_5}{l_8} + {d_6}{l_{13}}$ are also integrable functions.
    Using the comparison principle yields that
    \vspace{-2mm}
    \begin{equation}\label{comparison principle}
    	\begin{aligned}
    		&{e^{Ct}}{V_{11}}\left( t \right) \le {V_{11}}\left( 0 \right){e^{\int_0^t {F(s)ds} }}\\
    		&+ \int_0^t {{e^{Cs + \int_s^t {F(x)dx} }} [G(s) \!-\! \left( {{d_3}\mu  \!-\! H\left( s \right)} \right){z^2}(0, s)]ds} \\
    		&\le{V_{11}}\left( 0 \right){e^{\int_0^t {F(s)ds}}} + {e^{\int_0^t {F(s)ds} }}\int_0^t {{e^{Cs}}G(s)ds} \\
    		&+ {e^{\int_0^t {F(s)ds} }}\int_0^t {{e^{Cs}}\left( {{d_3}\mu  - H\left( s \right)} \right){z^2}(0,s)ds} .
    	\end{aligned}
    	\vspace{-2mm}
    \end{equation}
    Since ${F}\left( t \right)$ is integrable function, ${e^{\int_0^t {F(s)ds} }}$ is bounded, and ${G}\left( t \right)$ converges in finite time, then there $\exists {T_g},\forall t > T_g \to G(t) = 0$, implying that $\int_0^t {{e^{Cs}}G(s)ds} $ is also bounded. Additionally, ${H\left( t \right)}$ converges to 0 in finite time, so that $\exists {T_H}$, ${d_3}\mu  - H\left( t \right) > 0$ and ${d_3}\mu  - H\left( t \right) \to {d_3}\mu$ while $t > {T_H}$. This implies that $\int_0^t {{e^{Cs}}\left( {{d_3}\mu  - H\left( s \right)} \right){z^2}(0,s)ds}$ must also be bounded. Hence, ${e^{Ct}}{V_{11}}\left( t \right)$ is bounded as $t \to \infty$, which means that ${e^{Ct}}{V_{11}}\left( t \right)$ converges to zero for finite time.
    
    Now therefore, $\left\| n \right\|,\left\| z \right\|,\left\| {{f_2}} \right\|,\left\| {{h_2}} \right\| \in {L_2} \cap {L_\infty }$ are established. Then $\left\| \psi  \right\|,\left\| \varphi  \right\| \in {L_2} \cap {L_\infty }$ can be obtained directly by formulas \eqref{real status boundary} and \eqref{estimated error boundary}, as $\psi (x, \cdot ),\varphi (x, \cdot ) \to 0$ for all $x \in [0,1]$.
	
	\end{IEEEproof}
	
	% more
	In short, this section analyses the design process of the controller incorporating neural operator, giving a rigorous analytical proof of the stability of the closed-loop system. Combined with the properties of the DeepONet approximate kernel in Section 3, it is shown that by applying a kernel function of sufficient approximate accuracy in the adaptive control process, convergence of the closed-loop system state can be guaranteed.
	
	Then, we will verify the effectiveness of the proposed NN-adaptive control scheme through simulation experiments and analyze the role of the well-trained neural operators in greatly accelerating the control.

\section{Simulation}
	We consider the simulation of system \eqref{original_system} with the parameter settings
	$\lambda  = 1, \mu  = 1,  q = 1$, 
	and the simulated actual values of the unknown parameters as:
	\begin{equation} \label{actual values}
		\setlength{\abovedisplayskip}{2pt}
		\setlength{\belowdisplayskip}{2pt}
		{m_1} = 0.4,  \quad {m_2} = 0.6, \quad {m_3} = 1.0, \quad {m_4} = 0.8.
	\end{equation}
	 Additionally, we set the prior parameter bounds to 
	${{\bar m}_1} = {{\bar m}_2} = {{\bar m}_3} = {{\bar m}_4} = 4$,
	and the adaptive law gains to
	${\rho _1} = {\rho _2} = {\rho _3} = {\rho _4} = 40$.
	Then we set the initial states of the system as
	$\psi (x,0) = 1,\varphi (x,0) =2 \sin (10x) + 3$,
	the initial value of the estimated parameters as
	${{\hat m}_1}(0) = {{\hat m}_2}(0) = {{\hat m}_3}(0) = {{\hat m}_4}(0) = 0.4$,
	and the initial value of the boundary filters as
	${n_1}(0) = {n_2}(0) = {n_3}(0) = {z_1}(0) = {z_2}(0) = {z_3}(0) = 2$.
	
	A finite difference method (FDM) with 50×2000 grid points uniformly distributed in the [0,1]×[0,10] domain, with spatial step of dx = 0.02 and temporal step of dt = 0.005, is employed to simulate the system dynamics.
	
\begin{figure}[!t] 
		\centering
		%\subfloat[parameter estimation]
		{\includegraphics[width=0.8\linewidth]{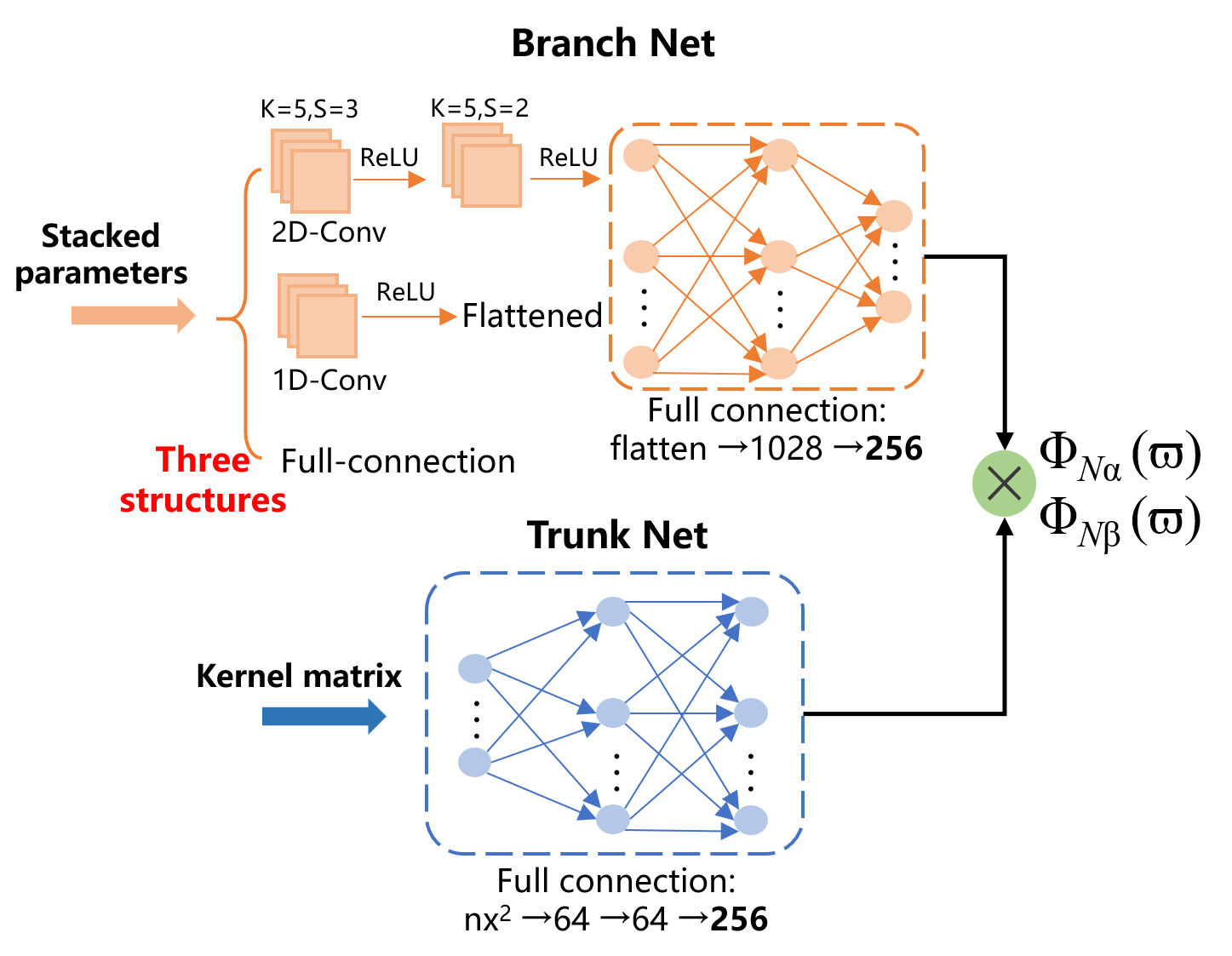}} 
		\caption{The structure of the DeepONet network used for the simulation. Three structures of the branch network are set up for comparison.}
		\label{deeponet structure}
	\end{figure}

	We use DeepONet $\Phi _{\mathcal{N}\alpha}$ and $\Phi _{\mathcal{N}\beta}$ to represent the approximation of solution operators $\Phi _\alpha$ and $\Phi _\beta$, which map the system parameters $m_i$ to the related kernel PDE solutions $\alpha$ and $\beta$. The branch network of our DeepONet comprises a CNN structure, including a two-layer 2-dimensional convolution,  a flattening layer, and two linear layers. Meanwhile, the trunk network is an multi-layer perceptron(MLP) with three fully connected layers. The overall structure of the DeepONet is shown in Fig.\ref{deeponet structure}. 
	
	We utilize DeepONet to approximate the kernel and now outline the training process for the neural operator. To generate an effective and diverse training dataset, we randomly select $N \!=\! 2000$ samples within the parameter boundary, i.e., $2000 \times {m_i} \! \sim \! \text{uniform}(- 4,4),\ i \! = \! 1,...,4$, and then stack the four groups of random parameters as sample inputs. The corresponding gain kernel is calculated by the FDM according to the kernel equation, serving as the sample output.
%	It should be noted that the FDM is highly sensitive to boundary conditions. Randomly selected parameters within the boundaries may lead to significant deviations in the computed kernel solution, leading to numerical dissipation or numerical dispersion effects. Such inaccuracies in the training sets can negatively impact the network training process.
%	We perform a simple filtering operation on the samples obtained to reduce the noise error in the numerical solution, thereby ensuring the effectiveness of the training set. Additionally, the initial boundary values and parameters in the actual simulation system are selected based on stability tests.
	
	We construct 2000 pairs of datasets with a 9:1 ratio divided into training and testing sets for supervised learning of NO.	
	It should be noted that if the range of parameter values is expanded, the number of samples may need to be increased to ensure effective coverage of the target region. This is essential to construct an exhaustive dataset anticipating the
	possible parameter-set ($m_1 \!\! \sim \!\! m_4$) encountered.
	The DeepONet operator is trained by iterative gradient descent to minimize the mean square loss function, utilizing an Adam optimizer and a StepLR learning rate controller. 
	
\begin{figure}[!t] 
		\centering
		\subfloat[The cost function of kernel $\alpha$]{\includegraphics[width=0.49\linewidth]{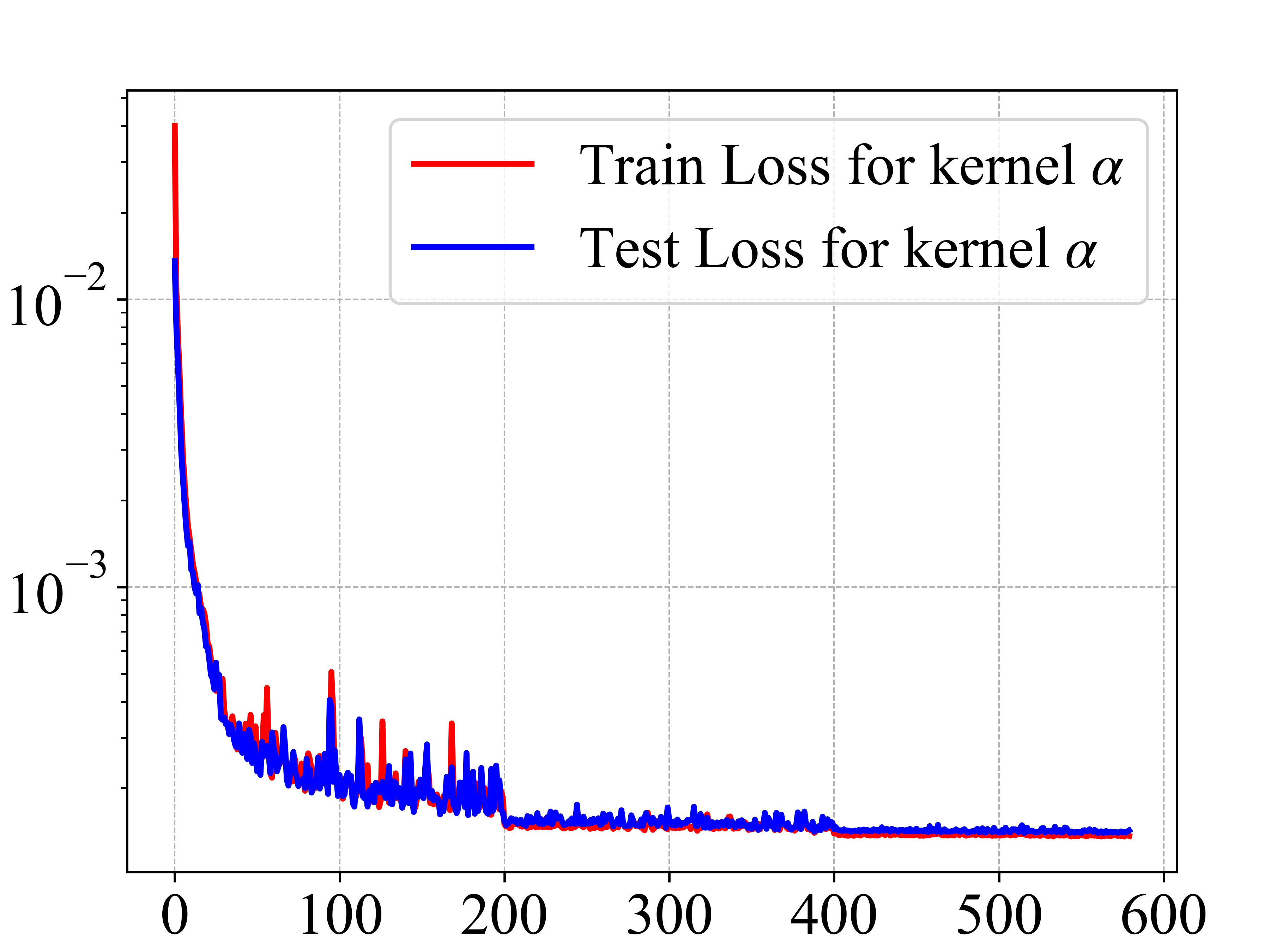}} \label{learning_loss_alpha}
		\subfloat[The cost function of kernel $\beta$]{\includegraphics[width=0.49\linewidth]{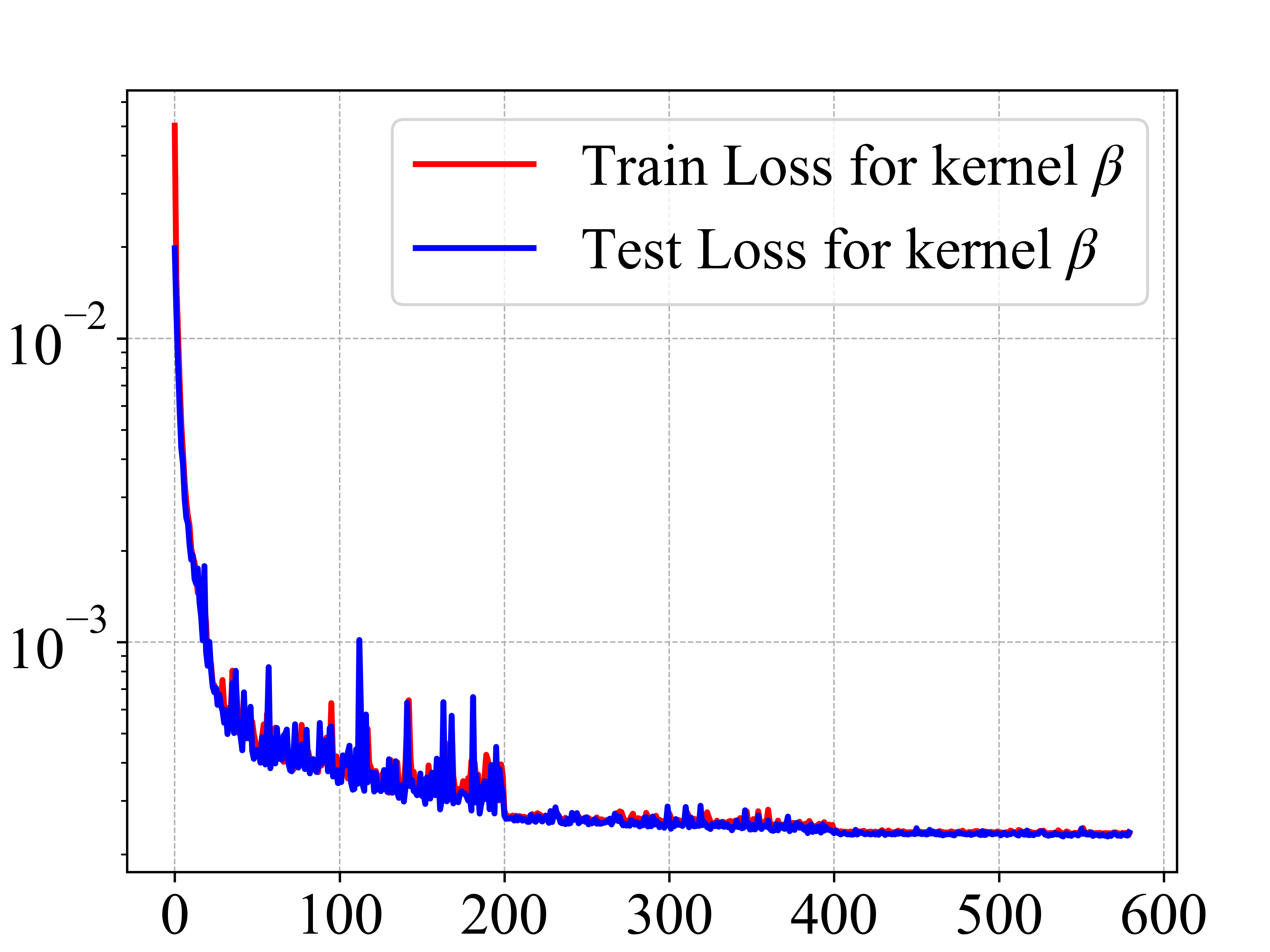}} 
		\label{learning_loss_beta}
		\caption{Training loss of the approximate kernels: the horizontal axis represents the number of iterations of the network, and the vertical axis represents the value of the loss function.}
		\label{train loss}
	\end{figure}

	The training loss of the backstepping kernel for neural network approximation is depicted in Fig.\ref{train loss}. It is observable that after 400 iterations, the training loss of both kernels reduces to the order of $10^{-4}$ and stabilizes. Moreover, the loss curves for both training and test sets are nearly identical, indicating excellent generalization performance of the network model.

\begin{figure}[!t]
		\centering
		\subfloat[Approximate kernel $\hat\alpha$]{\includegraphics[width=0.49\linewidth]{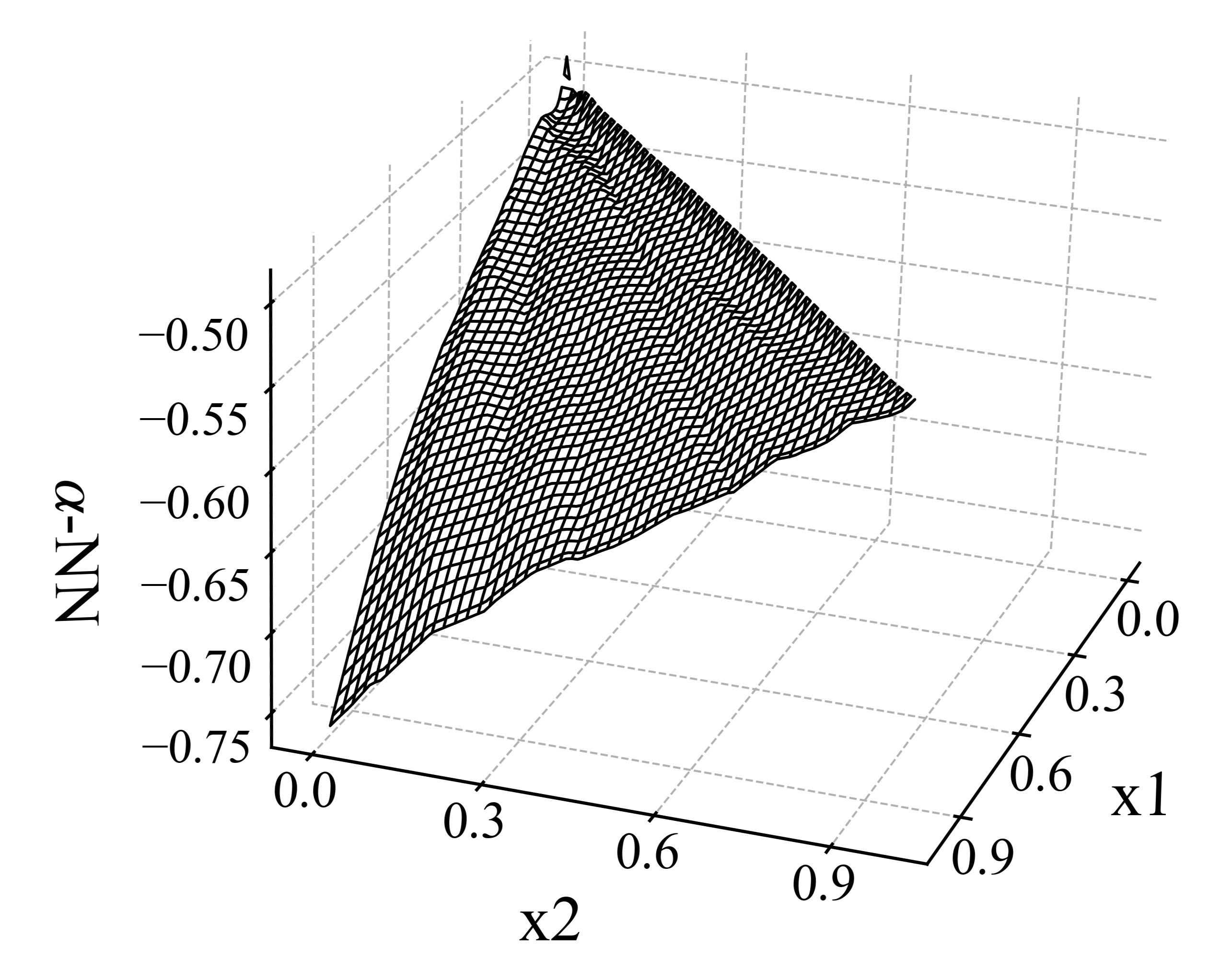}} \label{hat kernel alpha}
		\subfloat[Approximate kernel $\hat\beta$]{\includegraphics[width=0.49\linewidth]{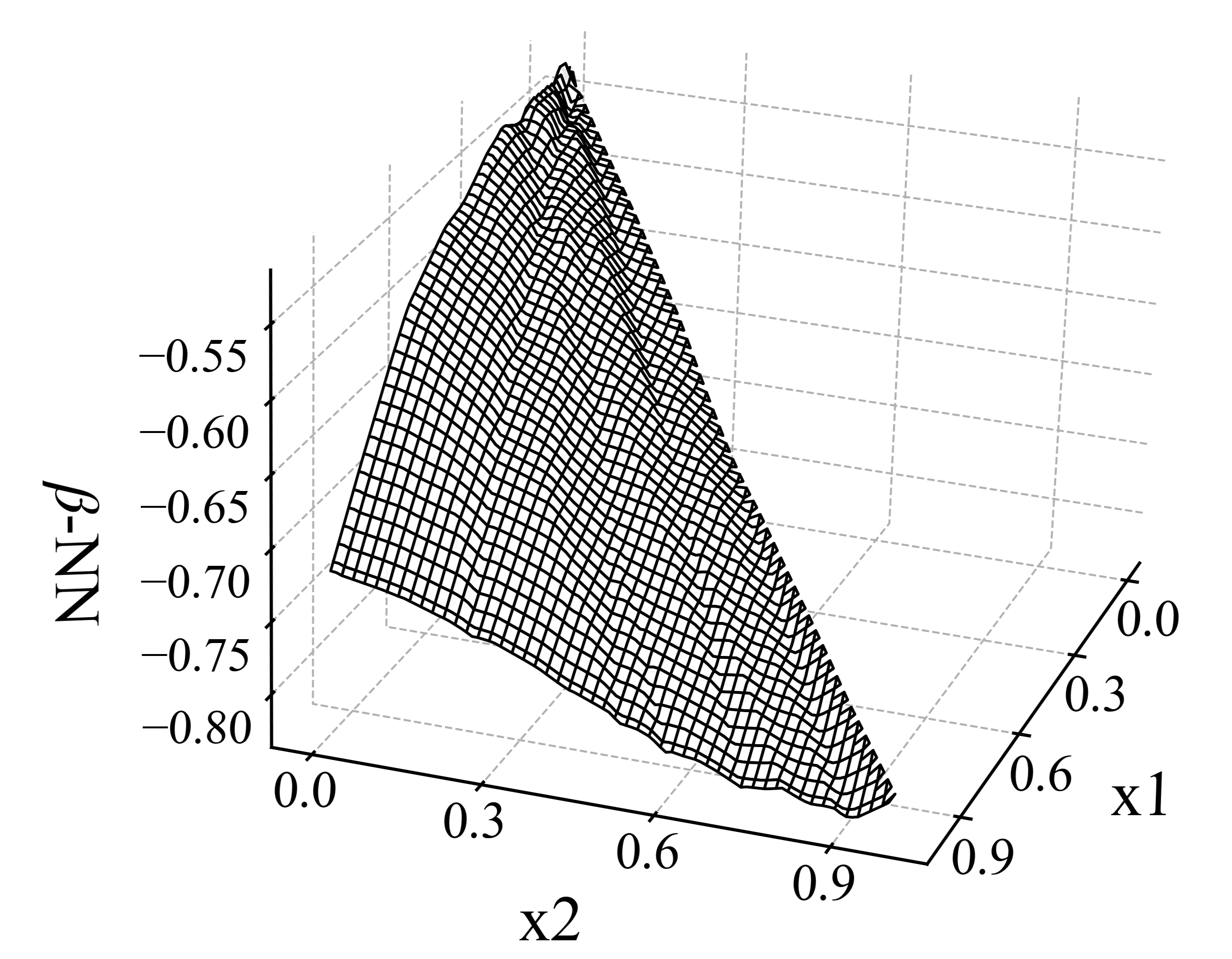}} \label{hat kernel beta}\\
		\subfloat[Numerical kernel $\alpha$]{\includegraphics[width=0.49\linewidth]{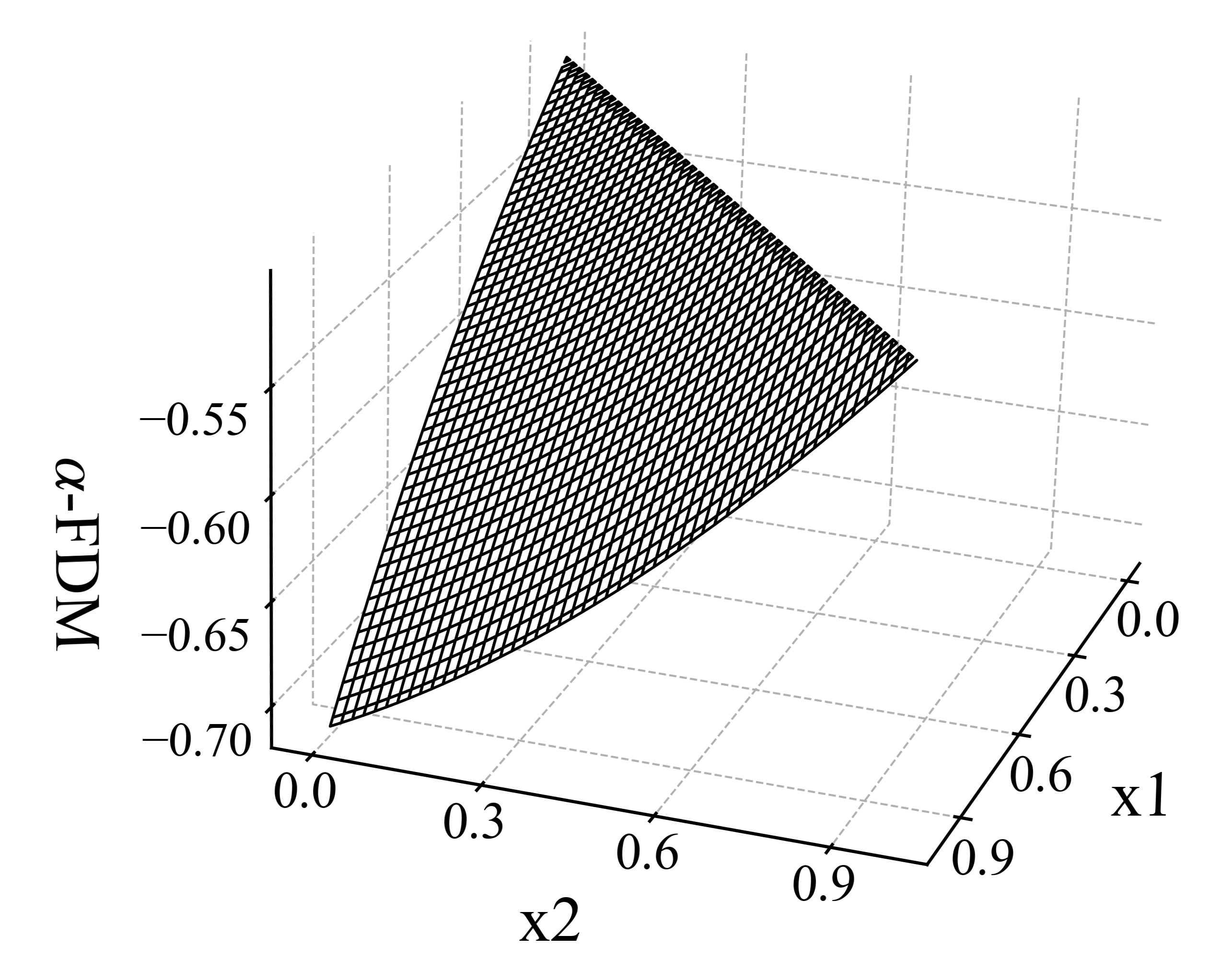}} \label{real kernel alpha} 
		\subfloat[Numerical kernel $\beta$]{\includegraphics[width=0.49\linewidth]{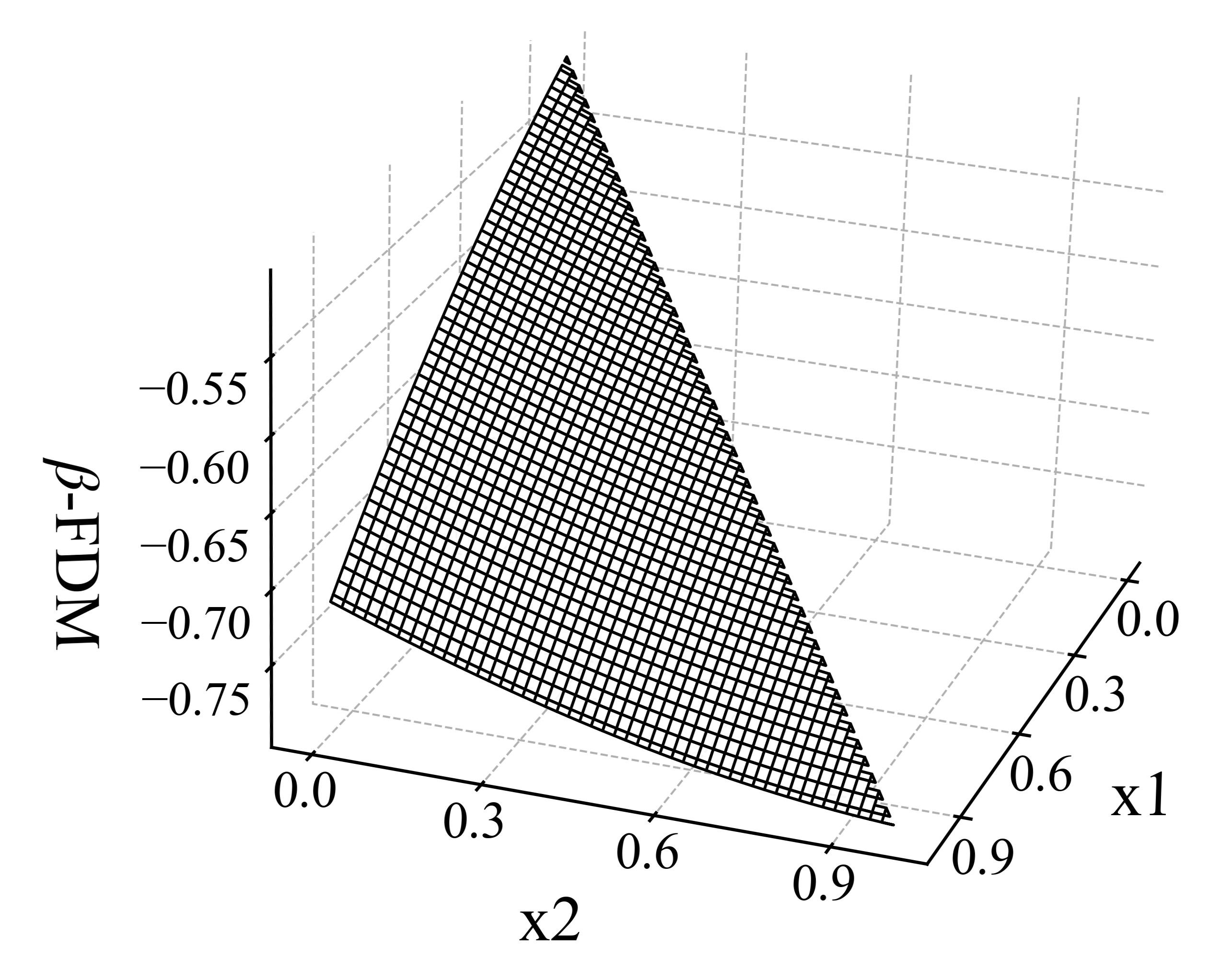}} \label{real kernel beta}\\
		\subfloat[estimation bias $\tilde\alpha$]{\includegraphics[width=0.49\linewidth]{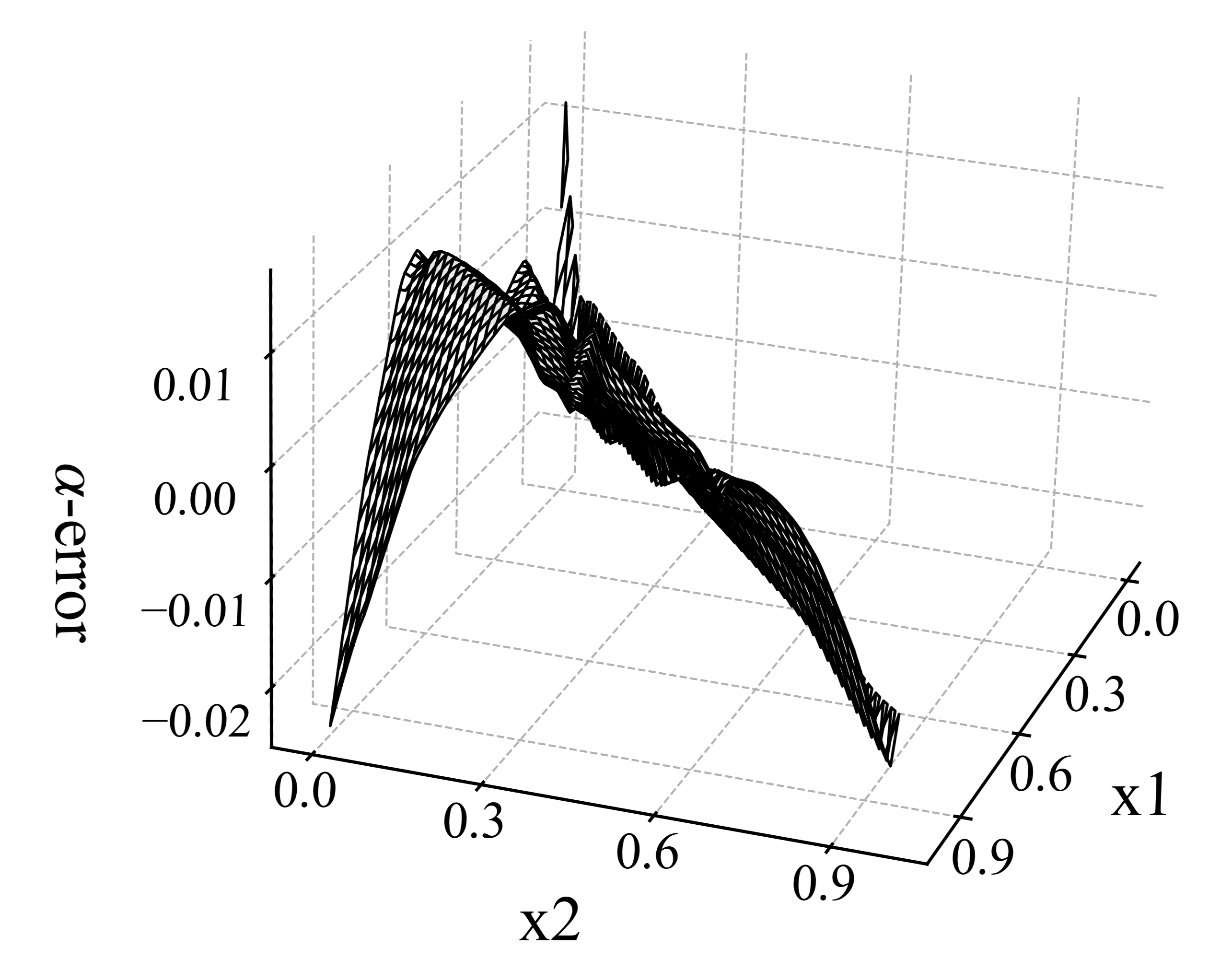}} \label{tilde kernel alphe}
		\subfloat[estimation bias $\tilde\beta$]{\includegraphics[width=0.49\linewidth]{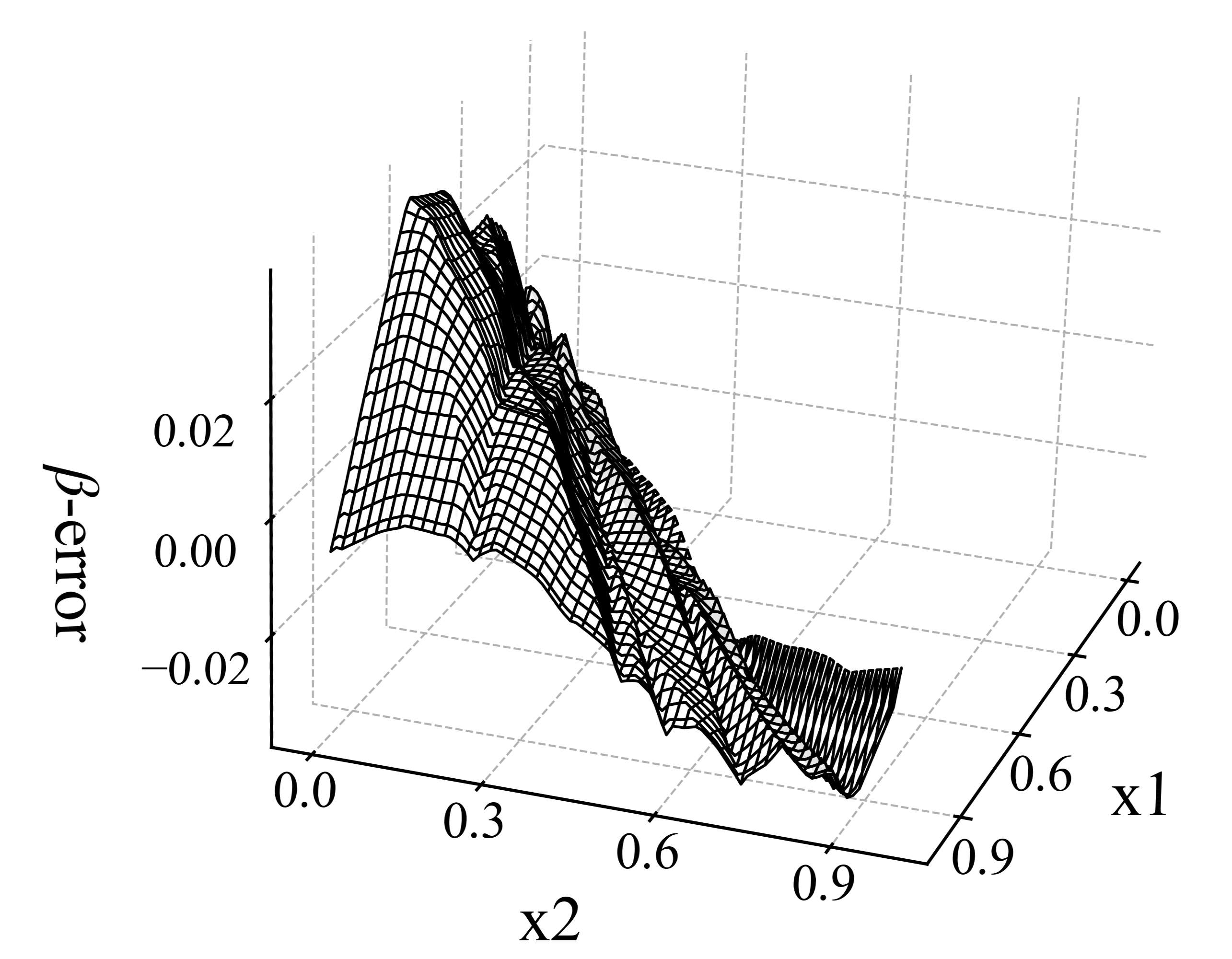}} \label{tilde kernel beta}
		\caption{Backstepping kernels. Exact solution obtained by the finite difference method (Top row). Predicted solution of the trained DeepONet (Middle row). Deviation of estimated value from actual value of corresponding kernel (Buttom row).}
		\label{kernel}
	\end{figure}

	Then, under the above parameter conditions \eqref{actual values}, the approximate backstepping kernels $\hat\alpha$ and $\hat\beta$ obtained by NO approximation are compared with the numerically solved kernels $\alpha$ and $\beta$ in equation \eqref{kernel equation}, as shown in Fig.\ref{kernel}. The relative estimation errors ${\tilde\alpha}/{\alpha}$ and ${\tilde\beta}/{\beta}$ are all within 5$\%$, indicating an acceptable prediction accuracy. 
	
\begin{figure}[!t] 
		\centering
		\subfloat[${\psi _{NO}}(x,t)$]
		{\includegraphics[width=0.49\linewidth]{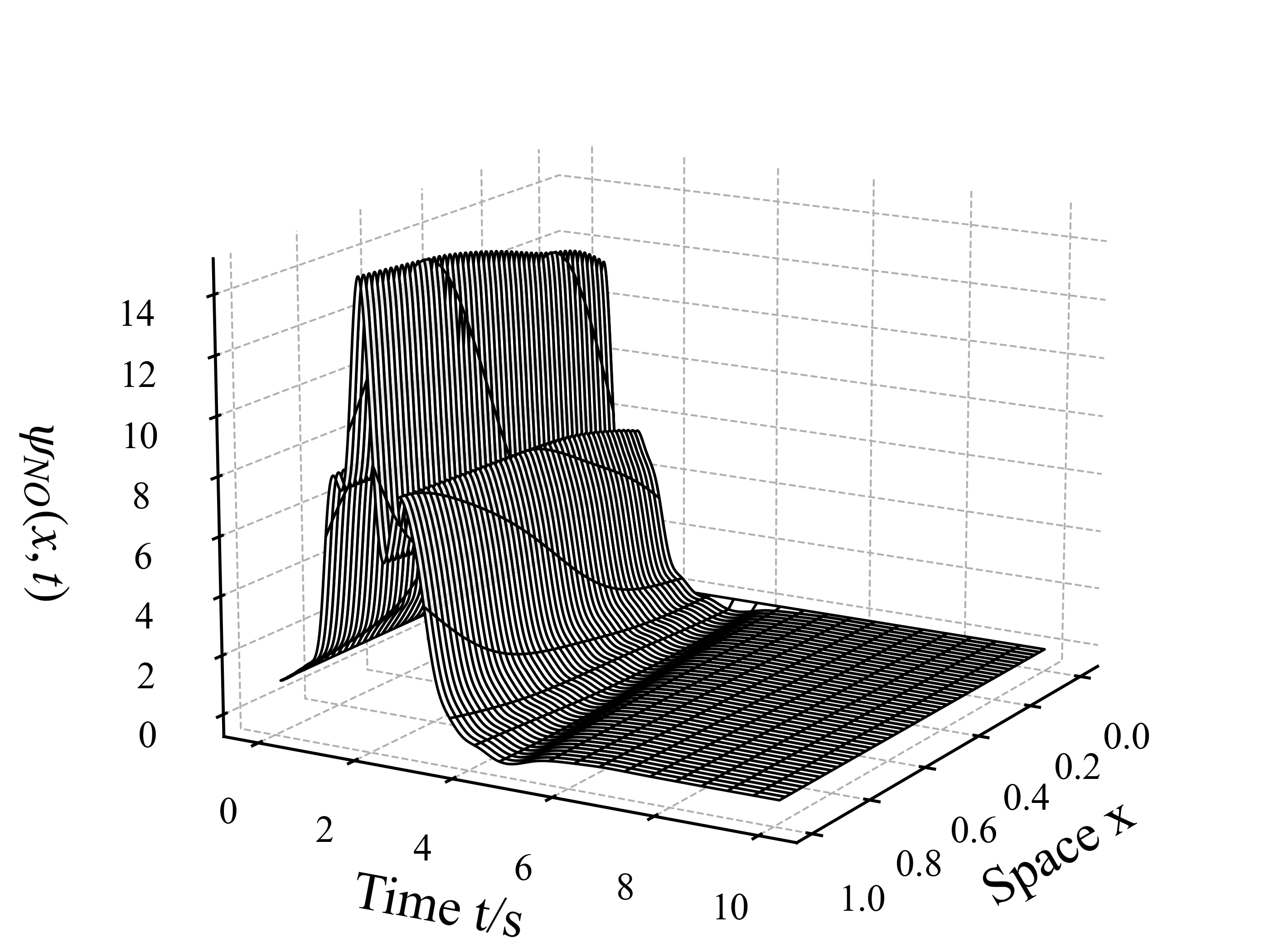}} \label{statu 1}
		\subfloat[${\varphi _{NO}}(x,t)$]
		{\includegraphics[width=0.49\linewidth]{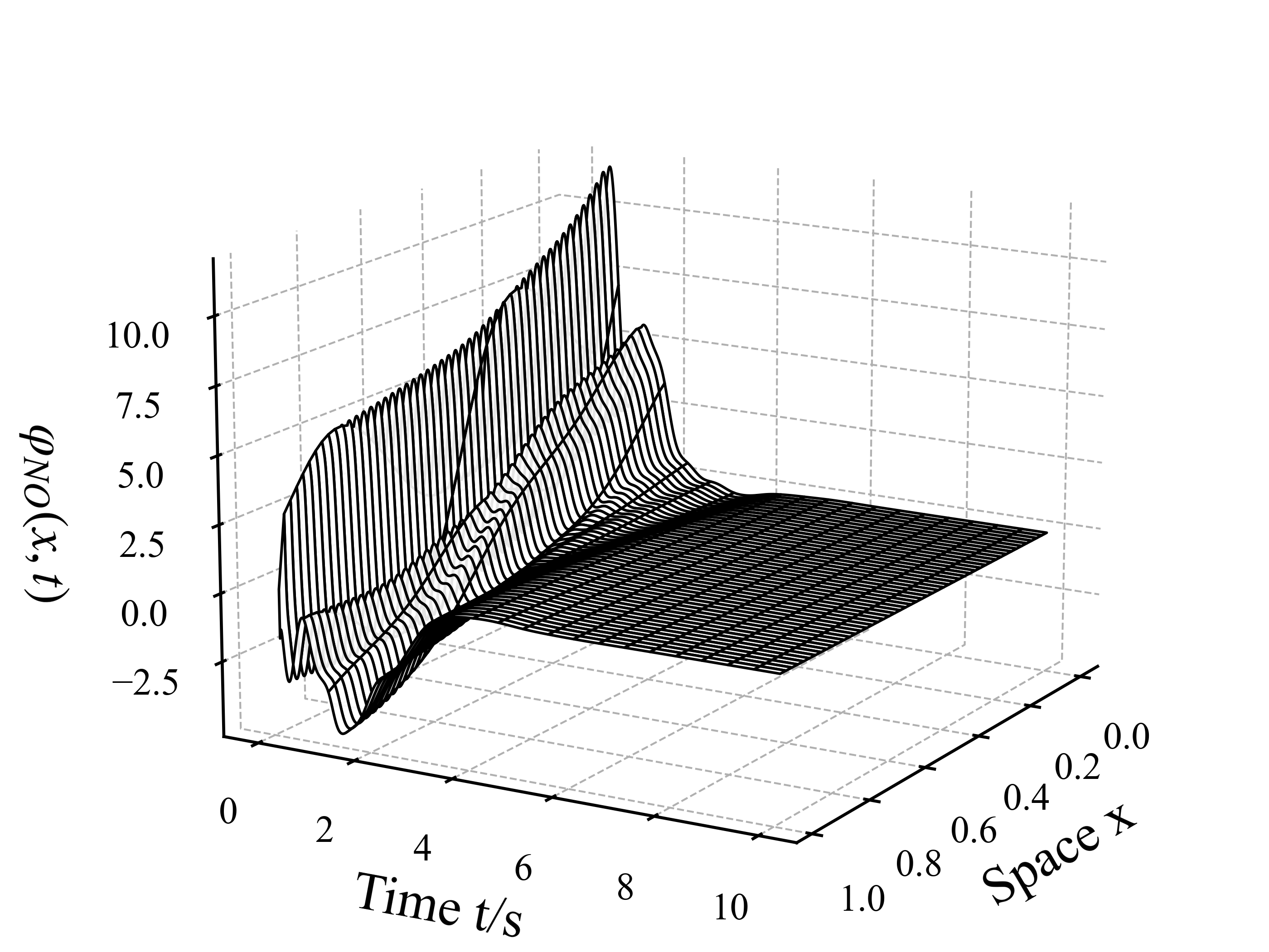}} \label{statu 2}
		\caption{NO-approximate kernels for backstepping design in adaptive stabilization control for a class of 2$\times$2 hyperbolic PDE systems.}
		\label{statu NO}
	\end{figure}

\begin{figure*}[htbp]
		\centering
		\begin{minipage}{0.3\linewidth}
			\centering
			%\subfloat[parameter estimation]
			\includegraphics[width=1\linewidth]{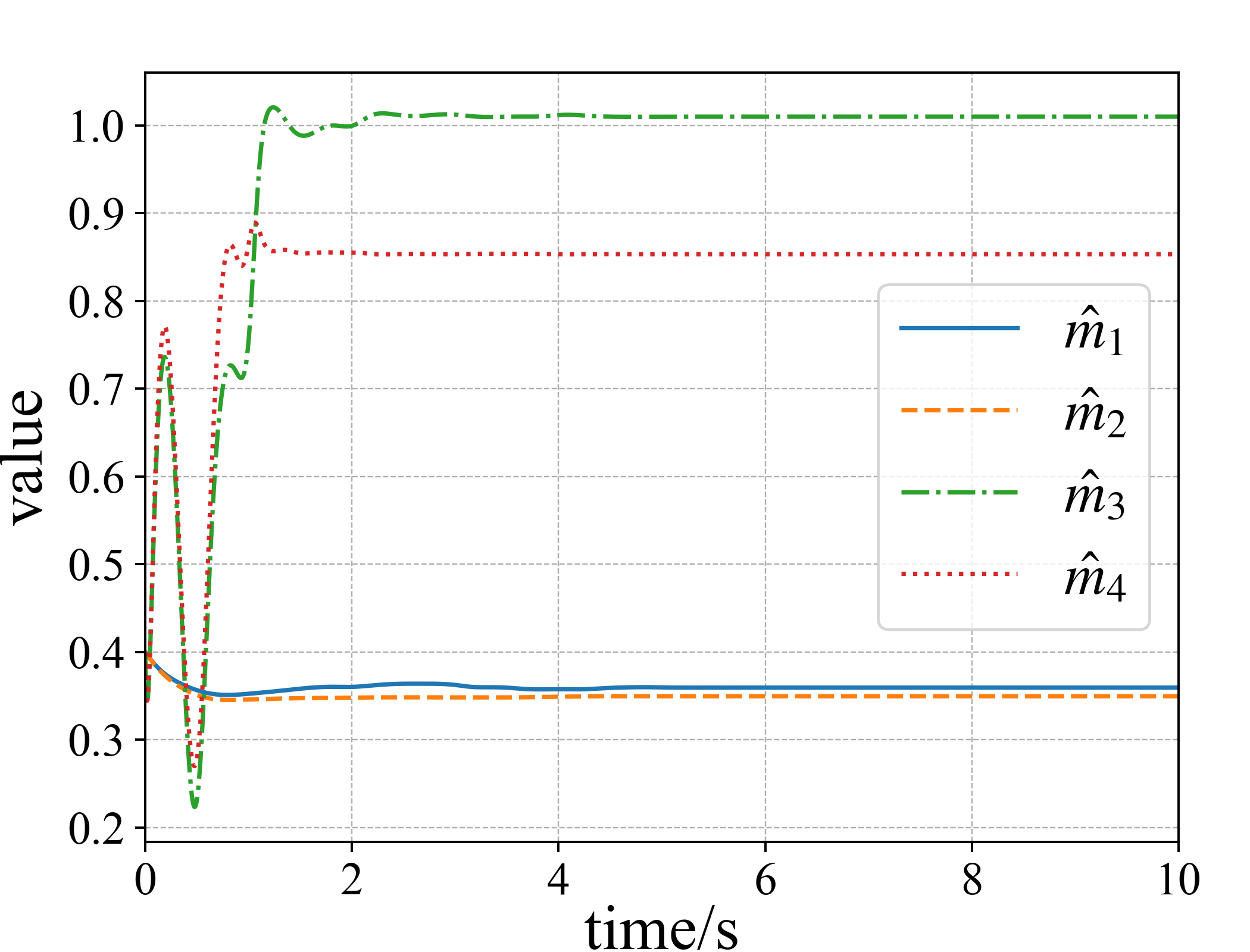}
			\caption{Estimation of unknown parameters of the original system \eqref{original_system} is obtained using the adaptive laws \eqref{adaptive laws}.}
			\label{parameter estimation}
		\end{minipage}
		\quad
		\begin{minipage}{0.3\linewidth}
			\centering
			\includegraphics[width=0.95\linewidth]{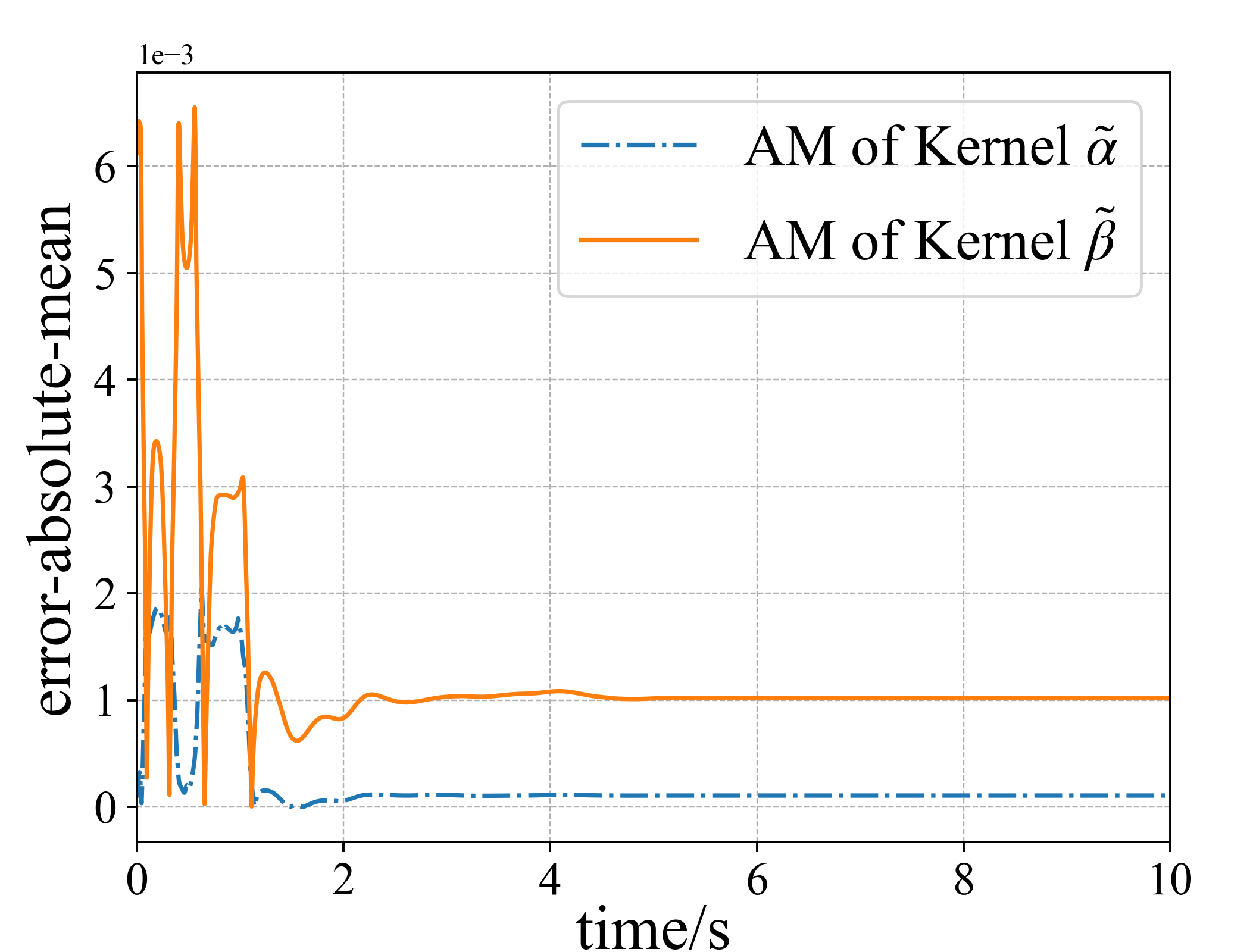}
			\caption{Absolute mean variation of the approximate kernel error matrix for the full process of adaptive control.}
			\label{kernel error variation}
		\end{minipage}
		\quad
		\begin{minipage}{0.3\linewidth}
			\centering
			\includegraphics[width=0.97\linewidth]{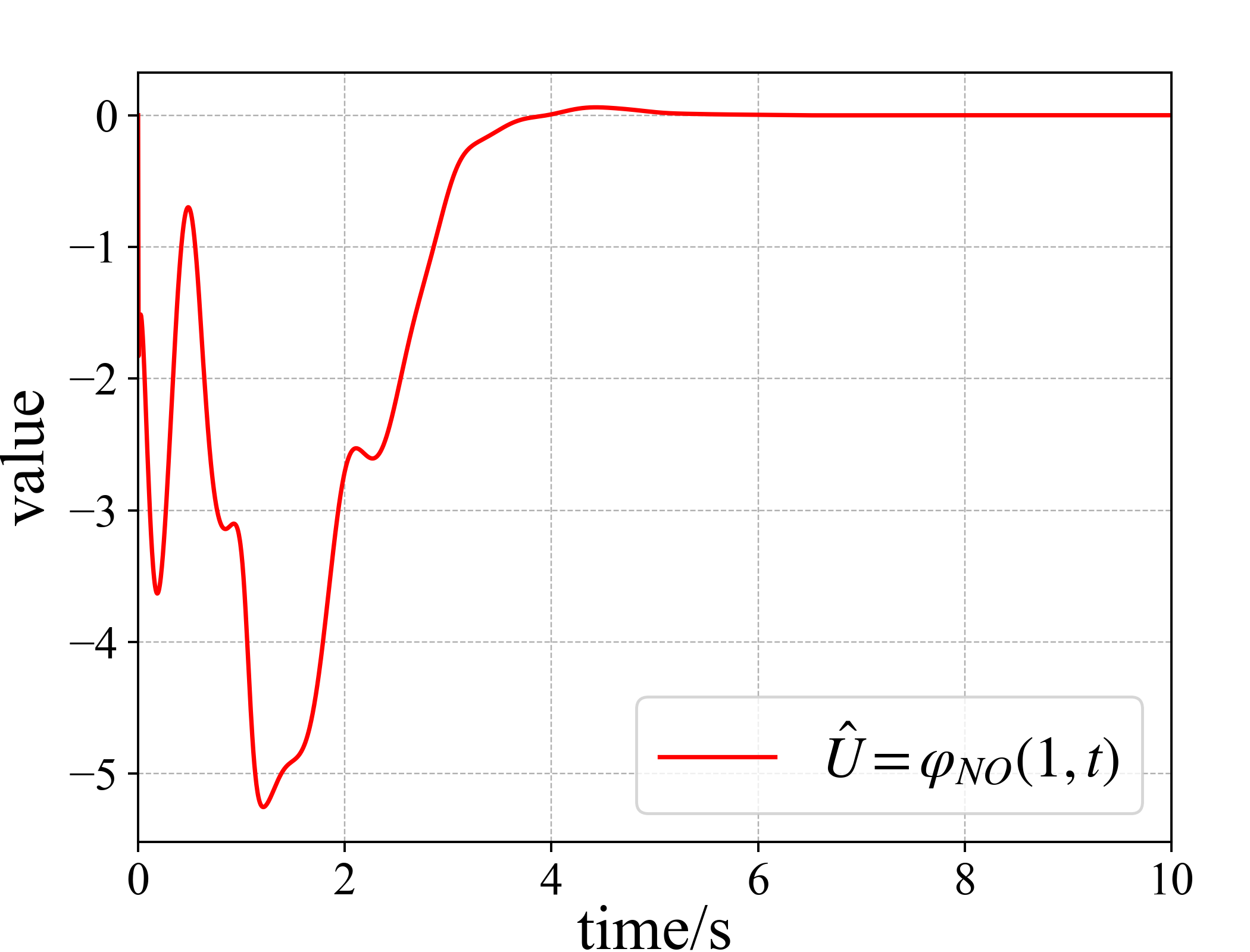}
			\caption{The adaptive controller $\hat{U}$ is obtained by applying approximated kernels from \eqref{approximate controller}.}
			\label{controller}
		\end{minipage}
	\end{figure*}

	In Fig.\ref{statu NO}, the adaptive controller with approximate NO kernel gain can also stabilize the origin system states to zero. 
	However, from 0 to 6 seconds, a noticeable discrepancy exists between the given initial parameters and the actual simulation parameters, as shown in Fig.\ref{parameter estimation}, leading to violent oscillations in the system states. After 6 seconds, the estimated parameters converge to constants and remain stable, resulting in a rapid decay of the system states.
	But there are still uneliminable biases between the estimated and set parameter values, 
%	For instance, the estimated ${\hat m}_2$ has a relative error of $40\%$ after stabilization from the graph, 
	which can be attributed to the limitation of the parameter updating law \eqref{adaptive laws}.
%	The system's state still converges, but not exponentially to zero.
	
	The single instance of the approximate kernel obtained corresponds to a two-dimensional solution space for a set of parameters, as illustrated in Fig.\ref{kernel}. For the complete adaptive control process, the kernel gain needs to be constantly recalculated according to the new parameter estimates. Therefore, the obtained operator kernel corresponds to two spatial dimensions and one temporal dimension. 
	At each moment, we only record the absolute mean of the corresponding kernel error matrix, as shown in Fig.\ref{kernel error variation}. The peak kernels error corresponds to the drastic changes in the system states. 
	In Fig.\ref{controller}, the kernel gains are approximated using neural operator that produces a bounded and convergent control quantity $\hat{U}$, which ensures the stability of the closed-loop system.
	
	Further, based on the initial qualifications mentioned earlier and after completing the effective training of DeepONet operators offline, we subject the system to adaptive control through backstepping design, utilizing the NO approximation kernels. We focusing on analysing the accelerated effect of combining neural operators in adaptive control, as presented in Table \ref{table1}-\ref{table2}.
	
	Table \ref{table1} displays the results of adaptive control simulations for the system under various spatial step size conditions. It records the time taken for both the numerical methods to compute and the network approximation to estimate the kernels over a total of 2000 steps.
	It's clear that as the discrete spatial step size decreases, that is, the sampling and solution precision increases, the computational time for the analytical kernel will increases exponentially. In contrast, the neural operator kernel is largely unaffected by the spatial step size. 
	We also calculated the absolute average error $\frac{1}{10}\int_0^{10} {\int_0^1 {\int_0^x {\left| {\tilde \alpha (x,\xi ,t)} \right| + \left| {\tilde \beta (x,\xi ,t)} \right|d\xi } dx} } dt$ between the neural operator kernels and the numerical kernels at different step sizes. 
	It can be seen that the error is maintained on the order of $10^{-3}$ and increases slightly with decreasing step size(since the scales of input and output data increase significantly with decreasing step size). For the adaptive control task, the effect of such a small error on system stability is negligible. 
	
	Furthermore, DeepONet offers a flexible framework for operator networks, enabling various structures in its sub-networks. Therefore, we adjusted the structure of the branch network accordingly and evaluated the acceleration performance of different NOs. As shown in Fig.\ref{deeponet structure}, we configured the branch network with three structures: 2D convolution by fully connected, 1D convolution by fully connected, and fully connected only (MLP). Then, we recorded the total approximation time in Table \ref{table2}.
	The data indicate that utilizing only the fully-connected layer is more efficient, being approximately 1.5x faster than employing the 2D-convolution network. 
	This suggests that the simpler the hierarchical structure of the network, the more efficient the operator learning, provided that the network can extract the necessary features of the operator.
%	On the other hand, the more complex the structure of the network learning ability is correspondingly stronger, its approximation accuracy is slightly improved but slows down the estimation speed
	
	Taken together, the speed advantage of the neural operator in solving the kernel gain is substantial compared to traditional methods. The NOs scheme improves real-time control and is particularly beneficial for the backstepping adaptive  control problems in PDE systems with uncertain parameters, where real-time solution of the kernel function is essential. 
	
%	This undoubtedly verifies the effectiveness of the proposed intelligent acceleration adaptive control scheme.
	
	\begin{table}[!h]
		\caption{Comparison of computation times for numerical and approximate kernels (2000 steps throughout the control time domain).}
		\centering
		\label{table1}
		\resizebox{0.98\hsize}{!}{
		\begin{tabular}{ccccc} 
			\toprule		
			\multirow{2}{*}{\diagbox{dx}{total-10s}{.}} & \multicolumn{2}{c}{Computation time/$s$} & \multirow{2}{*}{Speedup} & \multirow{2}{*}{Error} \\ 
			\cmidrule{2-3}
			& Numerical solver &NO \\	
			\midrule
			0.02 & 5.96 & 3.19 & 1.87x & $1.21 \!\! \times \!\! {10^{\!-\!3}}$\\
			0.01 & 24.12 & 3.13 & 7.70x & $1.48 \!\! \times \!\! {10^{\!-\!3}}$\\
			0.005 & 106.13 & 3.17 & 33.48x & $1.77 \!\! \times \!\! {10^{\!-\!3}}$ \\
			0.001 & 2717.3 & 3.21 & 846.51x & $2.38 \!\! \times \!\! {10^{\!-\!3}}$\\
			\bottomrule
		\end{tabular}
	}
	\end{table} 
	
	\begin{table}[!h]
		\caption{Three branch network structures, along with their respective approximate kernel computation times.}
		\centering
		\label{table2}
		\resizebox{0.95\hsize}{!}{
		\begin{tabular}{ccccc} 
			\toprule
			\multirow{2}{*}{\diagbox{dx}{time/s}{.}}&\multicolumn{3}{c}{Structure of the branch network}&FDM/MLP-NO\\
			\cmidrule{2-4}
			& Conv2D & Conv1D & MLP &Speedup\\ 
			\midrule
			0.02 & 3.19 & 2.48 & 2.04 & 2.92x\\
			0.01 & 3.13 & 2.46 & 2.03 & 11.88x\\
			0.005 & 3.17 & 2.63 & 2.12 & 50.06x\\
			0.001 & 3.21 & 2.65 & 2.14 & 1269.77x\\
			\bottomrule
		\end{tabular}	
	}
	\end{table}

\section{Conclusion}

In this paper, we integrate NOs with backstepping design for adaptive stabilizing control in a class of 2$\times$2 linear hyperbolic PDE systems containing constant unknown parameters in the domain. The kernel gain is approximated using the DeepONet operator, which accelerates the adaptive control process. 
Theoretical analysis and simulation experiments demonstrate the boundedness and convergence of the controlled system under approximate controller. Simulations also indicate that the neural operator can achieve significant speedup compared to numerical solution methods, making it suitable for implementing real-time adaptive PDE control.

\appendix{
	
% \subsection{Projection operator} \label{appendix A}

\subsection{Proof of target system \eqref{target system2}}\label{appendix B}
	For the equation \eqref{target system2 a}, 
	by substituting equation \eqref{approximate backstepping transformation a} into the dynamic system \eqref{dynamic system a}, then use the error system \eqref{adaptive errors} and replacing $\hat \varphi (x,t)$ with equation \eqref{approximate inverse transformation b}, we obtain
	\begin{equation}
		\begin{footnotesize}
			\setlength{\abovedisplayskip}{2pt}
			\setlength{\belowdisplayskip}{1pt}
			\begin{aligned}
				&{\partial _t}{f_2}(x) + \lambda {\partial _x}{f_2}(x) = {{\hat m}_1}\psi (x) + {{\hat m}_2}\varphi (x) + {n^T}(x)\dot {\hat M} \\
				&= {n^T}(x)\dot {\hat M} \!+\! {{\hat m}_1}[\hat \psi (x) \!+\! \hat e(x)] \!+\! {{\hat m}_2}[\hat \varphi (x) \!+\! \hat \tau (x)]\\
				&={n^T}(x)\dot {\hat M}\!+\!{{\hat m}_1}{f_2}(x) \!+\! {{\hat m}_1}\hat e(x) \!+\! {{\hat m}_2}\hat \tau (x) \!+\! {{\hat m}_2}h_2(x)\\
				&+{{\hat m}_2} \! \int_0^x \! {{{\hat \alpha }_I}(x,\xi )} {f_2}(\xi )d\xi \!+\! {{\hat m}_2} \! \int_0^x \! {{{\hat \beta }_I}(x,\xi )} {h_2}(\xi )d\xi .
			\end{aligned}
		\end{footnotesize}
	\end{equation}
	
	For the equation \eqref{target system2 b}, 
	after applying the shift transformation to equation \eqref{approximate backstepping transformation b}, $\hat \psi (x,t)$ is differentiated with respect to time $t$, followed by the replacement of ${{\partial _t}\hat \psi }$ and ${{\partial _t}\hat \varphi }$ with the dynamic equations \eqref{dynamic system a}-\eqref{dynamic system b}, then integrating the distribution, combined with \eqref{kernel error}, we obtain
	\begin{equation}
		\setlength{\abovedisplayskip}{2pt}
		\setlength{\belowdisplayskip}{1pt}
		\begin{footnotesize}
			\begin{aligned}
					&{\partial _t}\hat \varphi (x)= {\partial _t}{h_2}(x) + \int_0^x {{{\hat \alpha }_t}(x,\xi )\hat \psi (\xi )d\xi }  + \int_0^x {{{\hat \beta }_t}(x,\xi )\hat \varphi (\xi )d\xi } \\
					&+ \int_0^x {\hat \alpha (x,\xi )[{n^T}(\xi )\dot {\hat M}]d\xi }  + \int_0^x {\hat \beta (x,\xi )[{z^T}(\xi )\dot {\hat N}]d\xi } \\
					&+ \!\!\! \int_0^x \!\!\!{\hat \alpha (x,\! \xi )\![{{\hat m}_1}\psi (\xi ) \!\!+\! {{\hat m}_2}\varphi (\xi )]d\xi } \!\!+\!\!\! \int_0^x \!\!\!{\hat \beta (x,\! \xi )\![{{\hat m}_3}\psi (\xi ) \!\!+\! {{\hat m}_4}\varphi (\xi )]d\xi } \\
					&+ \!\lambda \!\!\int_0^x \!\!\!{[{\alpha _\xi }(x,\! \xi ) \!\!-\! {{\tilde \alpha }_\xi }(x,\! \xi )]\hat \psi (\xi )d\xi } \!\!-\! \mu \!\! \int_0^x \!\!\!{[{\beta _\xi }(x,\! \xi ) \!\!-\! {{\tilde \beta }_\xi }(x,\! \xi )]\hat \varphi (\xi )d\xi } \\
					&- \lambda [\alpha (x,x) - \tilde \alpha (x,x)]\hat \psi (x) + \lambda [\alpha (x,0) - \tilde \alpha (x,0)]\hat \psi (0)\\
					&+ \mu [\beta (x,x) - \tilde \beta (x,x)]\hat \varphi (x) - \mu [\beta (x,0) - \tilde \beta (x,0)]\hat \varphi (0).
			\end{aligned}
		\end{footnotesize}
	\end{equation}
	Similarly, differentiating $\hat \psi (x,t)$ with respect to space $x$, we obtain
	\begin{equation}
		\setlength{\abovedisplayskip}{2pt}
		\setlength{\belowdisplayskip}{1pt}
		\begin{footnotesize}
			\begin{aligned}
				&	{\partial _x}\hat \varphi (x) = {\partial _x}h(x) + \int_0^x {[{\alpha _x}(x,\xi ) - {{\tilde \alpha }_x}(x,\xi )]\hat \psi (\xi )d\xi } \\
				&+ \int_0^x {[{\beta _x}(x,\xi ) - {{\tilde \beta }_x}(x,\xi )]\hat \varphi (\xi )d\xi } \\
				&	+ \alpha (x,x)\hat \psi (x) - \tilde \alpha (x,x)\hat \psi (x) + \beta (x,x)\hat \varphi (x) - \tilde \beta (x,x)\hat \varphi (x).
			\end{aligned}
		\end{footnotesize}
	\end{equation}
	Substituting the aforementioned two formulas into \eqref{dynamic system b}, and utilizing the kernel equation set \eqref{kernel equation} to eliminate the middle term, and subsequently combining and sorting, we obtain
	\begin{equation}
		\setlength{\abovedisplayskip}{2pt}
		\setlength{\belowdisplayskip}{1pt}
		\begin{footnotesize}
			\begin{aligned}
				&{{\hat m}_3}\psi (x) + {{\hat m}_4}\varphi (x) + {z^T}(\xi )\dot {\hat N}= {\partial _t}h(x) - \mu {\partial _x}h(x)\\
				&+ \int_0^x {[{{\hat \alpha }_t}(x,\xi ) - \lambda {{\tilde \alpha }_\xi }(x,\xi ) + \mu {{\tilde \alpha }_x}(x,\xi ) + {{\hat m}_4}\alpha (x,\xi )} \\
				&\;\;\;\;\; - {{\hat m}_1}\tilde \alpha (x,\xi ) - {{\hat m}_3}\tilde \beta (x,\xi )]\hat \psi (\xi )d\xi \\
				&+ \!\!\int_0^x \!\!\!{[{{\hat \beta }_t}\!(x,\! \xi ) \!\!+\! \mu {{\tilde \beta }_\xi }\!(x,\! \xi ) \!\!+\! \mu {{\tilde \beta }_x}\!(x,\! \xi ) \!\!-\! {{\hat m}_2}\tilde \alpha \!(x,\! \xi )} \!\!+\! {{\hat m}_4}\hat \beta \!(x,\! \xi )]\hat \varphi \!(\xi )d\xi \\
				&+ \int_0^x {\hat \alpha (x,\xi )} [{{\hat m}_1}\hat e(\xi ) + {{\hat m}_2}\hat \tau (\xi )]d\xi +\int_0^x {\hat \alpha (x,\xi )} [{n^T}(\xi )\dot {\hat M}]d\xi\\
				&+ \int_0^x {\hat \beta (x,\xi )} [{{\hat m}_3}\hat e(\xi ) + {{\hat m}_4}\hat \tau (\xi )]d\xi +\int_0^x {\hat \beta (x,\xi )} [{z^T}(\xi )\dot {\hat N}]d\xi\\
				&+ (\lambda  \!+\! \mu )\tilde \alpha (x,x)\hat \psi (x) \!+\! \hat \psi (x){{\hat m}_3} \!+\! \lambda \hat \alpha (x,0)\hat \psi (0) \!-\! \mu \hat \beta (x,0)\hat \varphi (0).
			\end{aligned}
	     \end{footnotesize}
	\end{equation}
	Form \eqref{approximate inverse transformation}, \eqref{adaptive errors}, we get formula \eqref{target system2 b}.
	
	For the equation \eqref{target system2 c}, with \eqref{approximate backstepping transformation} and \eqref{dynamic system c}, we obtain:
	${f_2}(0) = \hat \psi (0) = q\varphi (0) = q\hat \varphi (0) + q\hat \tau (0)= q{h_2}(0) + q\hat \tau (0). $

	For the equation \eqref{target system2 d}, from \eqref{approximate backstepping transformation b}, \eqref{dynamic system d} and \eqref{approximate controller}, we obtain: 
	$ {h_2}(1) = \hat \varphi (1) - \hat U = U - \hat U. $
	
	\subsection{Proof of Theorem \ref{converge}}\label{appendix C}
	Firstly, we introduce some important properties to facilitate analysis.
	According to Cauchy–Schwarz's inequality, following inequalities which hold provided $\delta  \ge 1$
	\begin{subequations} \label{cauchy}
		\setlength{\abovedisplayskip}{3pt}
		\setlength{\belowdisplayskip}{3pt}
		\begin{small}
			\begin{align}
				&{I_\delta }\left[ {u\psi v} \right] \le \left\| u \right\|k \left\| \psi  \right\|_\delta ^2 + \frac{{\left\| u \right\|}}{k}\left\| v \right\|_\delta ^2 \label{cauchy a}, \\
				&{I_\delta }\left[\! {\psi (x) \! \int_0^x \! {u(x, \! \xi )v(\xi )d\xi } } \! \right] \! \le \! \left\| u \right\| \! k \!  \left\| \psi  \right\|_\delta ^2 \!+\! \frac{{\left\| u \right\|}}{k} \! {e^\delta }{\left\| v \right\|^2} \label{cauchy b}, \\
				&{I_{ - \delta }}\left[ {\psi (x)\int_0^x {u(x,\xi )v(\xi )d\xi } } \right] \le \bar uk  \left\| \psi  \right\|_{ - \delta }^2 + \frac{{\bar u}}{k}{\left\| v \right\|^2}. \label{cauchy c}
			\end{align}
		\end{small}
	\end{subequations}
	From \eqref{approximate kernel boundary} and \eqref{approximate backstepping transformation}, since the kernels
	$(\hat \alpha ,\hat \beta )$ and $({{\hat \alpha }_I},{{\hat \beta }_I})$ are bounded for every $t$, for some signals
	${u_1}(x),{u_2}(x),{u_3}(x)$ defined on ${T = \{ 0 \le x \le 1,t \ge 0\}}$ and satisfying
	${u_3}(x) = {\Gamma _2}[{u_1},{u_2}](x)$ and ${u_2}(x) = \Gamma _2^{ - 1}[{u_1},{u_3}](x)$, there exist positive constants $M1,M2,N1,N2$ such that
	\begin{subequations} 
		\setlength{\abovedisplayskip}{2pt}
		\setlength{\belowdisplayskip}{2pt}
		\begin{small}
			\begin{align} 
				&{u_3}= {\Gamma _2}[{u_1},{u_2}] \;\; \to \left\| {{u_3}} \right\| \le {M_1}\left\| {{u_1}} \right\| + {M_2}\left\| {{u_2}} \right\| ,\\
				&{u_2}= {\Gamma_2^{-1} }[{u_1},{u_3}] \to\left\| {{u_1}} \right\| \le {N_1}\left\| {{u_1}} \right\| + {N_2}\left\| {{u_3}} \right\|.
			\end{align}
		\end{small}
	\end{subequations}
	Then from \eqref{adaptive errors} and \eqref{approximate backstepping transformation}, we obtain
	\begin{subequations}\label{real status boundary}
		\setlength{\abovedisplayskip}{3pt}
		\setlength{\belowdisplayskip}{3pt}
		\begin{small}
			\begin{align}
				&\left\| \psi  \right\| = \left\| {\hat \psi  + \hat e} \right\| \le \left\| {\hat \psi } \right\| + \left\| {\hat e} \right\| = \left\| {{f_2}} \right\| + \left\| {\hat e} \right\| \label{real status boundary a}, \\
				&\begin{aligned} \label{real status boundary b}
					\left\| \varphi  \right\| &= \left\| {\hat \varphi  + \hat \tau } \right\| \le \left\| {\hat \varphi } \right\| + \left\| {\hat \tau } \right\| = \left\| {\Gamma _2^{ - 1}[{f_2},{h_2}]} \right\| + \left\| {\hat \tau } \right\| \\
					&\le {N_1}\left\| {{f_2}} \right\| + {N_2}\left\| {{h_2}} \right\| + \left\| {\hat \tau } \right\| .
				\end{aligned}
			\end{align}
		\end{small}
	\end{subequations}
	By equations \eqref{non-adaptive states}, \eqref{non-adaptive errors}, \eqref{adaptive errors} and \eqref{adaptive states}, the estimated errors satisfy
	\begin{subequations}\label{estimated error boundary}
		\setlength{\abovedisplayskip}{3pt}
		\setlength{\belowdisplayskip}{3pt}
		\begin{small}
		\begin{align}
			&\left\| {\hat e} \right\| = \left\| {e + \bar \psi  - \hat \psi } \right\| = \left\| {e + {n^T}\tilde M} \right\| \le \left\| e \right\| + \left| {\tilde M} \right|\left\| n \right\| \label{estimated error boundary a}, \\
			&\left\| {\hat \tau } \right\| = \left\| {\tau  + \bar \varphi  - \hat \varphi } \right\| = \left\| {\tau  + {z^T}\tilde N} \right\| \le \left\| \tau  \right\| + \left| {\tilde N} \right|\left\| z \right\|. \label{estimated error boundary b} 
		\end{align}
		\end{small}
	\end{subequations}
    After clarifying these properties, we proceed to analyze the Lyapunov sub-function candidates defined in Eqs \eqref{Lyapunov sub-function}.
    
    \subsubsection{${V_1}=\left\| {{n_1}} \right\|_{ - a}^2$ and ${V_2} = \left\| {{n_2}} \right\|_{ - a}^2$} 
    From filter \eqref{filter1 a}, equation \eqref{zn} and \eqref{MN}, employing the properties \eqref{property 1}, \eqref{real status boundary a} and \eqref{estimated error boundary a}, we have
    \vspace{-2mm}
    \begin{equation}\label{V1}
    	\begin{footnotesize}
    		\begin{aligned}
    			&{{\dot V}_1} = 2\int_0^1 {{e^{ - ax}}n_1^T{\partial _t}{n_1}} dx\\
    			&=  - 2\lambda {I_{ - a}}[n_1^T{\partial _x}{n_1}] + 2{I_{ - \delta }}[{n_1}\psi ]\\
    			&\le  - \lambda (n_1^2(1){e^{ - a}} \!-\! n_1^2(0) + a\left\| {{n_1}} \right\|_{ - a}^2) + \left\| {{n_1}} \right\|_{ - a}^2 + \left\| \psi  \right\|_{ - a}^2\\
    			&\le (1 - \lambda a)\left\| {{n_1}} \right\|_{ - a}^2 + 2{\left\| {{f_2}} \right\|^2} + 4{\left\| e \right\|^2} + 4{\left| {\tilde M} \right|^2}{\left\| n \right\|^2}.
    		\end{aligned}
    	\end{footnotesize}
    	\vspace{-2mm}
    \end{equation}
    
    By filter \eqref{filter1 b}, equation \eqref{zn} and \eqref{MN}, employing the properties \eqref{property 1}, \eqref{real status boundary b} and \eqref{estimated error boundary b}, we obtain
    \vspace{-2mm}
    \begin{equation}\label{V2}
    	\begin{footnotesize}
    		\begin{aligned}
    			&{{\dot V}_2} = 2\int_0^1 {{e^{ - ax}}n_2^T{\partial _t}{n_2}} dx\\
    			&=  - 2\lambda {I_{ - a}}[n_2^T{\partial _x}{n_2}] + 2{I_{ - a}}[{n_2}\varphi ]\\
    			&\le \! (1 \!-\! \lambda a) \! \left\| {{n_2}} \right\|_{ - a}^2 \!\!+\! 4N_1^2{\left\| {{f_2}} \right\|^2} \!\!+\! 4N_2^2{\left\| {{h_2}} \right\|^2} \!\!+\! 4{\left\| \tau  \right\|^2} \!\!+\! 4{\left| {\tilde N} \right|^2} \!\! {\left\| z \right\|^2}.
    		\end{aligned}
    	\end{footnotesize}
    	\vspace{-2mm}
    \end{equation}
    
    Since the forms  of $V_1$ and $V_2$ are essentially the same, we combine the stability analysis accordingly, and utilizing property \eqref{property 2} to obtain:
    \vspace{-2mm}
    \begin{equation}\label{V1 V2}
    	\begin{footnotesize}
    		\begin{aligned}
    			&{{\dot V}_1} + {{\dot V}_2} \le (1 - \lambda a)\left\| n \right\|_{ - a}^2 + 4{\left| {\tilde M} \right|^2}{\left\| n \right\|^2} + 4{\left| {\tilde N} \right|^2}{\left\| z \right\|^2}\\
    			&+ (2 + 4N_1^2){\left\| {{f_2}} \right\|^2} + 4N_2^2{\left\| {{h_2}} \right\|^2} + 4{\left\| e \right\|^2} + 4{\left\| \tau  \right\|^2}\\
    			&\le (1 - \lambda a)[{V_1} + {V_2}] + {l_1}[{V_1} + {V_2}] + {l_2}[{V_3} + {V_4}]\\
    			&+ {c_1}{e^a }{V_7} + {c_2}{V_8} + 4{e^a }{V_9} + 4{V_{10}},
    		\end{aligned}
    	\end{footnotesize}
    	\vspace{-2mm}
    \end{equation}
    where 
    \begin{footnotesize}
    	$ {l_1} = 4{e^a}{\left| {\tilde M} \right|^2}$, ${l_2} = 4{\left| {\tilde N} \right|^2}$
    \end{footnotesize}
    are integrable functions and 
    \begin{footnotesize}
    	${c_1} = 4N_1^2 + 2 $, ${c_2} = 4N_2^2$
    \end{footnotesize}
    are two positive constants.
    
    \subsubsection{${V_3} = \left\| {{n_3}} \right\|_{ - a}^2$}
    From filter \eqref{filter1 c} and equation \eqref{MN}, with the properties \eqref{property 1}, \eqref{real status boundary b} and \eqref{estimated error boundary a}, we have
    \vspace{-1mm}
	\begin{equation}\label{V3}
		\begin{footnotesize}
			\begin{aligned}
				&{{\dot V}_3} =  - 2\lambda {I_{ - a}}[n_3^T{\partial _x}{n_3}] =  - \lambda (n_3^2(1){e^{ - a}} \!-\! n_3^2(0) \!+\! a\left\| {{n_3}} \right\|_{ - a}^2)\\
				&\le  - \lambda a\left\| {{n_3}} \right\|_{ - a}^2 \!+\! \lambda {\varphi ^2}(0) \! \le \!  - \lambda a\left\| {{n_3}} \right\|_{ - a}^2 \!+\! 2\lambda {h_2}^2(0) \!+\! 2\lambda {{\hat \tau }^2}(0)\\
				& \le  - \lambda a{V_3} + {l_3}{\left| {\sigma (0)} \right|^2} + 2\lambda {h^2}(0) + 4\lambda {\tau ^2}(0),
			\end{aligned}
		\end{footnotesize}
		\vspace{-1mm}
	\end{equation}
	where 
	\begin{footnotesize}
		${l_3} = 4\lambda {\left| {\tilde N} \right|^2}$
	\end{footnotesize}
	denotes functions that can be integrated.\\
	
	\subsubsection{${V_4} = \left\| {{z_1}} \right\|_b^2$ and ${V_5} = \left\| {{z_2}} \right\|_b^2$} 
	From filter \eqref{filter2 a}, equation \eqref{zn} and \eqref{MN}, employing the properties \eqref{property 1}, \eqref{real status boundary a} and \eqref{estimated error boundary a}, we have
	\vspace{-2mm}
	\begin{equation}\label{V4}
		\begin{footnotesize}
			\begin{aligned}
				&{{\dot V}_4} = 2\int_0^1 {{e^{bx}}z_1^T{\partial _t}{z_1}} dx = 2\mu {I_b}[z_1^T{\partial _x}{z_1}] + 2{I_b}[{z_1}\psi ]\\
				&\le \mu (z_1^2(1){e^b} - z_1^2(0) - b\left\| {{z_1}} \right\|_b^2) + \left\| {{z_1}} \right\|_b^2 + \left\| \psi  \right\|_b^2\\
				&\le \!(1\! -\! \mu b)\!\! \left\| {{z_1}} \right\|_b^2 \!\! -\! \mu z_1^2(0) \!+\! 2{e^b}\!{\left\| {{f_2}} \right\|^2} \!\!+\! 4{e^b} \! {\left\| e \right\|^2} \!\!+\! 4{e^b} \! {\left| {\tilde M} \right|^2} \! \!{\left\|  n \right\|^2}.
			\end{aligned}
		\end{footnotesize}
		\vspace{-2mm}
	\end{equation}
	
	By filter \eqref{filter2 b}, equation \eqref{zn} and \eqref{MN}, employing the properties \eqref{property 1}, \eqref{real status boundary b} and \eqref{estimated error boundary b}, we obtain
	\begin{equation}\label{V5}
		\setlength{\abovedisplayskip}{2pt}
		\setlength{\belowdisplayskip}{1pt}
		\begin{footnotesize}
			\begin{aligned}
				&{{\dot V}_5} = 2\mu {I_b}[z_2^T{\partial _x}{z_2}] + 2{I_b}[{z_2}\varphi ]\\
				&\le (1 - \mu b)\left\| {{z_2}} \right\|_b ^2 - \mu z_2^2(0) + 4{e^b}N_1^2{\left\| {{f_2}} \right\|^2}\\
				&\;\; + 4{e^\alpha }N_2^2{\left\| {{h_2}} \right\|^2} + 4{e^b}{\left\| \tau  \right\|^2} + 4{e^b}{\left| {\tilde N} \right|^2}{\left\| z \right\|^2}.
			\end{aligned}
		\end{footnotesize}
	\end{equation}
	
	Since the forms  of $V_4$ and $V_5$ are essentially the same, we combine the stability analysis accordingly, and utilizing property \eqref{property 2} to obtain:
	\begin{equation}\label{V4 V5}
		\setlength{\abovedisplayskip}{2pt}
		\setlength{\belowdisplayskip}{1pt}
		\begin{footnotesize}
			\begin{aligned}
				&{{\dot V}_4} + {{\dot V}_5} \le (1 - \mu b)\left\| z \right\|_b^2 + 4{\left| {\tilde M} \right|^2}{e^b}{\left\| n \right\|^2} + 4{\left| {\tilde N} \right|^2}{e^b}{\left\| z \right\|^2} \\
				&+ \!(4N_1^2{e^b} \!\!+\!\! 2{e^b}){\left\| {{f_2}} \right\|^2} \!\!+\!\! 4N_2^2{e^b }{\left\| {{h_2}} \right\|^2} \!\!+\!\! 4{e^b}{\left\| e \right\|^2} \!\!+\!\! 4{e^b}{\left\| \tau  \right\|^2} \!\!-\!\! \mu {z^2}(0)\\
				&\le (1 - \mu b)[{V_4} + {V_5}] + {l_4}[{V_1} + {V_2}]+ {l_5}[{V_4} + {V_5}] \\
				&+ {c_3}{e^{b + a}}{V_7} + {c_4}{e^{b}}{V_8} + 4{e^{b + a}}{V_9} + 4{e^b}{V_{10}} - \mu {z^2}(0),
			\end{aligned}
		\end{footnotesize}
	\end{equation}
	where 
	\begin{footnotesize}
		$ {l_4} = 4{\left| {\tilde M} \right|^2}{e^{b + a}}$, ${l_5} = 4{\left| {\tilde N} \right|^2}{e^b}$
	\end{footnotesize}
	are integrable functions and 
	\begin{footnotesize}
		$ {c_3} = 2 + 4N_1^2$, ${c_4} = 4N_2^2$
	\end{footnotesize}
	are positive constants.
	
	\subsubsection{${V_6} = \left\| {{z_3}} \right\|_b^2$}
	Using \eqref{property 1} and filter \eqref{filter2 c}, we can get
	\begin{equation}
		\setlength{\abovedisplayskip}{2pt}
		\setlength{\belowdisplayskip}{1pt}
		\begin{footnotesize}
			\begin{aligned}
				{{\dot V}_6} &= 2\mu {I_b}[z_3^T{\partial _x}{z_3}] = \mu (z_3^2(1){e^b } - z_3^2(0) - b\left\| {{z_3}} \right\|_b^2)\\
				&=  - \mu b\left\| {{z _3}} \right\|_b^2 - \mu z_3^2(0) + \mu {e^b}{{\hat U}^2}.
			\end{aligned}
		\end{footnotesize}
	\end{equation}
	Considering formula \eqref{approximate controller}, \eqref{approximate kernel boundary} and \eqref{real status boundary b}, we analyze ${\hat U}^2$
	to obtain:
	\vspace{-2mm}
	\begin{equation}\label{hat U}
		\begin{footnotesize}
			\begin{aligned}
				&{{\hat U}^2} = {\left( {\int_0^1 {\hat \alpha (1,\xi )\hat \psi (\xi )d\xi }  + \int_0^1 {\hat \beta (1,\xi )\hat \varphi (\xi )d\xi } } \right)^2}\\
				&\le 2\int_0^1 {{{\hat \alpha }^2}(1,\xi )} {{\hat \psi }^2}(\xi )d\xi  + 2\int_0^1 {{{\hat \beta }^2}(1,\xi )} {{\hat \varphi }^2}(\xi )d\xi \\
				&\le \! 2{{\bar {\hat \alpha} }^2}{\left\| {\hat \psi } \right\|^2} \!\!\!+\!\! 2{{\bar {\hat \beta} }^2}{\left\| {\hat \varphi } \right\|^2} \!\le \! 2{{\bar {\hat \alpha} }^2}{\left\| {{f_2}} \right\|^2} \!\!+\! 4{{\bar {\hat \beta} }^2} \! N_1^2{\left\| {{f_2}} \right\|^2} \!\!+\! 4{{\bar {\hat \beta} }^2}\! N_2^2{\left\| {{h_2}} \right\|^2}
			\end{aligned}
		\end{footnotesize}
		\vspace{-2mm}
	\end{equation}
	Hence, by utilizing property \eqref{property 2}, we can obtain
	\vspace{-2mm}
	\begin{equation}\label{V6}
		\begin{footnotesize}
			\begin{aligned}
				&{{\dot V}_6} =  - \mu b\left\| {{z _3}} \right\|_b^2 - \mu z_3^2(0) + \mu {e^b}{{\hat U}^2} \le  - \mu b \left\| {{z_3}} \right\|_b^2 - \mu z_3^2(0) \\
				&+ 2\mu {e^b}({{\bar {\hat \alpha} }^2} + 2{{\bar {\hat \beta} }^2}N_1^2){\left\| {{f_2}} \right\|^2} + 4\mu {e^b}{{\bar {\hat \beta} }^2}N_2^2{\left\| {{h_2}} \right\|^2}\\
				&\le  - \mu b{V_6} + {c_5}{e^{b + a}}{V_7} + {c_6}{e^b}{V_8} ,
			\end{aligned}
		\end{footnotesize}
		\vspace{-2mm}
	\end{equation}
	where 
	\begin{footnotesize}
		${c_5} = 2\mu ({{\bar {\hat \alpha} }^2} + 2{{\bar {\hat \beta} }^2}N_1^2)$,
    	${c_6} = 4\mu {{\bar {\hat \beta} }^2}N_2^2$
    \end{footnotesize}
	are two positive constants.
	
	\subsubsection{${V_7} = \left\| {{f_2}} \right\|_{ - a}^2$}
	From dynamic \eqref{target system2 a}, we get
	\vspace{-2mm}
	\begin{equation}
		\begin{footnotesize}
			\begin{aligned}
				&{{\dot V}_7} = 2{I_{ - a}}[{f_2}{\partial _t}{f_2}]\\
				&= 2{I_{ - a}}[{f_2}( - \lambda {\partial _x}{f_2} + {{\hat m}_1}{f_2} + {{\hat m}_1}\hat e + {{\hat m}_2}{h_2} + {{\hat m}_2}\hat \tau  + {n^T}\dot {\hat M}\\
				&+ {{\hat m}_2}\int_0^x {{{\hat \alpha }_I}(x,\xi )} {f_2}(\xi )d\xi  + {{\hat m}_2}\int_0^x {{{\hat \beta }_I}(x,\xi )} {h_2}(\xi )d\xi )].
			\end{aligned}
		\end{footnotesize}
		\vspace{-2mm}
	\end{equation}
	With \eqref{target system2 c}, inequality \eqref{cauchy a}, \eqref{cauchy c} and \eqref{estimated error boundary}, boundary condition \eqref{hat parameter bounds} and \eqref{approximate kernel boundary}, we find
	\vspace{-2mm}
	\begin{equation} \label{V7}
		\begin{footnotesize}
			\begin{aligned}
				&{{\dot V}_7} \le  - \lambda a\left\| {{f_2}} \right\|_{ - a}^2 + 2\lambda {q^2}{h_2}^2(0) + 2\lambda {q^2}{{\hat \tau }^2}(0)\\
				&+ \bar m_1^2\left\| {{f_2}} \right\|_{ - a}^2 + \left\| {{f_2}} \right\|_{ - a}^2 + \bar m_1^2\left\| {{f_2}} \right\|_{ - a}^2 + \left\| {\hat e} \right\|_{ - a}^2 + \bar m_2^2\left\| {{f_2}} \right\|_{ - a}^2\\
				&+ \left\| {{h_2}} \right\|_{ - a}^2 + \bar m_2^2\left\| {{f_2}} \right\|_{ - a}^2 + \left\| {\hat \tau } \right\|_{ - a}^2 + \left\| {{f_2}} \right\|_{ - a}^2 + \left\| {{n^T}\dot {\hat M}} \right\|_{ - a}^2\\
				&+ \left\| {{f_2}} \right\|_{ - a}^2 + \bar m_2^2{{\bar {\hat \alpha} }_I}^2\left\| {{f_2}} \right\|_{ - a}^2 + \left\| {{h_2}} \right\|_{ - a}^2 + \bar m_2^2{{\bar {\hat \beta} }_I}^2\left\| {{f_2}} \right\|_{ - a}^2\\
				&\le  - (\lambda a - {c_7}){V_7} + {l_6}[{V_1} + {V_2}] + {l_7}[{V_4} + {V_5}]+ 2{V_8} + 2{e^\delta }{V_9} \\
				&+ 2{V_{10}}+ {l_8}{z^2}(0) + 2\lambda {q^2}{h_2}^2(0) + 4\lambda {q^2}{\tau ^2}(0),		
			\end{aligned}
		\end{footnotesize}
		\vspace{-2mm}
	\end{equation}
	where 
	\begin{footnotesize}
		${l_6} = {e^a}(2{\left| {\tilde M} \right|^2} + {\left| {\dot {\hat M}} \right|^2})$,
		${l_7} = 2{\left| {\tilde N} \right|^2}$,
		${l_8} = 4\lambda {q^2}{\left| {\tilde N} \right|^2}$
	\end{footnotesize}
	are integrable functions and 
	\begin{footnotesize}
		$ {c_7} = (2\bar m_1^2 + 2\bar m_2^2+\bar m_2^2{{\bar {\hat \alpha} }_I}^2+\bar m_2^2{{\bar {\hat \beta} }_I}^2+3)$
	\end{footnotesize}
	are positive constants.
	
	\subsubsection{${V_8} = \left\| {{h_2}} \right\|_b^2$}
	By using dynamic \eqref{target system2 b}, we obtain
	\vspace{-2mm}
	\begin{equation}\label{V8_1}
		\begin{footnotesize}
			\begin{aligned}
				&{{\dot V}_8} = 2{I_b}[{h_2}{\partial _t}{h_2}] = 2{I_b}[{h_2}(\mu {\partial _x}{h_2})]\\
				&- 2{I_b}[{h_2}\int_0^x {\left[ {{{\hat \alpha }_t}(x,\xi ) - \lambda {{\tilde \alpha }_\xi }(x,\xi ) + \mu {{\tilde \alpha }_x}(x,\xi ) + {{\hat m}_4}\alpha (x,\xi )} \right.} \\
				&\;\;\;\;\;\;\;\;\;\;\;\;\;\;\; - {{\hat m}_1}\tilde \alpha (x,\xi ) - {{\hat m}_3}\tilde \beta (x,\xi )]\hat \psi (\xi )d\xi ]\\
				&- \cdots \\
				&+ 2\mu {I_b}[{h_2}\hat \beta (x,0)h_2 (0)] + 2{{\hat m}_3}{I_b}[{h_2}\hat e] + 2{{\hat m}_4}{I_b}[h_2 \hat \tau ] \\
				&+ 2{{\hat m}_4}{I_b}[{h_2}\hat \varphi (x)] + 2{I_b}[h_2 {z^T}(x )\dot {\hat N}].
			\end{aligned}
		\end{footnotesize}
		\vspace{-2mm}
	\end{equation}
	Then, by employing inequalities \eqref{cauchy a}, as well as boundary conditions \eqref{hat parameter bounds} and \eqref{approximate kernel boundary}, we have
	\vspace{-2mm}
	\begin{equation}\label{V8_2}
		\begin{footnotesize}
			\begin{aligned}
				&{{\dot V}_8} \le  - \mu {h_2}^2(0) + (25 - \mu b + {\lambda ^2}{q^2}{{\bar {\hat \alpha} }^2} + {\mu ^2}{{\bar {\hat \beta} }^2})\left\| {{h_2}} \right\|_b^2\\
				&+ (\bar {\hat \alpha} _t^2 + {\lambda ^2}\bar {tilde \alpha} _\xi ^2 + {\mu ^2}\bar {\tilde \alpha} _x^2 + \bar m_4^2{{\bar \alpha }^2} + \bar m_1^2{{\bar {\tilde \alpha} }^2}\\
				&\;\;\;\;\;\; + \bar m_3^2{{\bar {\tilde \beta} }^2} + {\lambda ^2}{{\bar {\tilde \alpha} }^2} + {\mu ^2}{{\bar {\tilde \alpha} }^2}){e^b}{\left\| {{f_2}} \right\|^2}\\
				&+ (\bar {\hat \beta} _t^2 + {\mu ^2}\bar {\tilde \beta} _\xi ^2 + {\mu ^2}\bar {\tilde \beta} _x^2 + \bar m_4^2{{\bar {\hat \beta} }^2} + \bar m_2^2{{\bar {\tilde \alpha} }^2} + \bar m_4^2){e^b}{\left\| {\hat \varphi } \right\|^2}\\
				&+ (\bar m_1^2{{\bar {\hat \alpha} }^2} \!+\! \bar m_3^2{{\bar {\hat \beta} }^2} \!+\! \bar m_3^2){e^b}{\left\| {\hat e} \right\|^2} \!+\! (\bar m_2^2{{\bar {\hat \alpha} }^2} \!+\! \bar m_4^2{{\bar {\hat \beta} }^2} \!+\! \bar m_4^2){e^b}{\left\| {\hat \tau } \right\|^2}\\
				&+ {e^b}{{\bar {\hat \alpha} }^2}{\left| {\dot {\hat M}} \right|^2}{\left\| n \right\|^2} \!\!+\! ({{\bar {\hat \beta} }^2} \!\!+\! 1){e^b}{\left| {\dot {\hat N}} \right|^2}{\left\| z \right\|^2} \!\!+\! 3{e^b}{h_2}^2(0) \!+\! 2{e^b}{{\hat \tau }^2}(0).
			\end{aligned}
		\end{footnotesize}
		\vspace{-2mm}
	\end{equation}
	From equations \eqref{real status boundary} and \eqref{estimated error boundary}, using property \eqref{property 2}, we can obtain
	\vspace{-2mm}
	\begin{equation}\label{V8_3}
		\begin{footnotesize}
			\begin{aligned}
				&{{\dot V}_8} \le - (\mu b - {c_8})\left\| h \right\|_b^2{V_8} + {l_9}[{V_1} + {V_2}] + {l_{10}}[{V_4} + {V_5}]\\
				&+ {l_{11}}{V_7} + {l_{12}}{V_8} + {c_9}{e^{b + a}}{V_9} + {c_{10}}{e^b}{V_{10}}\\
				&+ {l_{13}}{{z^2}(0)} + {l_{14}}{h_2^2}(0) - \mu {h^2_2}(0) + 4{e^b}{\tau ^2}(0),
			\end{aligned}
		\end{footnotesize}
		\vspace{-2mm}
	\end{equation}
	where 
	\vspace{-2mm}
	\begin{equation}
		\begin{footnotesize}
			\begin{aligned}
				{l_9} &= {e^{b + a}}[(2\bar m_1^2{{\bar {\hat \alpha} }^2} + 2\bar m_3^2{{\bar {\hat \beta} }^2} + 2\bar m_3^2){\left| {\tilde M} \right|^2} + {{\bar {\hat \alpha} }^2}{\left| {\dot {\hat M}} \right|^2}], \\
				{l_{10}} &= {e^b}[(2\bar m_2^2{{\bar {\hat \alpha} }^2} + 2\bar m_4^2{{\bar {\hat \beta} }^2} + 2\bar m_4^2){\left| {\tilde N} \right|^2} + ({{\bar {\hat \beta} }^2} + 1){\left| {\dot {\hat N}} \right|^2}], \\
				{l_{11}} &=\! {e^{b + a}}[\bar {\hat \alpha} _t^2 \!+\! {\lambda ^2}\bar {\tilde \alpha} _\xi ^2 \!+\! {\mu ^2}\bar {\tilde \alpha} _x^2 \!+\! \bar m_4^2{{\bar \alpha }^2} +\! \bar m_1^2{{\bar {\tilde \alpha} }^2} \!+\! \bar m_3^2{{\bar {\tilde \beta} }^2} \!+\!  {\lambda ^2}{{\bar {\tilde \alpha} }^2} \!+\! {\mu ^2}{{\bar {\tilde \alpha} }^2} \\
				&\;\;\;\;\;\; +\! 2N_1^2(\bar {\hat \beta} _t^2 \!+\! {\mu ^2}\bar {\tilde \beta} _\xi ^2 \!+\! {\mu ^2}\bar {\tilde \beta} _x^2 \!+\! \bar m_4^2{{\bar {\hat \beta} }^2} \!+\! \bar m_2^2{{\bar {\tilde \alpha} }^2} + \bar m_4^2)],\\
				{l_{12}} &= 2N_2^2(\bar {\hat \beta} _t^2 + {\mu ^2}\bar {\tilde \beta} _\xi ^2 + {\mu ^2}\bar {\tilde \beta} _x^2),\quad {l_{13}} = 4{\left| {\tilde N} \right|^2}{e^b}, \quad
				{l_{14}} = 3{e^b}
			\end{aligned}
			\nonumber
		\end{footnotesize}
		\vspace{-2mm}
	\end{equation}
	\vspace{-2mm}
	are integrable functions and 
	\begin{equation}
		\begin{footnotesize}
			\begin{aligned}
				{c_8} &= 25 + {\lambda ^2}{q^2}{{\bar {\hat \alpha} }^2} + {\mu ^2}{{\bar {\hat \beta} }^2} + 2N_2^2(\bar m_4^2{{\bar {\hat \beta} }^2} + \bar m_2^2{{\bar {\tilde \alpha} }^2} + \bar m_4^2),\\
				{c_9} &= 2(\bar m_1^2{{\bar {\hat \alpha} }^2} + \bar m_3^2{{\bar {\hat \beta} }^2} + \bar m_3^2), \quad {c_{10}} = 2(\bar m_2^2{{\bar {\hat \alpha} }^2} + \bar m_4^2{{\bar {\hat \beta} }^2} + \bar m_4^2)
			\end{aligned}
			\nonumber
		\end{footnotesize}
		\vspace{-2mm}
	\end{equation}
	are positive constants.
	
	\subsubsection{${V_9} = \left\| e \right\|_{ - a}^2$ and ${V_{10}} = \left\| \tau  \right\|_b^2$}
	By using \eqref{property 1}, \eqref{real error bound 1} and \eqref{real error bound 2}, we obtain
    \begin{subequations}\label{V9_10}
    	\begin{footnotesize}
    	\begin{align}
    	&\begin{aligned}
    		{{\dot V}_9} &=  - 2\lambda {I_{ - a}}[e{e_x}] \!=\!  -\! \lambda a\left\| e \right\|_{ - a}^2 \!-\! \lambda {e^{ - a}}{e^2}(1) \!+\! \lambda {e^2}(0)\\
    		&\le  - \lambda a\left\| e \right\|_{ - a}^2 \!=\!  - \lambda a{V_9}
    	\end{aligned}\\
    	&\begin{aligned}
    		{{\dot V}_{10}} &= 2\mu {I_b}[\tau {\tau _x}] \!= \!- \mu b\left\| \tau  \right\|_b^2 \!+\! \mu {e^b}{\tau ^2}(1) \!-\! \mu {\tau ^2}(0),\\
    		&\le  - \mu b\left\| \tau  \right\|_b^2 \!=\!  - \mu b{V_{10}}.
    	\end{aligned}
        \end{align}
    	\end{footnotesize}
    \end{subequations}
}

%\section*{Acknowledgment}

%\section*{Acknowledgment}

\ifCLASSOPTIONcaptionsoff
  \newpage
\fi

\bibliographystyle{IEEEtran}
\bibliography{ref.bib}

\vspace{-22pt}
\begin{IEEEbiography}[{\includegraphics[width=1in,height=1.25in,clip,keepaspectratio]{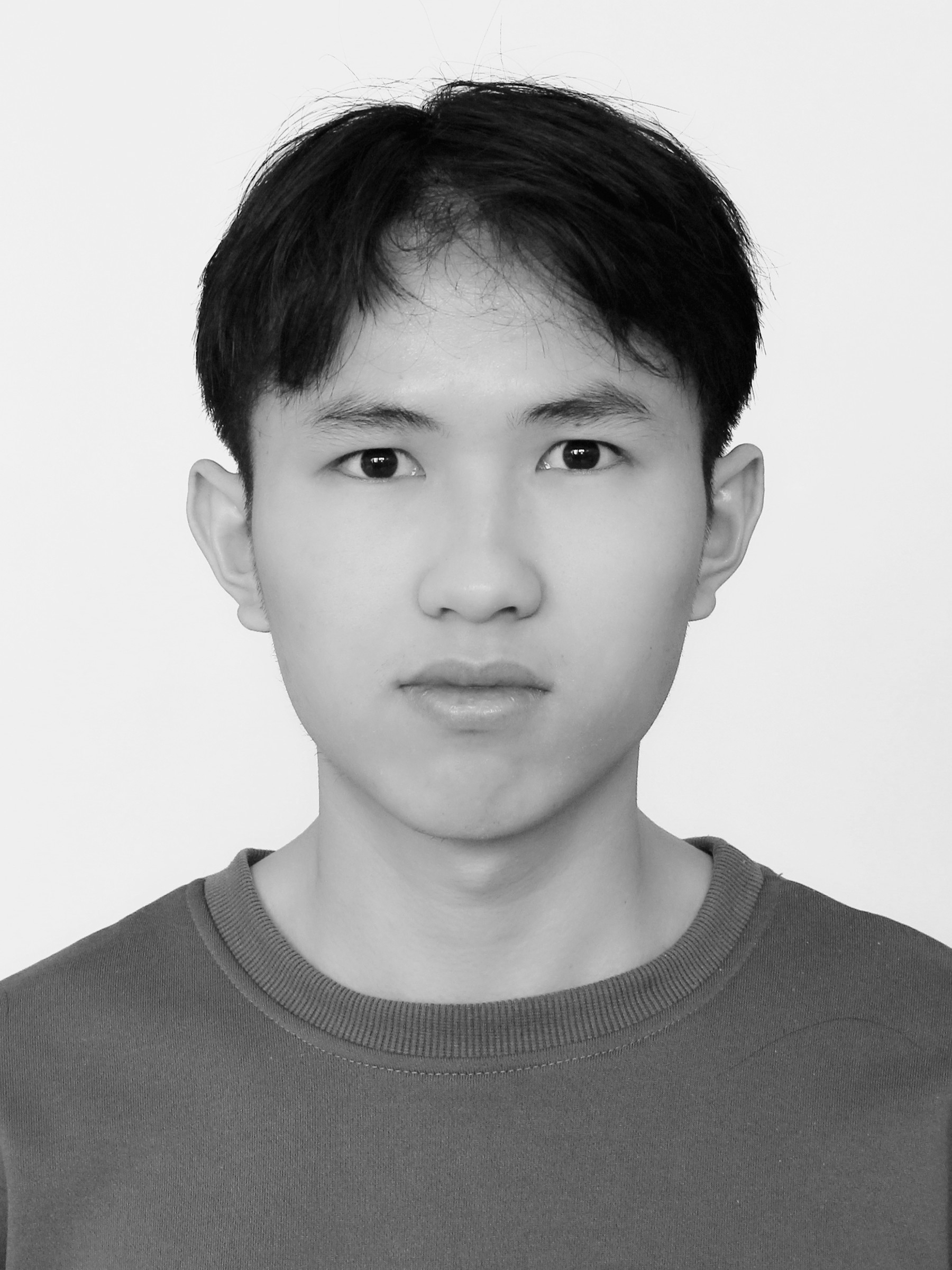}}]{Xianhe Zhang}
	received the B.Eng. degree in automation from the Nanjing Agricultural University, Nanjing, China, in 2023. He is currently pursuing the M.E. degree in control theory and control engineering with the school of Automation, Central South University, Changsha, China.
	
	His research involved adaptive control of dynamical systems.
\end{IEEEbiography}

\begin{IEEEbiography}[{\includegraphics[width=1in,height=1.25in,clip,keepaspectratio]{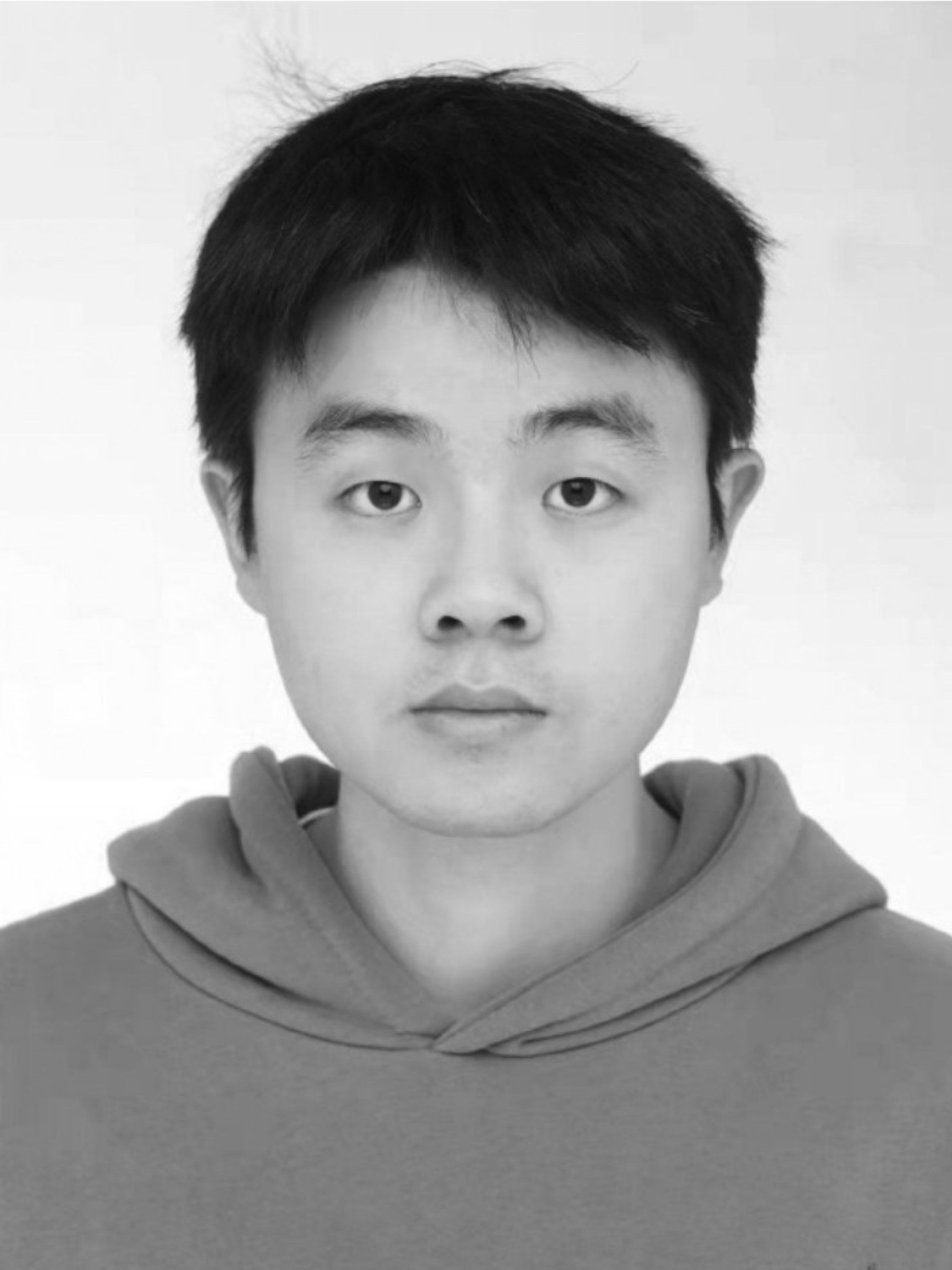}}]{Xiao Yu}
	received the B.Eng. degree in automation from Central South University, Changsha, China, in 2021. He is currently pursuing the Ph.D. degree in control theory and control engineering with the school of Automation, Central South University, Changsha, China.
	
	His research interests include the internal model control of infinite-dimensional systems and adaptive neural-network control.
\end{IEEEbiography}

\begin{IEEEbiography}[{\includegraphics[width=1in,height=1.25in,clip,keepaspectratio]{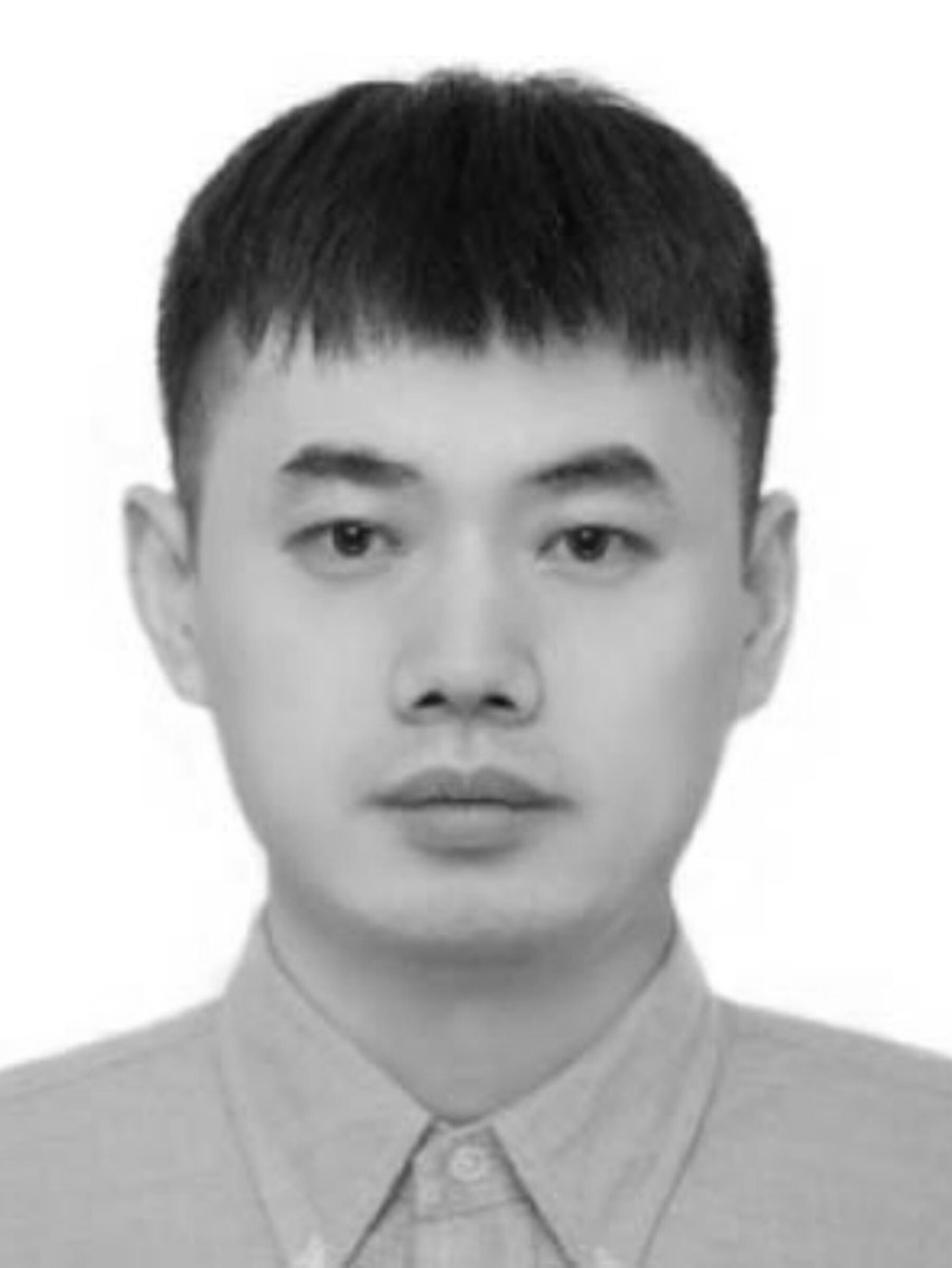}}]{Xiaodong Xu}
	(Member, IEEE) received the B.E. degree in process control from the Beijing Institute of Technology, Beijing, China, and the Ph.D. degree in process control from the University of Alberta, Edmonton, AB, Canada, in 2010 and 2017, respectively.
	
	He is currently a Full Professor with the School of Automation, Central South University, Changsha, China. His research involved robust/optimal control and fault estimation of infinite-dimensional systems including energy systems.
\end{IEEEbiography}

\begin{IEEEbiography}[{\includegraphics[width=1in,height=1.25in,clip,keepaspectratio]{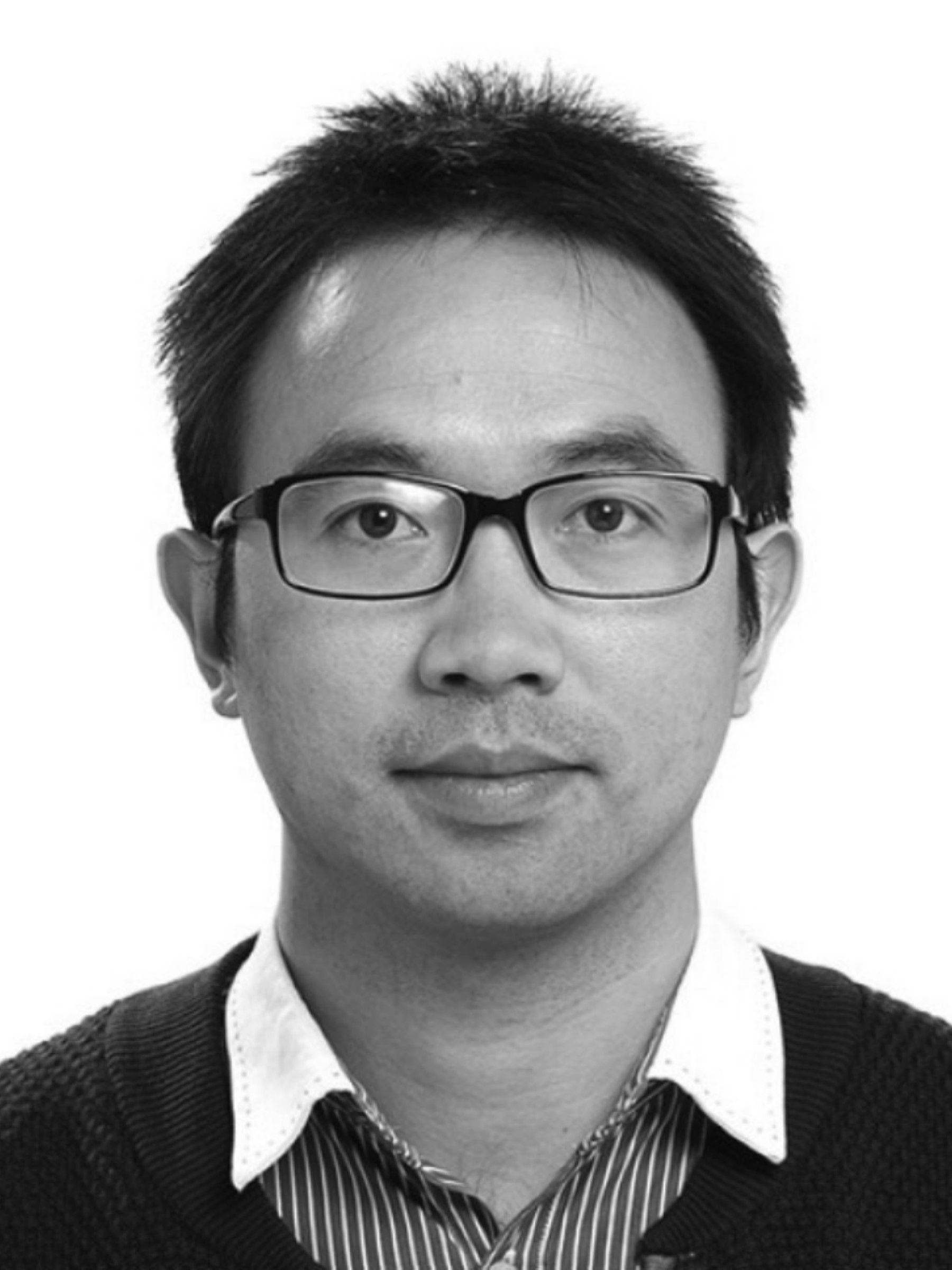}}]{Biao Luo}
	(Senior Member, IEEE) received the Ph.D. degree in control science and engineering from Beihang University, Beijing, China, in 2014.
	
	From 2014 to 2018, he was an Associate Professor and an Assistant Professor with the Institute of Automation, Chinese Academy of Sciences, Beijing. He is currently a Full Professor with the School of Automation, Central South University, Changsha, China. His current research interests include distributed parameter systems, intelligent control, reinforcement learning, and computational intelligence.
\end{IEEEbiography}

\vfill

\end{document}

\title{Intelligent acceleration adaptive control of linear $2\times2$ hyperbolic PDE systems}
	
	\author{Xianhe~Zhang, Yu~Xiao, Xiaodong~Xu,\IEEEmembership{ Member,~IEEE}, Biao Luo,\IEEEmembership{ Senior Member,~IEEE} \vspace{0mm}
		
		\thanks{
			Xianhe Zhang, Yu Xiao, Xiaodong Xu, Biao Luo, Zhiyong Chen are with the School of Automation, Central South University, Changsha, China (e-mail:
			zxh\_csu@csu.edu.cn, yu\_xiao@csu.edu.cn, xx1@ualberta.ca, biao.luo@hotmail.com).}
		%		\thanks{Zhiyong Chen is with the School of Electrical Engineering and Computing, University of Newcastle, Australia  (e-mail:
			%			zhiyong.chen@newcastle.edu.au).}
	}


	
	\maketitle
	
	\begin{abstract}
	
	\end{abstract}

	\begin{IEEEkeywords}
		Adaptive control, Intelligent acceleration, Hyperbolic PDE, DeepONet, Neural operator.
	\end{IEEEkeywords}

	\section{Introduction}
	\subsection{Notation}
	Throughout the paper, $\mathbb{R}^n$ denotes the  $n$ dimensional space, and  the corresponding Euclidean norms are denoted $\left|\cdot\right|$.
	For a signal $\psi (x,t)$ defined on $0 \le x \le 1,t \ge 0$, $\left\| \psi \right\|$ denotes the $L_2$-norm 
	$\left\| {\psi(t)} \right\| = \sqrt {\int_0^1 {{\psi^2}(x,t)dx} }$.
	For a time-varying signal $\psi\left(x,t\right)\in \mathbb{R}$, the $L_p$-norm is denoted by
	$\psi  \in {L_p} \Leftrightarrow {\left\| \psi  \right\|_p} = {(\int_0^\infty  {{{\left| {\psi (t)} \right|}^p}dt} )^{\frac{1}{p}}} < \infty $.
	For the vector signal $\psi(x)$ for all $x\in\left[0,1\right]$,
	we introduce the following operator 
	${I_\delta }[\psi] = \int_0^1 {{e^{\delta x}}\psi(x)dx} $
	with the derived norm 
	${I_\delta }[\psi^T{\psi}] = {\left\| \psi \right\|_\delta^2 } = \int_0^1 {{e^{\delta x}}{\psi^T}(x)\psi(x)dx} $.
	The following properties can be derived
	\begin{equation} \label{property 1}
		\setlength{\abovedisplayskip}{2pt}
		\setlength{\belowdisplayskip}{1pt}
		{I_\delta }[\psi{\psi_x}] = \frac{1}{2}({e^\delta }{\psi^2}(1) - {\psi^2}(0) - \delta {\left\| \psi \right\|_\delta^2 })
	\end{equation}
	\begin{equation} \label{property 2}
		\setlength{\abovedisplayskip}{1pt}
		\setlength{\belowdisplayskip}{2pt}
		\left\|  \cdot  \right\|_{ - \delta }^2 \le {\left\|  \cdot  \right\|^2} \le {e^\delta }\left\|  \cdot  \right\|_{ - \delta }^2 \le \left\|  \cdot  \right\|_\delta ^2,\delta  \ge 1
	\end{equation}
	
	Moreover, use the shorthand notation ${\psi _x} = {\partial _x}\psi  = \frac{{\partial \psi }}{{\partial x}}$ and ${\psi _t} = {\partial _t}\psi  = \frac{{\partial \psi }}{{\partial t}}$.
	
	\begin{figure}[!t] 
		\centering
		%\subfloat[parameter estimation]
		{\includegraphics[width=0.7\linewidth]{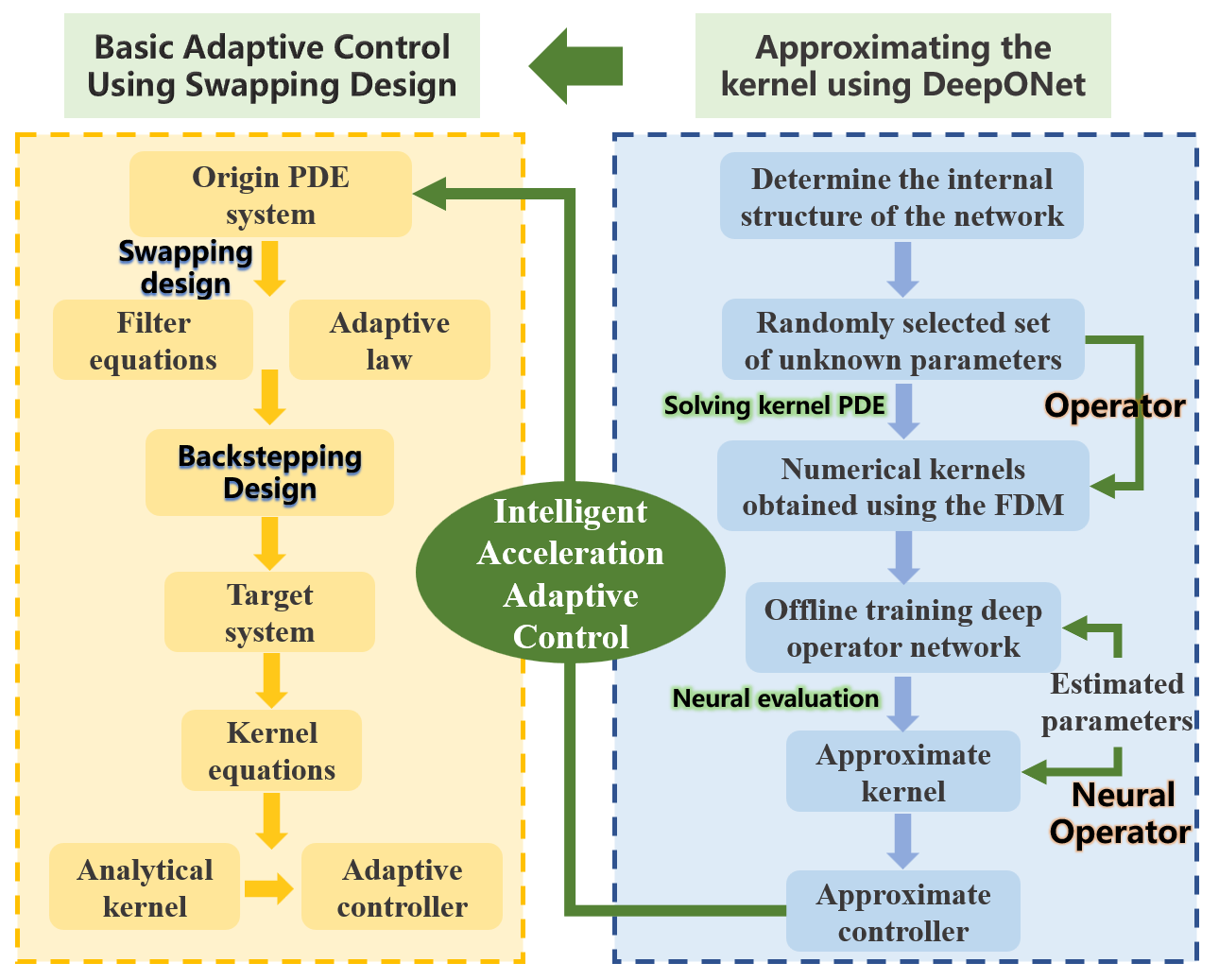}} 
		\caption{The intelligent acceleration adaptive control scheme combined with DeepONet. The approximated kernel obtained from the off-line trained operator network instead of the analytic kernel in numerical computation.}
		\label{accelerated control}
	\end{figure}

	\section{Problem statement}
	Consider the following linear $2\times2$ hyperbolic PDE systems with constant in-domain coefficients
	\begin{subequations}\label{original_system}
		\setlength{\abovedisplayskip}{2pt}
		\setlength{\belowdisplayskip}{2pt}
		\begin{align}
			{\psi _t}(x,t) + \lambda {\psi_x}(x,t) &= {m_1}\psi(x,t) + {m_2}\varphi(x,t)\\
			{\varphi_t}(x,t) - \mu {\varphi_x}(x,t) &= {m_3}\psi(x,t) + {m_4}\varphi(x,t)\\
			\psi(0,t) &= q\varphi(0,t)\\
			\varphi (1,t) &= U(t)
		\end{align}
	\end{subequations}
	defined for $0 \le x \le 1,t \ge 0$, where $\psi(x,t)$,$\varphi(x,t)$ are the system states, $U$ is the control input, and $0 < \lambda  \in \mathbb{R}$, $0 < \mu  \in \mathbb{R}$ are known transport speeds. The boundary parameter $q$ is known constant,  while the coupling coefficient ${m_1}, {m_2}, {m_3}, {m_4}$ are unknown. The initial condition $\psi (x,0) = {\psi _0}(x)$, $\varphi (x,0) = {\varphi _0}(x)$ are assumed to satisfy ${\psi _0}(x),{\varphi _0}(x) \in {L_2}$. Since the system bounds are not strictly defined, we assume that there have some bounds on ${m_i}$ defined as ${{{\bar m}_i}}$ so that 
	\begin{equation} \label{parameter bounds}
		\setlength{\abovedisplayskip}{1pt}
		\setlength{\belowdisplayskip}{1pt}
		\left| {{m_i}} \right| \le {{\bar m}_i},i = 1,2,3,4 .
	\end{equation}
	
	The control goal is to design an adaptive state controller that achieves regulation of the system states $\psi$ and $\varphi$ to zero equilibrium profile, while ensuring that all additional signals remain within certain bounds.
	
	\section{Basic adaptive control using swapping design} \label{basic adaptive control}
	
	\subsection{Filter equation}
	Firstly, we introduce the parameter filters that are defined on $0 \le x \le 1,t \ge 0$
	\begin{subequations}\label{filter1}
		\begin{align}
			{\partial_t}{n_1}(x,t)+\lambda {\partial_x}{n_1}(x,t) &= \psi (x,t)&&{n_1}(0,t) = 0  \label{filter1 a} \\
			{\partial_t}{n_2}(x,t)+\lambda {\partial_x}{n_2}(x,t) &= \varphi (x,t)&&{n_2}(0,t) = 0 \label{filter1 b} \\
			{\partial_t}{n_3}(x,t)+\lambda {\partial_x}{n_3}(x,t) &= 0&&{n_3}(0,t) = \varphi (0,t) \label{filter1 c}
		\end{align}
	\end{subequations}
	\begin{subequations}\label{filter2}
		\begin{align}
			{\partial _t}{z_1}(x,t)-\mu {\partial_x}{z_1}(x,t)&= \psi (x,t)&&{z_1}(1,t) = 0  \label{filter2 a} \\
			{\partial _t}{z_2}(x,t)-\mu {\partial_x}{z_2}(x,t)&= \varphi (x,t)&&{z_2}(1,t) = 0 \label{filter2 b} \\
			{\partial _t}{z_3}(x,t)-\mu {\partial_x}{z_3}(x,t)&= 0&&{z_3}(1,t) = U(t) \label{filter2 c}
		\end{align}
	\end{subequations}
	where the initial states conditions ${n_i}(x,0) = n_i^0(x),{z_i}(x,0) = z_i^0(x),i = 1,2,3$ satisfy $n_1^0,n_2^0,n_3^0,z_1^0,z_2^0,z_3^0 \in {L_2}$.
	
	The following aggregated symbols are defined for  analysis easily
	\begin{subequations}\label{zn}
		\begin{align}
			n(x,t) = {[{n_1}(x,t),{n_2}(x,t)]^T}\\
			z(x,t) = {[{z_1}(x,t),{z_2}(x,t)]^T},
		\end{align}
	\end{subequations}
	\begin{align} \label{MN}
		M = {[{m_1},{m_2}]^T},N = {[{m_3},{m_4}]^T}.
	\end{align}
	Then we introduce the non-adaptive estimate of the system states 
	\begin{subequations}\label{non-adaptive states}
		\begin{align}
			\bar \psi (x,t) &= {n^T}(x,t)M + {n_3}(x,t)q\\
			\bar \varphi (x,t) &= {z^T}(x,t)N + {z_3}(x,t),
		\end{align}
	\end{subequations}
	With parameter filters \eqref{filter1} and \eqref{filter2}, we obtain
	\begin{subequations}\label{non-adaptive system}
		\begin{align}
			{{\bar \psi }_t}(x,t) + \lambda {{\bar \psi }_x}(x,t) &= {m_1}\psi (x,t) + {m_2}\varphi (x,t)\\
			{{\bar \psi }_t}(x,t) - \mu {{\bar \varphi }_x}(x,t) &= {m_3}\psi (x,t) + {m_4}\varphi (x,t)\\
			\bar \psi (0,t) &= q\varphi (0,t)\\
			\bar \varphi (1,t) &= U(t),
		\end{align}
	\end{subequations}
	The system \eqref{non-adaptive system} states estimation errors
	\begin{subequations}\label{non-adaptive errors}
		\begin{align}
			e(x,t) &= \psi (x,t) - \bar \psi (x,t)\\
			\tau (x,t) &= \varphi (x,t) - \bar \varphi (x,t),
		\end{align}
	\end{subequations}
	From \eqref{original_system}, we can straightforwardly obtain
	\begin{subequations}\label{non-adaptive errors system}
		\begin{align}
			{e_t}(x,t) + \lambda {e_x}(x,t) &= 0\\
			{\tau _t}(x,t) - \mu {\tau _x}(x,t) &= 0\\
			e(0,t) &= 0 \label{real error bound 1}\\ 
			\tau (1,t) &= 0 \label{real error bound 2}, 
		\end{align}
	\end{subequations}
	When $ t > \max \{ {\lambda ^{ - 1}},{\mu ^{ - 1}}\} $, the system errors will converge to zero. 
	
	\subsection{Parameter estimation law}
	Then, we introduce the adaptive estimate
	\begin{subequations}\label{adaptive states}
		\begin{align}
			\hat \psi (x,t) &= {n^T}(x,t)\hat M + {n_3}(x,t)q\\
			\hat \varphi (x,t) &= {z^T}(x,t)\hat N + {z_3}(x,t),
		\end{align}
	\end{subequations}
	where $\hat M = {[{{\hat m}_1},{{\hat m}_2}]^T}$,$\hat N = {[{{\hat m}_3},{{\hat m}_4}]^T}$ denote the estimates of coupling coefficient. So that the prediction errors ${\hat e}$,${\hat \tau }$ denotes as
	\begin{subequations}\label{adaptive errors}
		\begin{align}
			&\hat e(x,t) = \psi (x,t) - \hat \psi (x,t)\\
			&\hat \tau (x,t) = \varphi (x,t) - \hat \varphi (x,t).
		\end{align}
	\end{subequations}
	Introducing the following adaptive laws 
	\begin{subequations}\label{adaptive laws}
		\begin{small}
			\begin{align}
				&{{\dot {\hat m}}_1} \!=\! {\mathcal{P}_{{{\bar m}_1}}}\left\{ {{\rho _1}\frac{{\int_0^1 {{n_1}(x)\hat e(x)dx} }}{{1 + {{\left\| {{n_1}} \right\|}^2} + {{\left\| {{n_2}} \right\|}^2} + {{\left\| {{n_3}} \right\|}^2}}}} \right\}\\
				&{{\dot {\hat m}}_2} \!=\! {\mathcal{P}_{{{\bar m}_2}}}\left\{ {{\rho _2}\frac{{\int_0^1 {{n_2}(x)\hat e(x)dx} }}{{1 + {{\left\| {{n_1}} \right\|}^2} + {{\left\| {{n_2}} \right\|}^2} + {{\left\| {{n_3}} \right\|}^2}}}} \right\}\\
				&{{\dot {\hat m}}_3} \!=\! {\mathcal{P}_{{{\bar m}_3}}} \!\! \left\{ \! {{\rho _3} \!\! \left[ \! {\frac{{\int_0^1 {{z_1}(x)\hat \tau (x)dx} }}{{1 \!+\! {{\left\| {{z_1}} \right\|}^2} \!+\! {{\left\| {{z_2}} \right\|}^2}}} \!+\! \frac{{{z_1}(0)\hat \tau (0)}}{{1 \!+\! {{\left| {{z_1}(0)} \right|}^2} \!+\! {{\left| {{z_2}(0)} \right|}^2}}}} \! \right]} \!\! \right\}\\
				&{{\dot {\hat m}}_4} \!=\! {\mathcal{P}_{{{\bar m}_4}}} \!\! \left\{ \! {{\rho _4} \!\! \left[ \! {\frac{{\int_0^1 {{z_2}(x)\hat \tau (x)dx} }}{{1 \!+\! {{\left\| {{z_1}} \right\|}^2} \!+\! {{\left\| {{z_2}} \right\|}^2}}} \!+\! \frac{{{z_2}(0)\hat \tau (0)}}{{1 \!+\! {{\left| {{z_1}(0)} \right|}^2} \!+\! {{\left| {{z_2}(0)} \right|}^2}}}} \! \right]} \!\! \right\}
			\end{align}
		\end{small}
	\end{subequations}
	where ${\rho _1},{\rho _2},{\rho _3},{\rho _4}$ are the given positive gains, 
	we assume that the parameters estimated by the above equation also satisfy bounds \eqref{parameter bounds}, that expressed as
	\begin{equation} \label{hat parameter bounds}
		\left| {{\hat{m}_i}} \right| \le {{\bar m}_i},i = 1,2,3,4 .
	\end{equation}
	the adaptive laws operator $\mathcal{P}$  defined as 
	\begin{equation}
		{\mathcal{P}_{{{\bar m}_i}}}({\theta _i},{{\hat m}_i}) = \left\{ \begin{array}{l}
			0\;\;\; \text{if}\;\left| {{{\hat m}_i}} \right| \ge {{\bar m}_i}\;\& \;{\theta _i}{{\hat m}_i} \ge 0\\
			{\theta _i}\;\;\; \text{others},
		\end{array} \right.
	\end{equation}
	Here we introduce a lemma to illustrate the properties of this operator.
	\begin{lemma} \label{lemma adaptive law operator}
		The adaptive law \eqref{adaptive laws} using the operator symbol $\mathcal{P}$ has the following properties for all $t \ge 0$ : (Proof in \cite[appendix A]{anfinsen2018adaptive})
		\begin{subequations}\label{lemma adaptive law}
			\begin{align}
				1)\quad &{\cal P}_{\bar m}^2(\theta ,\hat m) \le {\theta ^2}\\
				2)\quad &- \tilde m{{\cal P}_{\bar m}}(\theta ,\hat m) \le  - {{\tilde m}^T}\theta , \quad \tilde m = m - \hat m \\
				3)\quad &\left| {{{\hat m}_i}} \right| \le \bar m, \quad i = 1...n,
			\end{align}
		\end{subequations}
	\end{lemma}
	
	\subsection{Adaptive control}
	Associative initial system \eqref{original_system}, parameter filters \eqref{filter1} and \eqref{filter2}, the dynamic states estimation system can be formulated as
	\begin{small}
		\begin{subequations}\label{dynamic system}
			\setlength{\abovedisplayskip}{5pt}
			\setlength{\belowdisplayskip}{5pt}
			\begin{align}
				{{\hat \psi }_t}(x,t) \!+\! \lambda {{\hat \psi }_x}(x,t) &={{\hat m}_1}\psi (x,t) \!+\! {{\hat m}_2}\varphi (x,t) \!+\! {n^T}(x,t)\dot {\hat M} \label{dynamic system a} \\ 
				{{\hat \varphi }_t}(x,t) \!-\! \mu{{\hat \varphi }_x}(x,t)&={{\hat m}_3}\psi (x,t)\!+\!{{\hat m}_4}\varphi(x,t) \!+\! {z^T}(x,t)\dot {\hat N} \label{dynamic system b} \\
				\hat \psi (0,t) &= q\varphi (0,t) \label{dynamic system c} \\
				\hat \varphi (1,t) &= U(t). \label{dynamic system d}
			\end{align}
		\end{subequations}
	\end{small}
	Then we apply the backstepping transformation to map the dynamic into target system, using the boundary control law
	\begin{small}
		\begin{equation}
			\setlength{\abovedisplayskip}{3pt}
			\setlength{\belowdisplayskip}{3pt}
			\begin{aligned}
				U(t) \!=\! \int_0^1 \!{\alpha (1,\xi ,t)} \hat \psi (\xi ,t)d\xi  \!+\! \int_0^1 \!{\beta (1,\xi ,t)} \hat \varphi (\xi ,t)d\xi ,
			\end{aligned}
		\end{equation}
	\end{small}
	where $\alpha$ and $\beta$ represent the resolvent kernel of the Volterra integral transformations, and satisfy the following kernel equation (\textit{Here t is omitted for brevity}). Note that the numerical kernel is obtained by solving the following PDE, contrasted with the approximated kernel.
	\begin{small}
		\begin{subequations}\label{kernel equation}
			\setlength{\abovedisplayskip}{2pt}
			\setlength{\belowdisplayskip}{2pt}
			\begin{align}
				\mu {\alpha _x}(x,\xi ) - \lambda {\alpha _\xi }(x,\xi ) &= ({{\hat m}_1} - {{\hat m}_4})\alpha (x,\xi ) + {{\hat m}_3}\beta (x,\xi )\\
				\mu {\beta _x}(x,\xi ) + \mu {\beta _\xi }(x,\xi ) &= {{\hat m}_2}\alpha (x,\xi )\\
				\beta (x,0) &= q\frac{\lambda }{\mu }\alpha (x,0)\\
				\alpha (x,x) &=  - \frac{{{{\hat m}_3}}}{{\lambda  + \mu }},
			\end{align}
		\end{subequations}
	\end{small}
	%\resizebox{0.9\hsize}{!}{$
		defined on
		\vspace{-1mm}
		\begin{align}\label{region T}
			T=\left\{ {(x,\xi )\left| {0 \le \xi  \le x \le 1} \right.} \right\} \times \left\{ {t \ge 0} \right\}.
			\vspace{-1mm}
		\end{align} 
		Moreover, the backstepping kernel is invertible with inverse
		\begin{small}
			\begin{subequations}\label{inverse kernel}
				\setlength{\abovedisplayskip}{5pt}
				\setlength{\belowdisplayskip}{5pt}
				\begin{align}
					{\alpha _I}(x,\xi,t ) = \alpha (x,\xi,t ) + \int_\xi ^x {\beta (x,s,t){\alpha _I}(s,\xi,t )ds} \\
					{\beta _I}(x,\xi,t ) = \beta (x,\xi,t ) + \int_\xi ^x {\beta (x,s,t){\beta _I}(s,\xi,t )ds} .
				\end{align}
			\end{subequations}
		\end{small}
		With the above kernel function, we define the following backstepping transformation 
		\begin{small}
			\begin{subequations}\label{backstepping transformation}
				\setlength{\abovedisplayskip}{5pt}
				\setlength{\belowdisplayskip}{5pt}
				\begin{align}
					&f(x,t) = \hat \psi (x,t)\\
					&\begin{aligned}
						h(x,t) = \Gamma [\hat \psi ,\hat \varphi ](x) &= \hat \varphi (x,t)- \int_0^x {\alpha (x,\xi ,t)} \hat \psi (\xi ,t)d\xi  \\
						& - \int_0^x {\beta (x,\xi ,t)} \hat \varphi (\xi ,t)d\xi ,
					\end{aligned}
				\end{align}
			\end{subequations}
		\end{small}
		the inverse form of this transformation is expressed as
		\begin{small}
			\begin{subequations}\label{backstepping transformation inverse}
				\setlength{\abovedisplayskip}{5pt}
				\setlength{\belowdisplayskip}{5pt}
				\begin{align}
					&\hat \psi (x,t) = f(x,t)\\
					&\begin{aligned}
						\hat \varphi (x,t) = {\Gamma ^{ - 1}}[f,h](x) &\!=\! h(x,t) \!+\! \int_0^x \! {{\alpha _I}(x,\xi ,t)} f(\xi ,t)d\xi \\
						&\!+ \int_0^x {{\beta _I}(x,\xi ,t)} h(\xi ,t)d\xi ,
					\end{aligned}
				\end{align}
			\end{subequations}
		\end{small}
		where $\Gamma$ and ${\Gamma^{ -1}}$ denote the transformation operators.
		
		To illustrate the boundedness of the operators defined above, we introduce a lemma as follows.
		\begin{lemma} \label{operator boundary}
			For every time $t \ge 0$, we introduce the following properties for the operators defined in \eqref{region T}
			\begin{small}
				\begin{subequations}\label{kernel boundary}
					\begin{align}
						&\left| {\alpha (x,\xi ,t)} \right| \le \bar \alpha ,   \;\;\;\;\ 
						\left| {\beta (x,\xi ,t)} \right| \le \bar \beta, \label{lemma2a} \\
						&\left| {{\alpha _I}(x,\xi ,t)} \right| \le {{\bar \alpha }_I}, \;\;
						\left| {{\beta _I}(x,\xi ,t)} \right| \le {{\bar \beta }_I}, \label{lemma2b} \\
						&\left| {{\alpha _t}(x,\xi ,t)} \right|,  \   
						\left| {{\beta _t}(x,\xi ,t)} \right| \in {L_2}. \label{lemma2d}
					\end{align}
				\end{subequations}
			\end{small}
		\end{lemma}  
		\begin{IEEEproof}
			For every time $t$, it has been proved in \cite{coron2013local} that $\alpha\left(x,\xi,t\right)$ and $\beta\left(x,\xi,t\right)$ are bounded and unique, and meet the following inequalities for all $(x,\xi ,t) \in T$:
			\begin{small}
				\begin{subequations}
					\setlength{\abovedisplayskip}{3pt}
					\setlength{\belowdisplayskip}{3pt}
					\begin{align}
						\left| {\alpha\left( {x,\xi ,t} \right)} \right|& \le {F_1}\left( {{{\hat m}_1}, \ldots {{\hat m}_4}} \right),  \label{proof2a}\\
						\left| {\beta\left( {x,\xi ,t} \right)} \right| &\le {F_2}\left( {{{\hat m}_1}, \ldots {{\hat m}_4}} \right), \label{proof2b}\\
						\left| {{\alpha_t}\left( {x,\xi ,t} \right)} \right| &\le {G_1}\left( {\left| {{{\dot {\hat m}}_1}} \right| + \left| {{{\dot{ \hat m}}_2}} \right| + \left| {{{\dot {\hat m}}_3}} \right| + \left| {{{\dot {\hat m}}_4}} \right|} \right), \label{proof2c}\\
						\left| {{\beta_t}\left( {x,\xi ,t} \right)} \right| &\le {G_2}\left( {\left| {{{\dot {\hat m}}_1}} \right| + \left| {{{\dot {\hat m}}_2}} \right| + \left| {{{\dot {\hat m}}_3}} \right| + \left| {{{\dot {\hat m}}_4}} \right| } \right), \label{proof2d}
					\end{align}
				\end{subequations}
			\end{small}
			where $F_1\left(\cdot\right)$ and $F_2\left(\cdot\right)$ are continuous functions concerning the variables therein. Since $\hat m_i,i=1,\ldots,4$ are compact due to the projection ${\mathcal{P}_{{{\bar m}_i}}}$, $F_1\left(\cdot\right)$ and $F_2\left(\cdot\right)$ are bounded. Let $\bar \alpha$ and $\bar \beta$ be the maximum values of $F_1\left(\cdot\right)$ and $F_2\left(\cdot\right)$, respectively,  then \eqref{lemma2a} must hold. 
			With the successive approximation method, the solution $\left({\alpha}_I,\beta_I\right)$ to equation \eqref{inverse kernel} is bounded whose bounds $\bar {\alpha}_I$ and $\bar {\beta}_I$ are determined by $\bar {\alpha}$, $\bar {\beta}$.
			Again using the results in \cite{coron2013local}, the time derivative $\left(\alpha_t,\beta_t\right)$ satisfies \eqref{proof2c} and \eqref{proof2d}. Note  $\dot{ \hat{m}}_i , i=1...4$ is bounded convergent, hence $\left| {\alpha_t} \right|$, $\left| {\beta_t} \right|$ are bounded and convergent as well, \eqref{lemma2d} is established.
		\end{IEEEproof}
		
		From \eqref{backstepping transformation} differentiates the time, inserts the dynamic system equation \eqref{dynamic system}, the kernel function transformation \eqref{kernel equation}, \eqref{tran equations} and the boundary filters 
		\eqref{filter1}, \eqref{filter2}, we obtain target system as
		\begin{small}
			\begin{subequations}\label{target system}
				\begin{align}
					&\begin{aligned} \label{target system a}
						&{\partial _t}f(x,t) + \lambda {\partial _x}f(x,t) = {n^T}(x,t)\dot {\hat M} \\
						&+ {{\hat m}_1}f(x,t) + {{\hat m}_1}\hat e(x,t) + {{\hat m}_2}h(x,t) + {{\hat m}_2}\hat \tau (x,t)\\
						&+ \int_0^x {\sigma (x,\xi ,t)} f(\xi ,t)d\xi  + \int_0^x {\vartheta (x,\xi ,t)} h(\xi ,t)d\xi 
					\end{aligned}\\
					&\begin{aligned} \label{target system b}
						&{\partial _t}h(x,t) - \mu {\partial _x}h(x,t) \\
						&= {{\hat m}_4}h(x,t) - \lambda \alpha (x,0)q\hat \tau (0,t) \\
						&- \! \int_0^x \! {{\alpha _t}(x,\xi ,t)} f(\xi ,t)d\xi  \!-\! \int_0^x \! {{\beta _t}(x,\xi ,t)} {{\rm \Gamma}^{ - 1}}[f,h](x)d\xi  \\
						&+ \! \Gamma [{{\hat m}_1}\hat e + {{\hat m}_2}\hat \tau ,{{\hat m}_3}\hat e + {{\hat m}_4}\hat \tau ](x) + \Gamma [{n^T}\dot {\hat M},{z^T}\dot {\hat N}](x),
					\end{aligned}\\
					&f(0,t) = qh(0,t) + q\hat \tau (0,t) \label{target system c}\\
					&h(1,t) = 0 \label{target system d}.
				\end{align}
			\end{subequations}
		\end{small}
		
		Finally, based on the Lyapunov stability theory, the stability of the closed-loop system can be proved by considering the boundary filter \eqref{filter1}, \eqref{filter2}, the adaptive update rate of parameters \eqref{adaptive laws}, the dynamic system equation \eqref{dynamic system}, the kernel function transformation \eqref{kernel equation}, and the adaptive control gain. The system states converges to zero, and other signals are all bounded and integrable (The detailed proof process can be found in \cite[Proof of lemma 8]{anfinsen2018adaptive}).
		
		\section{Approximated kernel by neural operator}
		\begin{equation}
			\setlength{\abovedisplayskip}{3pt}
			\setlength{\belowdisplayskip}{3pt}
			\begin{aligned}\label{deeponet}
				{\Phi _{\cal N}}(u)(y): = \sum\limits_{k = 1}^p {{g^{\cal N}}({u_m};\theta _1^{(k)}){f^{\cal N}}(y;\theta _2^{(k)})},
			\end{aligned}
		\end{equation}
		
		\begin{theorem}\label{universal approximation theorem}(Universal approximation of continuous operators \cite[Theorem 2.1]{lu2021learning}).
			Let ${D_x} \subset {R^{{d_x}}}$, ${D_y} \subset {R^{{d_y}}}$ be compact sets of vectors $x \in {D_x}$ and $y \in {D_y}$. Let $D_u:{D_x} \to {D_u} \subset {R^{{d_u}}}$ is a compact set and $u \in {D_u}$, which be sets of continuous functions $u(x_i)$. let ${D_\Phi }:{D_y} \to {D_\Phi } \subset {R^{{d_\Phi }}}$ be sets of continuous functions of $\Phi (u)(y)$. Assume that $\Phi :{D_u} \to {D_\Phi }$ is a continuous operator. Then, for any $\varepsilon  > 0$, there exist ${m^*},{p^*} \in \mathbb{N}$ such that for each $m \ge {m^*},p \ge {p^*}$, there exist weight coefficients ${\theta _1^{(k)}}$, ${\theta _2^{(k)}}$, neural networks ${{g^{\cal N}}}$, ${{f^{\cal N}}}$, $k = 1,...,p$ and ${x_i} \in {D_x},i = 1,2,...,m$, there exists a neural operator ${\Phi _N}$ such that the following holds:
			\begin{equation}
				\setlength{\abovedisplayskip}{3pt}
				\setlength{\belowdisplayskip}{3pt}
				\begin{aligned}\label{universal approximation}
					\left| {\Phi (u)(y) - {\Phi _N}(u)}(y) \right| < \varepsilon, \;\;\;\;\forall u \in {D_u},y \in {D_y}
				\end{aligned}	
			\end{equation}
			
		\end{theorem}
		where $\textbf{u}_m = {\left[ {u\left( {{x_1}} \right){\rm{, }}u\left( {{x_2}} \right){\rm{, }}...{\rm{, }}u\left( {{x_m}} \right)} \right]^T}$.

		\subsection{Backstepping kernel with neural operators} 
		
		Define the following notations, $ \varpi  = \{ {m_1},{m_2},{m_3},{m_4}\} $, its boundary $\bar \varpi \ge \left| \varpi \right|$, then denote by ${\Phi _\alpha }$, ${\Phi _\beta }$ the operators that maps $\varpi$ to the kernel $\alpha$, $\beta$ that satisfies \eqref{kernel equation}. Expressed as 
		\begin{subequations}\label{kernel maps}
			\setlength{\abovedisplayskip}{5pt}
			\setlength{\belowdisplayskip}{5pt}
			\begin{align}
				&{\Phi _\alpha }:\varpi \to \alpha ,\;\; \alpha = {\Phi _\alpha }(\varpi),\\
				&{\Phi _\beta }:\varpi \to \beta ,\;\; \beta = {\Phi _\beta }(\varpi).
			\end{align}
		\end{subequations}
		According to the proof of Lemma \ref{operator boundary}, the backstepping kernel operator $\alpha$ and $\beta$ satisfy \eqref{proof2a}-\eqref{proof2b}, with their boundaries $\bar{\alpha}$ and $\bar{\beta}$ determined by the continuous functions $F_1\left(\cdot\right)$ and $F_2\left(\cdot\right)$ of the variable $\varpi$. 
		It was shown in \cite{bhan2023neural} that when the boundary of $\varpi$ exists and is Lipschitz continuous, we obtain:
		\begin{subequations}\label{Lipschitz}
			\setlength{\abovedisplayskip}{5pt}
			\setlength{\belowdisplayskip}{5pt}
			\begin{align}
				&{\left\| {{\alpha _1} - {\alpha _2}} \right\|_ \infty } < {C_{\bar \alpha }}{\left\| {{\varpi _1} - {\varpi _2}} \right\|_ \infty }\\
				&{\left\| {{\beta _1} - {\beta _2}} \right\|_ \infty } < {C_{\bar \beta }}{\left\| {{\varpi _1} - {\varpi _2}} \right\|_ \infty }
			\end{align}
		\end{subequations}
		This leads to the operator ${\Phi _\alpha }$, ${\Phi _\beta }$ possessing the Lipschitz property
		\begin{subequations}\label{operator Lipschitz}
			\setlength{\abovedisplayskip}{5pt}
			\setlength{\belowdisplayskip}{5pt}
			\begin{align}
				{\left\| {{\Phi _\alpha }({\varpi _1}) - {\Phi _\alpha }({\varpi _2})} \right\|_\infty } < {C_{\bar \alpha }}{\left\| {{\varpi _1} - {\varpi _2}} \right\|_\infty }\\
				{\left\| {{\Phi _\beta }({\varpi _1}) - {\Phi _\beta }({\varpi _2})} \right\|_\infty } < {C_{\bar \beta }}{\left\| {{\varpi _1} - {\varpi _2}} \right\|_\infty },
			\end{align}
		\end{subequations}
		where ${C_{\bar \alpha }}$ and ${C_{\bar \beta }}$ denote the Lipschitz constants,  which can be expressed as expressions for the boundaries $\bar \alpha$, $\bar \beta$ of the variable function.
		
		Then combined with DeepONet universal approximation theorem (The theorem \ref{universal approximation theorem}), we obtain
		\begin{subequations}\label{kernel maps approximation}
			\setlength{\abovedisplayskip}{5pt}
			\setlength{\belowdisplayskip}{5pt}
			\begin{align}
				\left| {{\Phi _\alpha }(\varpi ) - {\Phi _{\mathcal{N}\alpha }}(\varpi )} \right| &< \varepsilon, \\
				\left| {{\Phi _\beta }(\varpi ) - {\Phi _{\mathcal{N}\beta }}(\varpi )} \right| &< \varepsilon,
			\end{align}
		\end{subequations}
		where ${\varpi _m} = {(\varpi ({x_1}),\varpi ({x_2}),...,\varpi ({x_m}))^T}$ with corresponding ${x_i} \in {D_x}$, for $i \!=\! 1,2,...,m$.

		The approximated kernels here are denoted as
		\vspace{-1mm}
		\begin{equation}\label{hat kernel maps}
			\begin{aligned}
				{\hat{\alpha} } = {\Phi _{\mathcal{N}\alpha }}(\varpi),
				\quad \quad 
				{\hat{ \beta} } = {\Phi _{\mathcal{N}\beta }}(\varpi).
			\end{aligned}
		\end{equation}
		
		\section{Approximate kernels are applied in \\ adaptive control}
		In this section, we utilize approximate kernels \eqref{hat kernel maps} instead of analytical computations for control using swapping design. (Some formulas have been defined in section \ref{basic adaptive control}). 
		\subsection{Adaptive control}
		
		There will be an approximation error between the approximate kernel and the exact kernel, which we define as
		\vspace{-1mm}
		\begin{equation}\label{kernel error}
			\begin{aligned}
				\tilde \alpha  = \alpha  - \hat \alpha, \quad \tilde \beta  = \beta  - \hat \beta .
			\end{aligned}
			\vspace{-1mm}
		\end{equation} 
		The sum of the boundary integrals of the operator kernel forms the approximate controller, denoted as
		\begin{small}
			\begin{equation}
				\begin{aligned}\label{approximate controller}
					\hat U(t) \!=\! \int_0^1 \!{\hat \alpha (1,\xi ,t)} \hat \psi (\xi ,t)d\xi  \!+\! \int_0^1 \!{\hat \beta (1,\xi ,t)} \hat \varphi (\xi ,t)d\xi.
				\end{aligned}
			\end{equation}
		\end{small}
		
		Obviously, the DeepONet operator kernel brings the approximation error term, resulting in an accumulation of errors between the approximate controller and the original controller. Then, the new transformation terms and transformation operators are obtained through PDE backstepping design.
		
		As follows, we consider the new backstepping transformation with the approximate kernel 
		\begin{small}
			\begin{subequations}\label{approximate backstepping transformation}
				\begin{align}
					&{f_2}(x,t) = \hat \psi (x,t) \label{approximate backstepping transformation a}\\
					&\begin{aligned}\label{approximate backstepping transformation b}
						{h_2}(x,t) = {\Gamma _2}[\hat \psi ,\hat \varphi ](x) &= \hat \varphi(x,t) - 		\int_0^x  {\hat \alpha (x,\xi ,t)} \hat \psi (\xi ,t)d\xi \\
						&- \int_0^x {\hat \beta (x,\xi ,t)} \hat \varphi (\xi ,t)d\xi ,
					\end{aligned}	
				\end{align}
			\end{subequations}
		\end{small}
		and its inverse transformation
		\begin{subequations}\label{approximate inverse transformation}
			\setlength{\abovedisplayskip}{2pt}
			\setlength{\belowdisplayskip}{2pt}
			\begin{footnotesize}
				\begin{align}
					&\hat \psi (x,t) = {f_2}(x,t) \label{approximate inverse transformation a}\\
					&\begin{aligned}\label{approximate inverse transformation b}
						\hat \varphi (x,t) = \Gamma _2^{ - 1}[{f_2},{h_2}](x) &= {h_2}(x,t) \!+\! \int_0^x \! {{{\hat \alpha }_I}(x,\xi ,t)} {f_2}(\xi ,t)d\xi \\
						&+ \int_0^x {{{\hat \beta }_I}(x,\xi ,t)} {h_2}(\xi ,t)d\xi  ,
					\end{aligned}
				\end{align}
			\end{footnotesize}
		\end{subequations}
		where ${\Gamma _2}$ and $\Gamma _2^{ - 1}$ are defined as new conversion operators.
		Similar to \eqref{inverse kernel}, the approximate inversion kernel ${{\hat \alpha }_I}$, ${{\hat \alpha }_I}$ satisfies the following form 
		\begin{subequations}\label{approximate inverse kernel}
			\setlength{\abovedisplayskip}{2pt}
			\setlength{\belowdisplayskip}{2pt}
			\begin{small}
				\begin{align}
					{\hat{\alpha} _I}(x,\xi,t ) = \hat{\alpha} (x,\xi,t ) + \int_\xi ^x {\hat{\beta} (x,s,t){\hat{\alpha} _I}(s,\xi,t )ds} \\
					{\hat{\beta} _I}(x,\xi,t ) = \hat{\beta} (x,\xi,t ) + \int_\xi ^x {\hat{\beta} (x,s,t){\hat{\beta} _I}(s,\xi,t )ds} .
				\end{align}
			\end{small}
		\end{subequations}
		
		The error introduced by the approximate kernel adds further perturbation terms to the target system, making its composition more complex and posing challenges to theoretical stability analysis.
		Considering the property \eqref{kernel maps approximation} obtained from the universal approximation theorem, the error in the operator kernel can be confined to a small domain. This suggests that the approximation operator has a similar boundary as in Eqs. \eqref{kernel boundary}. Formulated as the following lemma.
		\begin{lemma}
			For all $(x,\xi ,t) \in T$ at \eqref{region T}, the approximate kernel is bounded and the approximate errors are square integrable.
			\begin{subequations}\label{approximate kernel boundary}
				\begin{small}
					\begin{align}
						&\left| {\tilde \alpha } \right| \le \bar {\tilde \alpha} , \quad \left| {{{\tilde \alpha }_I}} \right| \le {{\bar {\tilde \alpha} }_I}, \quad \left| {\hat \alpha } \right| \le \bar {\hat \alpha} , \quad \left| {{{\hat \alpha }_I}} \right| \le {{\bar {\hat \alpha} }_I}\\
						&\left| {\tilde \beta } \right| \le \bar {\tilde \beta} , \quad \left| {{{{\tilde \beta} }_I}} \right| \le {{\bar {\tilde \beta} }_I}, \quad \left| {\hat \beta } \right| \le \bar {\hat \beta} , \quad \left| {{{\hat \beta }_I}} \right| \le {{\bar {\hat \beta} }_I}\\
						&\left| {{\partial _t}\tilde \alpha } \right|, \quad \left| {{\partial _t}{{\tilde \alpha }_I}} \right|, \quad \left| {{\partial _t}\hat \alpha } \right|, \quad \left| {{\partial _t}{{\hat \alpha }_I}} \right| \in {L_2}\\
						&\left| {{\partial _t}\tilde \beta } \right|, \quad \left| {{\partial _t}{{\tilde \beta }_I}} \right|, \quad \left| {{\partial _t}\hat \beta } \right|, \quad \left| {{\partial _t}{{\hat \beta }_I}} \right| \in {L_2}.
					\end{align}
				\end{small}
			\end{subequations}
		\end{lemma}
		\begin{IEEEproof}
			Similar to the proof of Lemma \ref{operator boundary}.
		\end{IEEEproof}
		
		%	\begin{lemma}\label{lemma target system}
			The backstepping transformation \eqref{approximate backstepping transformation} and the approximate controller \eqref{approximate controller}, with the approximate backstepping kernels satisfying \eqref{hat kernel maps} and \eqref{approximate inverse kernel}, map the dynamics \eqref{dynamic system} into the new target system as follows (\textit{Here t is omitted for brevity}).
			\begin{small} 
				\begin{subequations}\label{target system2}
					\begin{align}
						&\begin{aligned}\label{target system2 a}
							&{\partial _t}{f_2}(x) + \lambda {\partial _x}{f_2}(x) = {n^T}(x)\dot {\hat M}\\
							&\ +{{\hat m}_1}{f_2}(x) + {{\hat m}_1}\hat e(x) + {{\hat m}_2}\hat \tau (x) + {{\hat m}_2}h_2(x)\\
							&\ +{{\hat m}_2}\int_0^x {{{\hat \alpha }_I}(x,\xi )} {f_2}(\xi )d\xi + {{\hat m}_2}\int_0^x {{{\hat \beta }_I}(x,\xi )} {h_2}(\xi )d\xi  ,
						\end{aligned}\\
						&\begin{aligned}\label{target system2 b}
							&{\partial _t}{h_2}(x) - \mu {\partial _x}{h_2}(x)\\
							&\ = [{{\hat m}_3}\hat e(x) + {{\hat m}_4}\hat \tau (x) + {{\hat m}_4}\Gamma _2^{ - 1}[{f_2},{h_2}](x)] + [{z^T}(x) \dot {\hat N}]\\
							&\ - \int_0^x {\left[ {{{\hat \alpha }_t}(x,\xi ) - \lambda {{\tilde \alpha }_\xi }(x,\xi ) + \mu {{\tilde \alpha }_x}(x,\xi ) + {{\hat m}_4}\alpha (x,\xi )} \right.} \\
							&\;\;\;\;\;\;\;\;\; \left. { - {{\hat m}_1}\tilde \alpha (x,\xi ) - {{\hat m}_3}\tilde \beta (x,\xi )} \right]{f_2}(\xi )d\xi  \\
							&\ - \int_0^x {\left[ {{{\hat \beta }_t}(x,\xi ) + \mu {{\tilde \beta }_\xi }(x,\xi ) + \mu {{\tilde \beta }_x}(x,\xi ) + {{\hat m}_4}\hat \beta (x,\xi )} \right.} \\
							&\;\;\;\;\;\;\;\;\; \left. { - {{\hat m}_2}\tilde \alpha (x,\xi)} \right]\Gamma _2^{ - 1}[{f_2},{h_2}](\xi)d\xi \\
							&\ - \int_0^x {\hat \alpha (x,\xi)} [{{\hat m}_1}\hat e(\xi) + {{\hat m}_2}\hat \tau (\xi)]d\xi \\
							&\ - \int_0^x {\hat \beta (x,\xi)} [{{\hat m}_3}\hat e(\xi) + {{\hat m}_4}\hat \tau (\xi)]d\xi \\
							&\ - \int_0^x {\hat \alpha (x,\xi)[} {n^T}(\xi)\dot {\hat M}]d\xi  - \int_0^x {\hat \beta (x,\xi)} [{z^T}(\xi)\dot {\hat N}]d\xi \\
							&\ - (\lambda  + \mu )\tilde \alpha (x,x){f_2}(x) - \lambda \hat \alpha (x,0)[q{h_2}(0) + q\hat \tau (0)] \\
							&\ + \mu \hat \beta (x,0){h_2}(0)  ,
						\end{aligned} \\
						& {f_2}(0) = q{h_2}(0) + q\hat \tau (0)  \label{target system2 c}, \\
						& {h_2}(1) = U - \hat U  \label{target system2 d},
					\end{align}
				\end{subequations}
			\end{small}
			where $n$ and $z$ represent the boundary filters, $\hat e$ and $\hat \tau$ denote the state estimation error as defined in \eqref{adaptive errors}. ${\dot {\hat M}}$ and ${\dot {\hat N}}$ signify the derivatives of the time-varying parameters of the system, satisfying
			$	\dot {\hat M}(t) = {[{{\dot {\hat m}}_1}(t),{{\dot {\hat m}}_2}(t)]^T},\dot {\hat N}(t) = {[{{\dot {\hat m}}_3}(t),{{\dot {\hat m}}_4}(t)]^T} $.
			%	\end{lemma}
		\begin{IEEEproof}
			The detailed derivation of the target system \eqref{target system2} obtained using the approximate backstepping kernel is given in the Appendix \ref{appendix B}.
		\end{IEEEproof}

		\subsection{Convergence analysis}
		By bounding the operator error, we can then treat the NN as a disturbance effect on the adaptive controller and analyze the resulting stability properties of the closed-loop system that using the approximate controller.
		% The resulting kernel error may affect the controller and system state, 
		Specifically, we establish Theorem \ref{converge} based on the Lyapunov stability principle, as demonstrated in proof \ref{converge proof}. Further details are provided in Appendix \ref{appendix C}.
		\begin{theorem}\label{converge}
			Consider the origin system \eqref{original_system} and the state estimates ${\hat \psi }$,${\hat \varphi }$ defined from \eqref{adaptive states} using the boundary filters \eqref{filter1}-\eqref{filter2} and the adaptive control law \eqref{adaptive laws}. Consider the approximated controller $\hat U(t) = \int_0^1 {\hat \alpha (1,\xi ,t)} \hat \psi (\xi ,t)d\xi  + \int_0^1 {\hat \beta (1,\xi ,t)} \hat \varphi (\xi ,t)d\xi$, where $\hat\alpha,\hat\beta$ are the approximated kernels \eqref{hat kernel maps} obtained from DeepONet. Then all signals within the closed-loop system exhibit boundedness and integrability in the L2-sense. Furthermore, $\psi (x, \cdot ),\varphi (x, \cdot ) \in {L_\infty } \cap {L_2}$ and $\psi (x, \cdot ),\varphi (x, \cdot ) \to 0$ for all $x \in [0,1]$.
		\end{theorem}
		
		\begin{IEEEproof} \label{converge proof}
			Consider the following Lyapunov candidate sub-functions as
			\begin{subequations} \label{Lyapunov sub-function}
				\begin{align}\label{V11}
					&{V_1} = \left\| {{n_1}} \right\|_{ - a}^2 \; {V_2} = \left\| {{n_2}} \right\|_{ - a}^2 \; {V_3} = \left\| {{n_3}} \right\|_{ - a}^2\\
					&{V_4} = \left\| {{z_1}} \right\|_b^2  \;\;\;\; {V_5} = \left\| {{z_2}} \right\|_b^2 \;\;\;\; {V_6} = \left\| {{z_3}} \right\|_b^2\\
					&{V_7} = \left\| {{f_2}} \right\|_{ - a}^2 \;\; {V_8} = \left\| {{h_2}} \right\|_b^2\\
					&{V_{9}} = \left\| e \right\|_{ - a}^2 \;\;\;\;  {V_{10}} = \left\| \tau  \right\|_b^2,
				\end{align}
			\end{subequations}
			where we set $a$ and $b$ as constants greater than one for ease of analysis below. Then, we consider Lyapunov candidate functions that are linear combinations of sub-functions
			\begin{equation}
				\setlength{\abovedisplayskip}{3pt}
				\setlength{\belowdisplayskip}{3pt}
				\begin{aligned}
					{V_{11}} &= {d_1}[{V_1} + {V_2}] + {d_2}{V_3} + {d_3}[{V_4} + {V_5}]\\
					&+ {d_4}{V_6} + {d_5}{V_7} + {d_6}{V_8} + {d_7}{V_9} + {d_8}{V_{10}},
				\end{aligned}
			\end{equation}
			where $d_i, i\!=\!1,...,8$ are positive constants. 
			From Appendix \ref{appendix C}, 
			the differentiation of all Lyapunov sub-functions is transformed into linear combinations of candidate functions. Their coefficient multipliers are transformed into constant ($c_i, \, i=1...10$) or bounded convergence ($l_i, \, i=1...14$) terms.
			Further integration leads to the following form
			\begin{equation} \label{V11 de}
				\setlength{\abovedisplayskip}{5pt}
				\setlength{\belowdisplayskip}{5pt}
				\begin{small}
					\begin{aligned}  
						&{{\dot V}_{11}} \le  - {d_1}(\lambda a - 1)[{V_1} + {V_2}] - {d_2}\lambda a{V_3}\\
						&- {d_3}(\mu b - 1)[{V_4} + {V_5}] - {d_4}\mu b{V_6}\\
						&- [{d_5}(\lambda a - {c_7}) - {d_1}{c_1}{e^a} - {d_3}{c_3}{e^{b + a}} - {d_4}{c_5}{e^{b + a}}]{V_7}\\
						&- [{d_6}(\mu b - {c_8}) - {d_1}{c_2} - {d_3}{c_4}{e^b} - {d_4}{c_6}{e^b} - 2{d_5}]{V_8}\\
						&- [{d_7}\lambda a - 4{d_1}{e^a} - 4{d_3}{e^{b + a}} - 2{d_5}{e^a} - {d_6}{c_9}{e^{b + a}}]{V_9}\\
						&- [{d_8}\mu b - 4{d_1} - 4{d_3}{e^b} - 2{d_5} - {d_6}{c_{10}}{e^b}]{V_{10}}\\
						&+ \left( {{d_1}{l_1} + {d_3}{l_4} + {d_5}{l_6} + {d_6}{l_9}} \right)[{V_1} + {V_2}] + {d_6}{l_{11}}{V_7}\\
						&+ ({d_1}{l_2} + {d_3}{l_5} + {d_5}{l_7} + {d_6}{l_{10}})[{V_4} + {V_5}] 
						+ {d_6}{l_{12}}{V_8}\\
						&- ( - 2{d_2}\lambda  - 2{d_5}\lambda {q^2} + {d_6}\mu )h_2^2(0)\\
						&+ (4{d_2}\lambda  + 4{d_5}\lambda {q^2} + 4{d_6}{e^b}){\tau ^2}(0)+ {d_6}{l_{14}}h_2^2(0)\\
						&- ({d_3}\mu  - {d_2}{l_3} - {d_5}{l_8} - {d_6}{l_{13}}){z^2}(0) ,
					\end{aligned}
				\end{small}
			\end{equation}
			the constant coefficient multipliers for these subsystems are set to
			\vspace{-1mm}
			\begin{equation} \label{d set}
				\begin{aligned}
					&{d_1} = {d_2} = {e^{ - a}}, \;\;\;\;\;\; {d_3}= {d_4} = {e^{ - b - a}}, \;\;\;\;\;\; {d_5} = 1,\\
					&{d_6} = \frac{{2{d_2}\lambda  + 2{d_5}\lambda {q^2}}}{\mu }, \;\;\;\;\;\; {d_7} = {e^{b + a}}, \;\;\;\;\;\; {d_8} = {e^b},
				\end{aligned}
				\vspace{-1mm}
			\end{equation}                                               
			and we choose
			\begin{subequations}
				\setlength{\abovedisplayskip}{5pt}
				\setlength{\belowdisplayskip}{5pt}
				\begin{footnotesize}
					\begin{align}
						& \begin{aligned}\label{b set}
							\setlength{\abovedisplayskip}{1pt}
							\setlength{\belowdisplayskip}{1pt}
							b \!> \! \max \! \left\{ {\!1,\! \frac{1}{\mu },\! \frac{{{e^{\! - \!a}}{(c_2+c_4+c_6)} \!\!+\! {d_6}{c_8}}}{{\mu {d_6}}},\!} \!\right.\left. \! {\frac{{8{e^{ \!-\! b \!-\! a}} \!\!+\! 2{e^{ \!-\! b}} \!\!+\! {d_6}{c_{10}}}}{\mu }} \!\! \right\} ,
						\end{aligned} \\
						& \begin{aligned}\label{a set}
							\setlength{\abovedisplayskip}{1pt}
							\setlength{\belowdisplayskip}{1pt}
							a \!>\! \max \left\{ {1,\frac{1}{\lambda },\frac{{{c_1} \!+\! {c_3} \!+\! {c_5} \!+\! {c_7}}}{\lambda },} \right.\left. {\frac{{8{e^{ - \delta  - \alpha }} \!+\! 2{e^{ - \alpha }} \!+\! {d_6}{h_9}}}{\lambda }} \right\},
						\end{aligned} 
					\end{align}
				\end{footnotesize}
			\end{subequations}
			then we can obtain
			\begin{small}
				\setlength{\abovedisplayskip}{5pt}
				\setlength{\belowdisplayskip}{5pt}
				\begin{equation}\label{V11 simp}
					\begin{aligned}
						{{\dot V}_{11}}\left( t \right) &\le  - C{V_{11}}\left( t \right) + {F}\left( t \right){V_{11}}\left( t \right) + {G}\left( t \right)\\
						&- \left( {{d_3}\mu  - {H}\left( t \right)} \right){z^2}(0).
					\end{aligned}
				\end{equation}
			\end{small}
			The above formula corresponds to equation \eqref{V11 de}, $C$ is a constant greater than 0, determined by the first eight coefficients of equation \eqref{V11 de}, ${F}\left( t \right)$ is integrable function corresponding to terms 9-12, and ${G}\left( t \right) = (4{d_2}\lambda  + 4{d_5}\lambda {q^2} + 4{d_6}{e^b}){\tau ^2}(0)+{d_6}{l_{14}}h_2^2(0)$ converges in finite time, ${H}\left( t \right) = {d_2}{l_3} + {d_5}{l_8} + {d_6}{l_{13}}$ are also integrable functions.
			Using the comparison principle yields that
			\vspace{-2mm}
			\begin{equation}\label{comparison principle}
				\begin{aligned}
					&{e^{Ct}}{V_{11}}\left( t \right) \le {V_{11}}\left( 0 \right){e^{\int_0^t {F(s)ds} }}\\
					&+ \int_0^t {{e^{Cs + \int_s^t {F(x)dx} }} [G(s) \!-\! \left( {{d_3}\mu  \!-\! H\left( s \right)} \right){z^2}(0, s)]ds} \\
					&\le{V_{11}}\left( 0 \right){e^{\int_0^t {F(s)ds}}} + {e^{\int_0^t {F(s)ds} }}\int_0^t {{e^{Cs}}G(s)ds} \\
					&+ {e^{\int_0^t {F(s)ds} }}\int_0^t {{e^{Cs}}\left( {{d_3}\mu  - H\left( s \right)} \right){z^2}(0,s)ds} .
				\end{aligned}
				\vspace{-2mm}
			\end{equation}
			Since ${F}\left( t \right)$ is integrable function, ${e^{\int_0^t {F(s)ds} }}$ is bounded, and ${G}\left( t \right)$ converges in finite time, then there $\exists {T_g},\forall t > T_g \to G(t) = 0$, implying that $\int_0^t {{e^{Cs}}G(s)ds} $ is also bounded. Additionally, ${H\left( t \right)}$ converges to 0 in finite time, so that $\exists {T_H}$, ${d_3}\mu  - H\left( t \right) > 0$ and ${d_3}\mu  - H\left( t \right) \to {d_3}\mu$ while $t > {T_H}$. This implies that $\int_0^t {{e^{Cs}}\left( {{d_3}\mu  - H\left( s \right)} \right){z^2}(0,s)ds}$ must also be bounded. Hence, ${e^{Ct}}{V_{11}}\left( t \right)$ is bounded as $t \to \infty$, which means that ${e^{Ct}}{V_{11}}\left( t \right)$ converges to zero for finite time.
			
			Now therefore, $\left\| n \right\|,\left\| z \right\|,\left\| {{f_2}} \right\|,\left\| {{h_2}} \right\| \in {L_2} \cap {L_\infty }$ are established. Then $\left\| \psi  \right\|,\left\| \varphi  \right\| \in {L_2} \cap {L_\infty }$ can be obtained directly by formulas \eqref{real status boundary} and \eqref{estimated error boundary}, as $\psi (x, \cdot ),\varphi (x, \cdot ) \to 0$ for all $x \in [0,1]$.
			
		\end{IEEEproof}

		\section{Simulation}

		\section{Conclusion}

		\appendix{
			
		\subsection{Proof of target system \eqref{target system2}}\label{appendix B}
		For the equation \eqref{target system2 a}, 
		by substituting equation \eqref{approximate backstepping transformation a} into the dynamic system \eqref{dynamic system a}, then use the error system \eqref{adaptive errors} and replacing $\hat \varphi (x,t)$ with equation \eqref{approximate inverse transformation b}, we obtain
		\begin{equation}
			\begin{small}
				\begin{aligned}
					&{\partial _t}{f_2}(x) + \lambda {\partial _x}{f_2}(x) \\
					&= {{\hat m}_1}\psi (x) + {{\hat m}_2}\varphi (x) + {n^T}(x)\dot {\hat M} \\
					&= {n^T}(x)\dot {\hat M} \!+\! {{\hat m}_1}[\hat \psi (x) \!+\! \hat e(x)] \!+\! {{\hat m}_2}[\hat \varphi (x) \!+\! \hat \tau (x)]\\
					&={n^T}(x)\dot {\hat M}\!+\!{{\hat m}_1}{f_2}(x) \!+\! {{\hat m}_1}\hat e(x) \!+\! {{\hat m}_2}\hat \tau (x) \!+\! {{\hat m}_2}h_2(x)\\
					&\;\;\;\ +{{\hat m}_2} \! \int_0^x \! {{{\hat \alpha }_I}(x,\xi )} {f_2}(\xi )d\xi \!+\! {{\hat m}_2} \! \int_0^x \! {{{\hat \beta }_I}(x,\xi )} {h_2}(\xi )d\xi ,
				\end{aligned}
			\end{small}
		\end{equation}
		
		For the equation \eqref{target system2 b}, 
		after applying the shift transformation to equation \eqref{approximate backstepping transformation b}, $\hat \psi (x,t)$ is differentiated with respect to time $t$, followed by the replacement of ${{\partial _t}\hat \psi }$ and ${{\partial _t}\hat \varphi }$ with the dynamic equations \eqref{dynamic system a}-\eqref{dynamic system b}, then integrating the distribution, combined with \eqref{kernel error}, we obtain
		\begin{equation}
			\begin{small}
				\begin{aligned}
					&{\partial _t}\hat \varphi (x)= {\partial _t}{h_2}(x) \\
					&+ \int_0^x {{{\hat \alpha }_t}(x,\xi )\hat \psi (\xi )d\xi }  + \int_0^x {{{\hat \beta }_t}(x,\xi )\hat \varphi (\xi )d\xi } \\
					&+ \int_0^x {\hat \alpha (x,\xi )[{n^T}(\xi )\dot {\hat M}]d\xi }  + \int_0^x {\hat \beta (x,\xi )[{z^T}(\xi )\dot {\hat N}]d\xi } \\
					&+ \int_0^x {\hat \alpha (x,\xi )[{{\hat m}_1}\psi (\xi ) + {{\hat m}_2}\varphi (\xi )]d\xi } \\
					&+ \int_0^x {\hat \beta (x,\xi )[{{\hat m}_3}\psi (\xi ) + {{\hat m}_4}\varphi (\xi )]d\xi } \\
					&+ \lambda \int_0^x {[{\alpha _\xi }(x,\xi ) - {{\tilde \alpha }_\xi }(x,\xi )]\hat \psi (\xi )d\xi }  \\
					&- \mu \int_0^x {[{\beta _\xi }(x,\xi ) - {{\tilde \beta }_\xi }(x,\xi )]\hat \varphi (\xi )d\xi } \\
					&- \lambda [\alpha (x,x) - \tilde \alpha (x,x)]\hat \psi (x) + \lambda [\alpha (x,0) - \tilde \alpha (x,0)]\hat \psi (0)\\
					&+ \mu [\beta (x,x) - \tilde \beta (x,x)]\hat \varphi (x) - \mu [\beta (x,0) - \tilde \beta (x,0)]\hat \varphi (0),
				\end{aligned}
			\end{small}
		\end{equation}
		Similarly, differentiating $\hat \psi (x,t)$ with respect to space $x$, we obtain
		\begin{equation}
			\hspace{-15mm}
			\begin{aligned}
				&	{\partial _x}\hat \varphi (x) = {\partial _x}h(x)\\
				&	+ \int_0^x {[{\alpha _x}(x,\xi ) - {{\tilde \alpha }_x}(x,\xi )]\hat \psi (\xi )d\xi } \\
				&	+ \int_0^x {[{\beta _x}(x,\xi ) - {{\tilde \beta }_x}(x,\xi )]\hat \varphi (\xi )d\xi } \\
				&	+ \alpha (x,x)\hat \psi (x) - \tilde \alpha (x,x)\hat \psi (x)\\
				&	+ \beta (x,x)\hat \varphi (x) - \tilde \beta (x,x)\hat \varphi (x),
			\end{aligned}
		\end{equation}
		Substituting the aforementioned two formulas into \eqref{dynamic system b}, and utilizing the kernel equation set \eqref{kernel equation} to eliminate the middle term, and subsequently combining and sorting, we obtain
		\begin{equation}
			\begin{aligned}
				&{{\hat m}_3}\psi (x,t) + {{\hat m}_4}\varphi (x,t) + {z^T}(\xi )\dot {\hat N} \\
				&= {\partial _t}h(x) - \mu {\partial _x}h(x)\\
				&+ \int_0^x {[{{\hat \alpha }_t}(x,\xi ) - \lambda {{\tilde \alpha }_\xi }(x,\xi ) + \mu {{\tilde \alpha }_x}(x,\xi ) + {{\hat m}_4}\alpha (x,\xi )} \\
				&\;\;\;\;\; - {{\hat m}_1}\tilde \alpha (x,\xi ) - {{\hat m}_3}\tilde \beta (x,\xi )]\hat \psi (\xi )d\xi \\
				&+ \int_0^x {[{{\hat \beta }_t}(x,\xi ) + \mu {{\tilde \beta }_\xi }(x,\xi ) + \mu {{\tilde \beta }_x}(x,\xi ) - {{\hat m}_2}\tilde \alpha (x,\xi )} \\
				&\;\;\;\;\; + {{\hat m}_4}\hat \beta (x,\xi )]\hat \varphi (\xi )d\xi \\
				&+ \int_0^x {\hat \alpha (x,\xi )} [{{\hat m}_1}\hat e(\xi ) + {{\hat m}_2}\hat \tau (\xi )]d\xi \\
				&+ \int_0^x {\hat \beta (x,\xi )} [{{\hat m}_3}\hat e(\xi ) + {{\hat m}_4}\hat \tau (\xi )]d\xi \\
				&+ [\int_0^x {\hat \alpha (x,\xi )} [{n^T}(\xi )\dot {\hat M}]d\xi  + \int_0^x {\hat \beta (x,\xi )} [{z^T}(\xi )\dot {\hat N}]\\
				&+ (\lambda  + \mu )\tilde \alpha (x,x)\hat \psi (x) + \hat \psi (x){{\hat m}_3}\\
				&+ \lambda \hat \alpha (x,0)\hat \psi (0) - \mu \hat \beta (x,0)\hat \varphi (0)
			\end{aligned}
		\end{equation}
		Form \eqref{approximate inverse transformation}, \eqref{adaptive errors}, we get formula \eqref{target system2 b}.
		
		For the equation \eqref{target system2 c}, with \eqref{approximate backstepping transformation} and \eqref{dynamic system c}, we obtain
		\begin{equation}
			\begin{aligned}
				{f_2}(0) &= \hat \psi (0) = q\varphi (0) = q\hat \varphi (0) + q\hat \tau (0)\\
				&= q{h_2}(0) + q\hat \tau (0)
			\end{aligned}
		\end{equation}
		
		For the equation \eqref{target system2 d}, from \eqref{approximate backstepping transformation b}, \eqref{dynamic system d} and \eqref{approximate controller}, we obtain
		\begin{equation}
			\begin{aligned}
				{h_2}(1) = \hat \varphi (1) - \hat U = U - \hat U
			\end{aligned}
		\end{equation}

		\subsection{Proof of Theorem \ref{converge}}\label{appendix C}
		Firstly, we introduce some important properties to facilitate analysis.
		According to Cauchy–Schwarz's inequality, following inequalities which hold provided $\delta  \ge 1$
		\begin{subequations} \label{cauchy}
			\begin{align}
				&{I_\delta }\left[ {u\psi v} \right] \le \left\| u \right\|k \left\| \psi  \right\|_\delta ^2 + \frac{{\left\| u \right\|}}{k}\left\| v \right\|_\delta ^2 \label{cauchy a} \\
				&{I_\delta }\left[\! {\psi (x) \! \int_0^x \! {u(x, \! \xi )v(\xi )d\xi } } \! \right] \! \le \! \left\| u \right\| \! k \!  \left\| \psi  \right\|_\delta ^2 \!+\! \frac{{\left\| u \right\|}}{k} \! {e^\delta }{\left\| v \right\|^2} \label{cauchy b} \\
				&{I_{ - \delta }}\left[ {\psi (x)\int_0^x {u(x,\xi )v(\xi )d\xi } } \right] \le \bar uk  \left\| \psi  \right\|_{ - \delta }^2 + \frac{{\bar u}}{k}{\left\| v \right\|^2}, \label{cauchy c}
			\end{align}
		\end{subequations}
		From \eqref{approximate kernel boundary} and \eqref{approximate backstepping transformation}, since the kernels
		$(\hat \alpha ,\hat \beta )$ and $({{\hat \alpha }_I},{{\hat \beta }_I})$ are bounded for every $t$, for some signals
		${u_1}(x),{u_2}(x),{u_3}(x)$ defined on ${T = \{ 0 \le x \le 1,t \ge 0\}}$ and satisfying
		${u_3}(x) = {\Gamma _2}[{u_1},{u_2}](x)$ and ${u_2}(x) = \Gamma _2^{ - 1}[{u_1},{u_3}](x)$, there exist positive constants $M1,M2,N1,N2$ such that
		\begin{subequations} 
			\begin{align} 
				&{u_3}= {\Gamma _2}[{u_1},{u_2}] \;\; \to \left\| {{u_3}} \right\| \le {M_1}\left\| {{u_1}} \right\| + {M_2}\left\| {{u_2}} \right\|\\
				&{u_2}= {\Gamma_2^{-1} }[{u_1},{u_3}] \to\left\| {{u_1}} \right\| \le {N_1}\left\| {{u_1}} \right\| + {N_2}\left\| {{u_3}} \right\|,
			\end{align}
		\end{subequations}
		Then from \eqref{adaptive errors} and \eqref{approximate backstepping transformation}, we obtain
		\begin{subequations}\label{real status boundary}
			\begin{small}
				\begin{align}
					&\left\| \psi  \right\| = \left\| {\hat \psi  + \hat e} \right\| \le \left\| {\hat \psi } \right\| + \left\| {\hat e} \right\| = \left\| {{f_2}} \right\| + \left\| {\hat e} \right\| \label{real status boundary a} \\
					&\begin{aligned} \label{real status boundary b}
						\left\| \varphi  \right\| &= \left\| {\hat \varphi  + \hat \tau } \right\| \le \left\| {\hat \varphi } \right\| + \left\| {\hat \tau } \right\| = \left\| {\Gamma _2^{ - 1}[{f_2},{h_2}]} \right\| + \left\| {\hat \tau } \right\| \\
						&\le {N_1}\left\| {{f_2}} \right\| + {N_2}\left\| {{h_2}} \right\| + \left\| {\hat \tau } \right\| .
					\end{aligned}
				\end{align}
			\end{small}
		\end{subequations}
		By equations \eqref{non-adaptive states}, \eqref{non-adaptive errors}, \eqref{adaptive errors} and \eqref{adaptive states}, the estimated errors satisfy
		\begin{subequations}\label{estimated error boundary}
			\begin{small}
				\begin{align}
					&\left\| {\hat e} \right\| = \left\| {e + \bar \psi  - \hat \psi } \right\| = \left\| {e + {n^T}\tilde M} \right\| \le \left\| e \right\| + \left| {\tilde M} \right|\left\| n \right\| \label{estimated error boundary a} \\
					&\left\| {\hat \tau } \right\| = \left\| {\tau  + \bar \varphi  - \hat \varphi } \right\| = \left\| {\tau  + {z^T}\tilde N} \right\| \le \left\| \tau  \right\| + \left| {\tilde N} \right|\left\| z \right\|. \label{estimated error boundary b} 
				\end{align}
			\end{small}
		\end{subequations}
		After clarifying these properties, we proceed to analyze the Lyapunov sub-function candidates defined in Eqs \eqref{Lyapunov sub-function}.
		
		\subsubsection{${V_1}=\left\| {{n_1}} \right\|_{ - a}^2$ and ${V_2} = \left\| {{n_2}} \right\|_{ - a}^2$} 
		From filter \eqref{filter1 a}, equation \eqref{zn} and \eqref{MN}, employing the properties \eqref{property 1}, \eqref{real status boundary a} and \eqref{estimated error boundary a}, we have
		\begin{equation}\label{V1}
			\begin{small}
				\begin{aligned}
					&{{\dot V}_1} = 2\int_0^1 {{e^{ - ax}}n_1^T{\partial _t}{n_1}} dx\\
					&=  - 2\lambda {I_{ - a}}[n_1^T{\partial _x}{n_1}] + 2{I_{ - \delta }}[{n_1}\psi ]\\
					&\le  - \lambda (n_1^2(1){e^{ - a}} \!-\! n_1^2(0) + a\left\| {{n_1}} \right\|_{ - a}^2) + \left\| {{n_1}} \right\|_{ - a}^2 + \left\| \psi  \right\|_{ - a}^2\\
					&\le (1 - \lambda a)\left\| {{n_1}} \right\|_{ - a}^2 - \lambda n_1^2(1){e^{ - \delta }} + 2{\left\| {\hat \psi } \right\|^2} + 2{\left\| {\hat e} \right\|^2}\\
					&\le (1 - \lambda a)\left\| {{n_1}} \right\|_{ - a}^2 + 2{\left\| {{f_2}} \right\|^2} + 2{(\left\| e \right\| + {n^T}\left| {\tilde M} \right|)^2}\\
					&\le (1 - \lambda a)\left\| {{n_1}} \right\|_{ - a}^2 + 2{\left\| {{f_2}} \right\|^2} + 4{\left\| e \right\|^2} + 4{\left| {\tilde M} \right|^2}{\left\| n \right\|^2}.
				\end{aligned}
			\end{small}
		\end{equation}
		
		By filter \eqref{filter1 b}, equation \eqref{zn} and \eqref{MN}, employing the properties \eqref{property 1}, \eqref{real status boundary b} and \eqref{estimated error boundary b}, we obtain
		\begin{equation}\label{V2}
			\begin{footnotesize}
				\begin{aligned}
					&{{\dot V}_2} = 2\int_0^1 {{e^{ - ax}}n_2^T{\partial _t}{n_2}} dx\\
					&=  - 2\lambda {I_{ - a}}[n_2^T{\partial _x}{n_2}] + 2{I_{ - a}}[{n_2}\varphi ]\\
					&\le  - \lambda (n_2^2(1){e^{ - a}} \!-\! n_2^2(0) + a\left\| {{n_2}} \right\|_{ - a}^2) + \left\| {{n_2}} \right\|_{ - a}^2 + \left\| \varphi  \right\|_{ - a}^2\\
					&\le \! (1 \!-\! \lambda a)\left\| {{n_2}} \right\|_{ - a}^2 \!-\! \lambda n_2^2(1){e^{ - \delta }} \!+\! 2{\left\| {\hat \varphi } \right\|^2} \!+\! 2{\left\| {\hat \tau } \right\|^2}\\
					&\le \! (1 \!-\! \lambda a)\left\| {{n_2}} \right\|_{ - a}^2 \!+\! 2{({N_1}\left\| {{f_2}} \right\| \!+\! {N_2}\left\| {{h_2}} \right\|)^2} \!+\! 2{(\left\| \tau  \right\| \!+\! {z^T} \! \left| {\tilde N} \right|)^2}\\
					&\le \! (1 \!-\! \lambda a) \! \left\| {{n_2}} \right\|_{ - a}^2 \!\!+\! 4N_1^2{\left\| {{f_2}} \right\|^2} \!\!+\! 4N_2^2{\left\| {{h_2}} \right\|^2} \!\!+\! 4{\left\| \tau  \right\|^2} \!\!+\! 4{\left| {\tilde N} \right|^2} \!\! {\left\| z \right\|^2}
				\end{aligned}\
			\end{footnotesize}
		\end{equation}
		
		Since the forms  of $V_1$ and $V_2$ are essentially the same, we combine the stability analysis accordingly, and utilizing property \eqref{property 2} to obtain:
		\begin{equation}\label{V1 V2}
			\begin{small}
				\begin{aligned}
					&{{\dot V}_1} + {{\dot V}_2}\\
					&\le (1 - \lambda a)\left\| n \right\|_{ - a}^2 + 4{\left| {\tilde M} \right|^2}{\left\| n \right\|^2} + 4{\left| {\tilde N} \right|^2}{\left\| z \right\|^2}\\
					&+ (2 + 4N_1^2){\left\| {{f_2}} \right\|^2} + 4N_2^2{\left\| {{h_2}} \right\|^2} + 4{\left\| e \right\|^2} + 4{\left\| \tau  \right\|^2}\\
					&\le (1 \!-\! \lambda a)[{V_1} \!+\! {V_2}] \!+\! 4{e^a }{\left| {\tilde M} \right|^2}[{V_1} \!+\! {V_2}] \!+\! 4{\left| {\tilde N} \right|^2}[{V_3} \!+\! {V_4}]\\
					&+ (2 + 4N_1^2){e^a }{V_7} + 4N_2^2{V_8} + 4{e^a }{V_9} + 4{V_{10}}\\
					&\le (1 - \lambda a)[{V_1} + {V_2}] + {l_1}[{V_1} + {V_2}] + {l_2}[{V_3} + {V_4}]\\
					&+ {c_1}{e^a }{V_7} + {c_2}{V_8} + 4{e^a }{V_9} + 4{V_{10}}.
				\end{aligned}
			\end{small}
		\end{equation}
		where 
		\begin{equation}\label{V1 V2 l}
			{l_1} = 4{e^a}{\left| {\tilde M} \right|^2} \;\;\;\; {l_2} = 4{\left| {\tilde N} \right|^2}
		\end{equation}
		are integrable functions and 
		\begin{equation}\label{V1 V2 c}
			{c_1} = (2 + 4N_1^2)  \;\;\;\; {c_2} = 4N_2^2
		\end{equation}
		are two positive constants.
		
		\subsubsection{${V_3} = \left\| {{n_3}} \right\|_{ - a}^2$}
		From filter \eqref{filter1 c} and equation \eqref{MN}, with the properties \eqref{property 1}, \eqref{real status boundary b} and \eqref{estimated error boundary a}, we have
		\begin{equation}\label{V3}
			\begin{small}
				\begin{aligned}
					&{{\dot V}_3} =  - 2\lambda {I_{ - a}}[n_3^T{\partial _x}{n_3}]\\
					&=  - \lambda (n_3^2(1){e^{ - a}} - n_3^2(0) + a\left\| {{n_3}} \right\|_{ - a}^2)\\
					&\le  - \lambda a\left\| {{n_3}} \right\|_{ - a}^2 + \lambda {\varphi ^2}(0)\\
					&\le  - \lambda a\left\| {{n_3}} \right\|_{ - a}^2 + 2\lambda {h_2}^2(0) + 2\lambda {{\hat \tau }^2}(0)\\
					&\le  - \lambda a\left\| {{n_3}} \right\|_{ - a}^2 + 2\lambda {h_2}^2(0) + 2\lambda {\left\| {\tau (0) + {z^T}(0)\tilde N} \right\|^2}\\
					&\le  - \lambda a\left\| {{n_3}} \right\|_{ - a}^2 + 2\lambda {h_2}^2(0) + 4\lambda {\tau ^2}(0) + 4\lambda {\left| {\tilde N} \right|^2}{\left| {\sigma (0)} \right|^2}\\
					& \le  - \lambda a{V_3} + {l_3}{\left| {\sigma (0)} \right|^2} + 2\lambda {h^2}(0) + 4\lambda {\tau ^2}(0).
				\end{aligned}
			\end{small}
		\end{equation}
		where 
		\begin{equation}\label{V3 l}
			{l_3} = 4\lambda {\left| {\tilde N} \right|^2}
		\end{equation}
		denotes functions that can be integrated.\\
		
		\subsubsection{${V_4} = \left\| {{z_1}} \right\|_b^2$ and ${V_5} = \left\| {{z_2}} \right\|_b^2$} 
		From filter \eqref{filter2 a}, equation \eqref{zn} and \eqref{MN}, employing the properties \eqref{property 1}, \eqref{real status boundary a} and \eqref{estimated error boundary a}, we have
		\begin{equation}\label{V4}
			\begin{small}
				\begin{aligned}
					&{{\dot V}_4} = 2\int_0^1 {{e^{bx}}z_1^T{\partial _t}{z_1}} dx = 2\mu {I_b}[z_1^T{\partial _x}{z_1}] + 2{I_b}[{z_1}\psi ]\\
					&\le \mu (z_1^2(1){e^b} - z_1^2(0) - b\left\| {{z_1}} \right\|_b^2) + \left\| {{z_1}} \right\|_b^2 + \left\| \psi  \right\|_b^2\\
					&\le (1 - \mu b)\left\| {{z_1}} \right\|_b^2 - \mu z_1^2(0) + {e^b}\left\| \psi  \right\|_b^2\\
					&\le (1 - \mu b)\left\| {{z_1}} \right\|_b^2 - \mu z_1^2(0) + 2{e^b}{\left\| {{f_2}} \right\|^2} + 2{e^b}{\left\| {\hat e} \right\|^2}\\
					&\le \!(1\! -\! \mu b)\!\! \left\| {{z_1}} \right\|_b^2 \!\! -\! \mu z_1^2(0) \!+\! 2{e^b}\!{\left\| {{f_2}} \right\|^2} \!\!+\! 4{e^b} \! {\left\| e \right\|^2} \!\!+\! 4{e^b} \! {\left| {\tilde M} \right|^2} \! \!{\left\|  n \right\|^2}.
				\end{aligned}
			\end{small}
		\end{equation}
		
		By filter \eqref{filter2 b}, equation \eqref{zn} and \eqref{MN}, employing the properties \eqref{property 1}, \eqref{real status boundary b} and \eqref{estimated error boundary b}, we obtain
		\begin{equation}\label{V5}
			\begin{aligned}
				&{{\dot V}_5} = 2\mu {I_b}[z_2^T{\partial _x}{z_2}] + 2{I_b}[{z_2}\varphi ]\\
				&\le \mu (z_2^2(1){e^b} - z_2^2(0) - b\left\| {{z_2}} \right\|_b^2) + \left\| {{z_2}} \right\|_b^2 + \left\| \varphi  \right\|_b^2\\
				&\le (1 - \mu b)\left\| {{z_2}} \right\|_b ^2 - \mu z_2^2(0) + 2{e^b}{\left\| {\hat \varphi } \right\|^2} + 2{e^b}{\left\| {\hat \tau } \right\|^2}\\
				&\le (1 - \mu b)\left\| {{z_2}} \right\|_b ^2 - \mu z_2^2(0) + 4{e^b}N_1^2{\left\| {{f_2}} \right\|^2}\\
				&\;\; + 4{e^\alpha }N_2^2{\left\| {{h_2}} \right\|^2} + 4{e^b}{\left\| \tau  \right\|^2} + 4{e^b}{\left| {\tilde N} \right|^2}{\left\| z \right\|^2}.
			\end{aligned}\
		\end{equation}
		
		Since the forms  of $V_4$ and $V_5$ are essentially the same, we combine the stability analysis accordingly, and utilizing property \eqref{property 2} to obtain:
		\begin{equation}\label{V4 V5}
			\begin{small}
				\begin{aligned}
					&{{\dot V}_4} + {{\dot V}_5} \\
					&\le (1 - \mu b)\left\| z \right\|_b^2 + 4{\left| {\tilde M} \right|^2}{e^b}{\left\| n \right\|^2}\\
					& + 4{\left| {\tilde N} \right|^2}{e^b}{\left\| z \right\|^2} + (4N_1^2{e^b} + 2{e^b}){\left\| {{f_2}} \right\|^2}\\
					& + 4N_2^2{e^b }{\left\| {{h_2}} \right\|^2} + 4{e^b}{\left\| e \right\|^2} + 4{e^b}{\left\| \tau  \right\|^2} - \mu {z^2}(0)\\
					&\le (1 - \mu b)[{V_4} + {V_5}] + {l_4}[{V_1} + {V_2}]\\
					& + {l_5}[{V_4} + {V_5}] + {c_3}{e^{b + a}}{V_7} + {c_4}{e^{b}}{V_8}\\
					& + 4{e^{b + a}}{V_9} + 4{e^b}{V_{10}} - \mu {z^2}(0).
				\end{aligned}
			\end{small}
		\end{equation}
		where 
		\begin{equation}\label{V4 V5 l}
			{l_4} = 4{\left| {\tilde M} \right|^2}{e^{b + a}} \;\;\;\; {l_5} = 4{\left| {\tilde N} \right|^2}{e^b}
		\end{equation}
		are integrable functions and 
		\begin{equation}\label{V4 V5 c}
			{c_3} = 2 + 4N_1^2 \;\;\;\; {c_4} = 4N_2^2
		\end{equation}
		are positive constants.
		
		\subsubsection{${V_6} = \left\| {{z_3}} \right\|_b^2$}
		Using \eqref{property 1} and filter \eqref{filter2 c}, we can get
		\begin{equation}
			\begin{small}
				\begin{aligned}
					{{\dot V}_6} &= 2\mu {I_b}[z_3^T{\partial _x}{z_3}] = \mu (z_3^2(1){e^b } - z_3^2(0) - b\left\| {{z_3}} \right\|_b^2)\\
					&=  - \mu b\left\| {{z _3}} \right\|_b^2 - \mu z_3^2(0) + \mu {e^b}{{\hat U}^2},
				\end{aligned}
			\end{small}
		\end{equation}
		Considering formula \eqref{approximate controller}, \eqref{approximate kernel boundary} and \eqref{real status boundary b}, we analyze ${\hat U}^2$
		to obtain:
		\begin{equation}\label{hat U}
			\begin{small}
				\begin{aligned}
					&{{\hat U}^2} = {\left( {\int_0^1 {\hat \alpha (1,\xi )\hat \psi (\xi )d\xi }  + \int_0^1 {\hat \beta (1,\xi )\hat \varphi (\xi )d\xi } } \right)^2}\\
					&\le 2\int_0^1 {{{\hat \alpha }^2}(1,\xi )} {{\hat \psi }^2}(\xi )d\xi  + 2\int_0^1 {{{\hat \beta }^2}(1,\xi )} {{\hat \varphi }^2}(\xi )d\xi \\
					&\le 2{{\bar {\hat \alpha} }^2}{\left\| {\hat \psi } \right\|^2} + 2{{\bar {\hat \beta} }^2}{\left\| {\hat \varphi } \right\|^2}\\
					&\le 2{{\bar {\hat \alpha} }^2}{\left\| {{f_2}} \right\|^2} + 2{{\bar {\hat \beta} }^2}{({N_1}\left\| {{f_2}} \right\| + {N_2}\left\| {{h_2}} \right\|)^2}\\
					&\le 2{{\bar {\hat \alpha} }^2}{\left\| {{f_2}} \right\|^2} + 4{{\bar {\hat \beta} }^2}N_1^2{\left\| {{f_2}} \right\|^2} + 4{{\bar \hat \beta }^2}N_2^2{\left\| {{h_2}} \right\|^2},
				\end{aligned}
			\end{small}
		\end{equation}
		Hence, by utilizing property \eqref{property 2}, we can obtain
		\begin{equation}\label{V6}
			\begin{aligned}
				&{{\dot V}_6} =  - \mu b\left\| {{z _3}} \right\|_b^2 - \mu z_3^2(0) + \mu {e^b}{{\hat U}^2}\\
				&\le  - \mu b \left\| {{z_3}} \right\|_b^2 - \mu z_3^2(0) \\
				&+ 2\mu {e^b}({{\bar {\hat \alpha} }^2} + 2{{\bar {\hat \beta} }^2}N_1^2){\left\| {{f_2}} \right\|^2} + 4\mu {e^b}{{\bar {\hat \beta} }^2}N_2^2{\left\| {{h_2}} \right\|^2}\\
				&\le  - \mu b{V_6} + {c_5}{e^{b + a}}{V_7} + {c_6}{e^b}{V_8} ,
			\end{aligned}
		\end{equation}
		where 
		\begin{equation}\label{V6 c}
			{c_5} = 2\mu ({{\bar {\hat \alpha} }^2} + 2{{\bar {\hat \beta} }^2}N_1^2)\;\;\;\;
			{c_6} = 4\mu {{\bar {\hat \beta} }^2}N_2^2 ,
		\end{equation}
		are two positive constants.
		
		\subsubsection{${V_7} = \left\| {{f_2}} \right\|_{ - a}^2$}
		From dynamic \eqref{target system2 a}, we get
		\begin{equation}
			\begin{footnotesize}
				\begin{aligned}
					&{{\dot V}_7} = 2{I_{ - a}}[{f_2}{\partial _t}{f_2}]\\
					&= 2{I_{ - a}}[{f_2}( - \lambda {\partial _x}{f_2} + {{\hat m}_1}{f_2} + {{\hat m}_1}\hat e + {{\hat m}_2}{h_2} + {{\hat m}_2}\hat \tau  + {n^T}\dot {\hat M}\\
					&+ {{\hat m}_2}\int_0^x {{{\hat \alpha }_I}(x,\xi )} {f_2}(\xi )d\xi  + {{\hat m}_2}\int_0^x {{{\hat \beta }_I}(x,\xi )} {h_2}(\xi )d\xi )]\\
					&=  - \lambda a\left\| {{f_2}} \right\|_{ - a}^2 - \lambda {f_2}^2(1){e^{ - a}} + \lambda {f_2}^2(0)\\
					&+ 2{{\hat m}_1}{I_{ - a}}[{f_2}{f_2} + {f_2}\hat e] + 2{{\hat m}_2}{I_{ - a}}[{f_2}{h_2} + {f_2}\hat \tau ]\\
					&+ 2{I_{ - a}}[{f_2}({n^T}\dot \hat M)] + 2{{\hat m}_2}{I_{ - a}}[{f_2}(\int_0^x {{{\hat \alpha }_I}(x,\xi )} {f_2}(\xi )d\xi )]\\
					&+ 2{{\hat m}_2}{I_{ - a}}[{f_2}(\int_0^x {{{\hat \beta }_I}(x,\xi )} {h_2}(\xi )d\xi )]
				\end{aligned}
			\end{footnotesize}
		\end{equation}
		With \eqref{target system2 c}, inequality \eqref{cauchy a}, \eqref{cauchy c} and \eqref{estimated error boundary}, boundary condition \eqref{hat parameter bounds} and \eqref{approximate kernel boundary}, we find
		\begin{equation} \label{V7}
			\begin{footnotesize}
				\begin{aligned}
					&{{\dot V}_7} \le  - \lambda a\left\| {{f_2}} \right\|_{ - a}^2 + 2\lambda {q^2}{h_2}^2(0) + 2\lambda {q^2}{{\hat \tau }^2}(0)\\
					&+ \bar m_1^2\left\| {{f_2}} \right\|_{ - a}^2 + \left\| {{f_2}} \right\|_{ - a}^2 + \bar m_1^2\left\| {{f_2}} \right\|_{ - a}^2 + \left\| {\hat e} \right\|_{ - a}^2 + \bar m_2^2\left\| {{f_2}} \right\|_{ - a}^2\\
					&+ \left\| {{h_2}} \right\|_{ - a}^2 + \bar m_2^2\left\| {{f_2}} \right\|_{ - a}^2 + \left\| {\hat \tau } \right\|_{ - a}^2 + \left\| {{f_2}} \right\|_{ - a}^2 + \left\| {{n^T}\dot {\hat M}} \right\|_{ - a}^2\\
					&+ \left\| {{f_2}} \right\|_{ - a}^2 + \bar m_2^2{{\bar {\hat \alpha} }_I}^2\left\| {{f_2}} \right\|_{ - a}^2 + \left\| {{h_2}} \right\|_{ - a}^2 + \bar m_2^2{{\bar {\hat \beta} }_I}^2\left\| {{f_2}} \right\|_{ - a}^2\\
					&\le ( - \lambda a + 2\bar m_1^2 + 2\bar m_2^2 + \bar m_2^2{{\bar {\hat \alpha} }_I}^2 + \bar m_2^2{{\bar {\hat \beta} }_I}^2 + 3)\left\| {{f_2}} \right\|_{ - a}^2\\
					&+ 2\lambda {q^2}{h_2}^2(0) + 4\lambda {q^2}({\tau ^2}(0) + {\left| {\tilde N} \right|^2}{z^2}(0))\\
					&+ 2({\left\| e \right\|^2} + {\left| {\tilde M} \right|^2}{\left\| n  \right\|^2}) + 2\left\| h \right\|_{ - a}^2\\
					&+ 2({\left\| \tau  \right\|^2} + {\left| {\tilde N} \right|^2}{\left\| z \right\|^2}) + {\left| {\dot {\hat M}} \right|^2}{\left\| n \right\|^2}\\
					&\le  - (\lambda a - {c_7}){V_7} + {l_6}[{V_1} + {V_2}] + {l_7}[{V_4} + {V_5}]+ 2{V_8} + 2{e^\delta }{V_9} \\
					&+ 2{V_{10}}+ {l_8}{z^2}(0) + 2\lambda {q^2}{h_2}^2(0) + 4\lambda {q^2}{\tau ^2}(0)		
				\end{aligned}
			\end{footnotesize}
		\end{equation}
		where 
		\begin{equation}\label{V7 l}
			{l_6} = {e^a}(2{\left| {\tilde M} \right|^2} + {\left| {\dot {\hat M}} \right|^2}) \;\;\;\; 
			{l_7} = 2{\left| {\tilde N} \right|^2} \;\;\;\;
			{l_8} = 4\lambda {q^2}{\left| {\tilde N} \right|^2}
		\end{equation}
		are integrable functions and 
		\begin{equation}\label{V7 c}
			{c_7} = 2\bar m_1^2 + 2\bar m_2^2+\bar m_2^2{{\bar {\hat \alpha} }_I}^2+\bar m_2^2{{\bar {\hat \beta} }_I}^2+3
		\end{equation}
		are positive constants.
		
		\subsubsection{${V_8} = \left\| {{h_2}} \right\|_b^2$}
		By using dynamic \eqref{target system2 b}, we obtain
		\begin{equation}\label{V8_1}
			\begin{footnotesize}
				\begin{aligned}
					&{{\dot V}_8} = 2{I_b}[{h_2}{\partial _t}{h_2}] = 2{I_b}[{h_2}(\mu {\partial _x}{h_2})]\\
					&- 2{I_b}[{h_2}\int_0^x {\left[ {{{\hat \alpha }_t}(x,\xi ) - \lambda {{\tilde \alpha }_\xi }(x,\xi ) + \mu {{\tilde \alpha }_x}(x,\xi ) + {{\hat m}_4}\alpha (x,\xi )} \right.} \\
					&\;\;\;\;\;\;\;\;\;\;\;\;\;\;\; - {{\hat m}_1}\tilde \alpha (x,\xi ) - {{\hat m}_3}\tilde \beta (x,\xi )]\hat \psi (\xi )d\xi ]\\
					&- 2{I_b}[{h_2}\int_0^x {[{{\hat \beta }_t}(x,\xi ) + \mu {{\tilde \beta }_\xi }(x,\xi ) + \mu {{\tilde \beta }_x}(x,\xi ) - {{\hat m}_2}\tilde \alpha (x,\xi )} \\
					&\;\;\;\;\;\;\;\;\;\;\;\;\;\;\; + {{\hat m}_4}\hat \beta (x,\xi )]\hat \varphi (\xi )d\xi ]\\
					& - 2{I_b}[{h_2}\int_0^x {\hat \alpha (x,\xi )} [{{\hat m}_1}\hat e(\xi ) + {{\hat m}_2}\hat \tau (\xi )]d\xi ]\\
					&- 2{I_b }[{h_2}\int_0^x {\hat \beta (x,\xi )} [{{\hat m}_3}\hat e(\xi ) + {{\hat m}_4}\hat \tau (\xi )]d\xi ]\\
					&- 2{I_b}[{h_2}\int_0^x {\hat \alpha (x,\xi )} [{n^T}(\xi )\dot {\hat M}]d\xi ] \\
					&- 2{I_a }[{h_2}\int_0^x {\hat \beta (x,\xi )} [{z^T}(\xi )\dot {\hat N}]]\\
					&- 2(\lambda  + \mu ){I_b}[{h_2}\tilde \alpha (x,x)\hat \psi (x)] \\
					&- 2\lambda q{I_b}[{h_2}\hat \alpha (x,0)({h_2}(0) + \hat \tau (0))]\\
					&+ 2\mu {I_b}[{h_2}\hat \beta (x,0)h_2 (0)] + 2{{\hat m}_3}{I_b}[{h_2}\hat e] + 2{{\hat m}_4}{I_b}[h_2 \hat \tau ] \\
					&+ 2{{\hat m}_4}{I_b}[{h_2}\hat \varphi (x)] + 2{I_b}[h_2 {z^T}(x )\dot {\hat N}], 
				\end{aligned}
			\end{footnotesize}
		\end{equation}
		Then, by employing inequalities \eqref{cauchy a}, as well as boundary conditions \eqref{hat parameter bounds} and \eqref{approximate kernel boundary}, we have:
		\begin{equation}\label{V8_2}
			\begin{footnotesize}
				\begin{aligned}
					&{{\dot V}_8} \le  - \mu {h_2}^2(0) - \mu b\left\| {{h_2}} \right\|_b^2 + \left\| {{h_2}} \right\|_b^2 + \bar \hat \alpha _t^2{e^b}{\left\| {{f_2}} \right\|^2} + \left\| {{h_2}} \right\|_b^2\\
					&+ {\lambda ^2}\bar {\tilde \alpha} _\xi ^2{e^b}{\left\| {{f_2}} \right\|^2} + \left\| {{h_2}} \right\|_b^2 + {\mu ^2}\bar {\tilde \alpha} _x^2{e^b}{\left\| {{f_2}} \right\|^2} + \left\| {{h_2}} \right\|_b^2\\
					&+ \bar m_4^2{{\bar \alpha }^2}{e^b}{\left\| {{f_2}} \right\|^2} + \left\| {{h_2}} \right\|_b^2 + \bar m_1^2{{\bar {\tilde \alpha} }^2}{e^b}{\left\| {{f_2}} \right\|^2} + \left\| {{h_2}} \right\|_b^2\\
					&+ \bar m_3^2{{\bar {\tilde \beta} }^2}{e^b}{\left\| {{f_2}} \right\|^2} + \left\| {{h_2}} \right\|_b^2 + \bar {\hat \beta} _t^2{e^b}{\left\| {\hat \varphi } \right\|^2} + \left\| {{h_2}} \right\|_b^2\\
					&+ {\mu ^2}\bar {\tilde \beta} _\xi ^2{e^b}{\left\| {\hat \varphi } \right\|^2} + \left\| {{h_2}} \right\|_b^2 + {\mu ^2}\bar {\tilde \beta} _x^2{e^b}{\left\| {\hat \varphi } \right\|^2} + \left\| {{h_2}} \right\|_b^2\\
					&+ \bar m_4^2{{\bar {\hat \beta} }^2}{e^b}{\left\| {\hat \varphi } \right\|^2} + \left\| {{h_2}} \right\|_b^2 + \bar m_2^2{{\bar {\tilde \alpha} }^2}{e^b}{\left\| {\hat \varphi } \right\|^2} + \left\| {{h_2}} \right\|_b^2\\
					&+ \bar m_1^2{e^b}{{\bar {\hat \alpha} }^2}{\left\| {\hat e} \right\|^2} + \left\| {{h_2}} \right\|_b^2 + \bar m_2^2{e^b}{{\bar {\hat \alpha} }^2}{\left\| {\hat \tau } \right\|^2} + \left\| {{h_2}} \right\|_b^2\\
					&+ \bar m_3^2{e^b}{{\bar {\hat \beta} }^2}{\left\| {\hat e} \right\|^2} + \left\| {{h_2}} \right\|_b^2 + \bar m_4^2{e^b}{{\bar {\hat \beta} }^2}{\left\| {\hat \tau } \right\|^2} + \left\| {{h_2}} \right\|_b^2\\
					&+ {e^b}{{\bar {\hat \alpha} }^2}{\left| {\dot {\hat M}} \right|^2}{\left\| n \right\|^2} + \left\| {{h_2}} \right\|_b^2 + {e^b}{{\bar {\hat \beta} }^2}{\left| {\dot {\hat N}} \right|^2}{\left\| z \right\|^2} + \left\| {{h_2}} \right\|_b^2\\
					&+ {\lambda ^2}{{\bar {\tilde \alpha} }^2}{e^b}{\left\| {{f_2}} \right\|^2} + \left\| {{h_2}} \right\|_b^2 + {\mu ^2}{{\bar {\tilde \alpha} }^2}{e^b}{\left\| {{f_2}} \right\|^2} + {\lambda ^2}{q^2}{{\bar {\hat \alpha} }^2}\left\| {{h_2}} \right\|_b^2\\
					&+ 2\left\| 1 \right\|_b^2({h_2}^2(0) + {{\hat \tau }^2}(0)) + {\mu ^2}{{\bar {\hat \beta} }^2}\left\| {{h_2}} \right\|_b^2 + \left\| 1 \right\|_b^2{h_2}^2(0)\\
					&+ \left\| {{h_2}} \right\|_b^2 + \bar m_3^2{e^b}{\left\| {\hat e} \right\|^2} + \left\| {{h_2}} \right\|_b^2 + \bar m_4^2{e^b}{\left\| {\hat \tau } \right\|^2}\\
					&+ \left\| {{h_2}} \right\|_b^2 + \bar m_4^2{e^b}{\left\| {\hat \varphi } \right\|^2} + \left\| {{h_2}} \right\|_b^2 + {e^b}{\left| {\dot {\hat N}} \right|^2}{\left\| z \right\|^2} \\
					&\le  - \mu {h_2}^2(0) + (25 - \mu b + {\lambda ^2}{q^2}{{\bar {\hat \alpha} }^2} + {\mu ^2}{{\bar {\hat \beta} }^2})\left\| {{h_2}} \right\|_b^2\\
					&+ (\bar {\hat \alpha} _t^2 + {\lambda ^2}\bar {tilde \alpha} _\xi ^2 + {\mu ^2}\bar {\tilde \alpha} _x^2 + \bar m_4^2{{\bar \alpha }^2} + \bar m_1^2{{\bar {\tilde \alpha} }^2}\\
					&\;\;\;\;\;\; + \bar m_3^2{{\bar {\tilde \beta} }^2} + {\lambda ^2}{{\bar {\tilde \alpha} }^2} + {\mu ^2}{{\bar {\tilde \alpha} }^2}){e^b}{\left\| {{f_2}} \right\|^2}\\
					&+ (\bar {\hat \beta} _t^2 + {\mu ^2}\bar {\tilde \beta} _\xi ^2 + {\mu ^2}\bar {\tilde \beta} _x^2 + \bar m_4^2{{\bar {\hat \beta} }^2} + \bar m_2^2{{\bar {\tilde \alpha} }^2}\\
					&\;\;\;\;\;\; + \bar m_4^2){e^b}{\left\| {\hat \varphi } \right\|^2}\\
					&+ (\bar m_1^2{{\bar {\hat \alpha} }^2} + \bar m_3^2{{\bar {\hat \beta} }^2} + \bar m_3^2){e^b}{\left\| {\hat e} \right\|^2}\\
					&+ (\bar m_2^2{{\bar {\hat \alpha} }^2} + \bar m_4^2{{\bar {\hat \beta} }^2} + \bar m_4^2){e^b}{\left\| {\hat \tau } \right\|^2}\\
					&+ {e^b}{{\bar {\hat \alpha} }^2}{\left| {\dot {\hat M}} \right|^2}{\left\| n \right\|^2} + ({{\bar {\hat \beta} }^2} + 1){e^b}{\left| {\dot {\hat N}} \right|^2}{\left\| z \right\|^2}\\
					&+ 3{e^b}{h_2}^2(0) + 2{e^b}{{\hat \tau }^2}(0)
				\end{aligned}
			\end{footnotesize}
		\end{equation}
		From equations \eqref{real status boundary} and \eqref{estimated error boundary}, using property \eqref{property 2}, we can obtain
		\begin{equation}\label{V8_3}
			\begin{footnotesize}
				\begin{aligned}
					&{{\dot V}_8} \le  - \mu {h_2}^2(0) + [25 - \mu b + {\lambda ^2}{q^2}{{\bar {\hat \alpha} }^2} + {\mu ^2}{{\bar {\hat \beta} }^2}\\
					&\;\;\;\; + 2N_2^2(\bar  {\hat \beta} _t^2 + {\mu ^2}\bar {\tilde \beta} _\xi ^2 + {\mu ^2}\bar {\tilde \beta} _x^2 + \bar m_4^2{{\bar  {\hat \beta} }^2} + \bar m_2^2{{\bar {\tilde \alpha} }^2} + \bar m_4^2)]\left\| h \right\|_b ^2\\
					&+ [\bar {\hat \alpha} _t^2 + {\lambda ^2}\bar {\tilde \alpha} _\xi ^2 + {\mu ^2}\bar {\tilde \alpha} _x^2 + \bar m_4^2{{\bar \alpha }^2} + \bar m_1^2{{\bar {\tilde \alpha} }^2} + \bar m_3^2{{\bar {\tilde \beta} }^2} + {\lambda ^2}{{\bar {\tilde \alpha} }^2} + {\mu ^2}{{\bar {\tilde \alpha} }^2}\\
					&\;\;\;\; + 2N_1^2(\bar  {\hat \beta} _t^2 + {\mu ^2}\bar {\tilde \beta} _\xi ^2 + {\mu ^2}\bar {\tilde \beta} _x^2 + \bar m_4^2{{\bar  {\hat \beta} }^2} + \bar m_2^2{{\bar {\tilde \alpha} }^2} + \bar m_4^2)]{e^b}{\left\| f \right\|^2}\\
					&+ 2(\bar m_1^2{{\bar {\hat \alpha} }^2} + \bar m_3^2{{\bar  {\hat \beta} }^2} + \bar m_3^2){e^b}{\left\| e \right\|^2} + 2(\bar m_2^2{{\bar {\hat \alpha} }^2} + \bar m_4^2{{\bar {\hat \beta} }^2} + \bar m_4^2){e^b}{\left\| \tau  \right\|^2}\\
					&+ [(2\bar m_1^2{{\bar {\hat \alpha} }^2} + 2\bar m_3^2{{\bar  {\hat \beta} }^2} + 2\bar m_3^2){\left| {\tilde M} \right|^2} + {{\bar {\hat \alpha} }^2}{\left| {\dot {\hat M}} \right|^2}]{e^b}{\left\| n \right\|^2}\\
					&+ [(2\bar m_2^2{{\bar {\hat \alpha} }^2} + 2\bar m_4^2{{\bar  {\hat \beta} }^2} + 2\bar m_4^2){\left| {\tilde N} \right|^2} + ({{\bar  {\hat \beta} }^2} + 1){\left| {\dot {\hat N}} \right|^2}]{e^b}{\left\| z \right\|^2}\\
					&+ 3{e^b}{h_2}^2(0) + 4{e^b}{\tau ^2}(0) + 4{e^b}{\left| {\tilde N} \right|^2}{\left| {z(0)} \right|^2}\\
					&\le  - (\mu b - {c_8})\left\| h \right\|_b^2{V_8} + {l_9}[{V_1} + {V_2}] + {l_{10}}[{V_4} + {V_5}]\\
					&+ {l_{11}}{V_7} + {l_{12}}{V_8} + {c_9}{e^{b + a}}{V_9} + {c_{10}}{e^b}{V_{10}}\\
					&+ {l_{13}}{{z^2}(0)} + {l_{14}}{h_2^2}(0) - \mu {h^2_2}(0) + 4{e^b}{\tau ^2}(0)
				\end{aligned}
			\end{footnotesize}
		\end{equation}
		where 
		\begin{subequations}\label{V8 l}
			\begin{small}
				\begin{align}
					&{l_9} = {e^{b + a}}[(2\bar m_1^2{{\bar {\hat \alpha} }^2} + 2\bar m_3^2{{\bar {\hat \beta} }^2} + 2\bar m_3^2){\left| {\tilde M} \right|^2} + {{\bar {\hat \alpha} }^2}{\left| {\dot {\hat M}} \right|^2}] \\
					&{l_{10}} = {e^b}[(2\bar m_2^2{{\bar {\hat \alpha} }^2} + 2\bar m_4^2{{\bar {\hat \beta} }^2} + 2\bar m_4^2){\left| {\tilde N} \right|^2} + ({{\bar {\hat \beta} }^2} + 1){\left| {\dot {\hat N}} \right|^2}] \\
					&\begin{aligned}
						{l_{11}} &=\! {e^{b + a}}[\bar {\hat \alpha} _t^2 \!+\! {\lambda ^2}\bar {\tilde \alpha} _\xi ^2 \!+\! {\mu ^2}\bar {\tilde \alpha} _x^2 \!+\! \bar m_4^2{{\bar \alpha }^2} \\
						&\;\;\;\;\;\;\;\; +\! \bar m_1^2{{\bar {\tilde \alpha} }^2} \!+\! \bar m_3^2{{\bar {\tilde \beta} }^2} \!+\!  {\lambda ^2}{{\bar {\tilde \alpha} }^2} \!+\! {\mu ^2}{{\bar {\tilde \alpha} }^2} \\
						&\;\;\;\;\;\;\;\; +\! 2N_1^2(\bar {\hat \beta} _t^2 \!+\! {\mu ^2}\bar {\tilde \beta} _\xi ^2 \!+\! {\mu ^2}\bar {\tilde \beta} _x^2 \!+\! \bar m_4^2{{\bar {\hat \beta} }^2} \!+\! \bar m_2^2{{\bar {\tilde \alpha} }^2} + \bar m_4^2)]					
					\end{aligned}\\
					&{l_{13}} = 4{\left| {\tilde N} \right|^2}{e^b} \;\;\;\;\;\;\;\; {l_{14}} = 3{e^b}
				\end{align}
			\end{small}
		\end{subequations}
		are integrable functions and 
		\begin{subequations}\label{V8 c}
			\begin{small}
				\begin{align}
					&{c_8} = 25 + {\lambda ^2}{q^2}{{\bar {\hat \alpha} }^2} + {\mu ^2}{{\bar {\hat \beta} }^2} + 2N_2^2(\bar m_4^2{{\bar {\hat \beta} }^2} + \bar m_2^2{{\bar {\tilde \alpha} }^2} + \bar m_4^2) \\
					&{c_9} = 2(\bar m_1^2{{\bar {\hat \alpha} }^2} + \bar m_3^2{{\bar {\hat \beta} }^2} + \bar m_3^2) \\
					&{c_{10}} = 2(\bar m_2^2{{\bar {\hat \alpha} }^2} + \bar m_4^2{{\bar {\hat \beta} }^2} + \bar m_4^2)					
				\end{align}
			\end{small}
		\end{subequations}
		are positive constants.
		
		\subsubsection{${V_9} = \left\| e \right\|_{ - a}^2$ and ${V_{10}} = \left\| \tau  \right\|_b^2$}
		By using \eqref{property 1}, \eqref{real error bound 1} and \eqref{real error bound 2}, we obtain
		\begin{equation}\label{V9}
			\begin{aligned}
				{{\dot V}_9} &=  - 2\lambda {I_{ - a}}[e{e_x}] \!=\!  -\! \lambda a\left\| e \right\|_{ - a}^2 \!-\! \lambda {e^{ - a}}{e^2}(1) \!+\! \lambda {e^2}(0)\\
				&\le  - \lambda a\left\| e \right\|_{ - a}^2 \!=\!  - \lambda a{V_9}
			\end{aligned}
		\end{equation}
		and
		\begin{equation}\label{V10}
			\begin{aligned}
				{{\dot V}_{10}} &= 2\mu {I_b}[\tau {\tau _x}] \!= \!- \mu b\left\| \tau  \right\|_b^2 \!+\! \mu {e^b}{\tau ^2}(1) \!-\! \mu {\tau ^2}(0)\\
				&\le  - \mu b\left\| \tau  \right\|_b^2 \!=\!  - \mu b{V_{10}}
			\end{aligned}
		\end{equation}
		
		}

		%\section*{Acknowledgment}
		
		%\section*{Acknowledgment}
		
		\ifCLASSOPTIONcaptionsoff
		\newpage
		\fi

		\bibliographystyle{IEEEtran}
		\bibliography{ref.bib}

		\vfill